\documentclass[12pt]{amsart}
\usepackage{amsmath,amsthm, amscd, amssymb, amsfonts}
\usepackage[dvips]{graphicx}
\usepackage[all]{xy}
\numberwithin{equation}{section}
\usepackage[english]{babel}
\usepackage[dvipsnames]{xcolor}
\newcommand{\bab}{\color{DarkOrchid}{}}
\newcommand{\eab}{\normalcolor{}}

\usepackage{mathtools}
\usepackage{mathtools}
\usepackage[english]{babel}
\usepackage{tikz}
\usetikzlibrary{matrix}
\newcommand{\mdn}{\medbreak\noindent}
\newcommand{\tet}{\theta}
\newcommand{\vvec}{\mathrm{Vec}}
\newcommand{\ca}{\mathcal A}
\newcommand{\innv}{\mathrm{Inv}}
\usepackage{color}

\newcommand{\End}{\mbox{\rm End\,}}
\theoremstyle{plain}

\newtheorem{theorem}{Theorem}[section]
\newtheorem{corollary}[theorem]{Corollary}
\newtheorem{proposition}[theorem]{Proposition}

\theoremstyle{definition}
\newtheorem{definition}[theorem]{Definition}

\theoremstyle{remark}
\newtheorem{remark}[theorem]{Remark}

\renewcommand{\1}{\textbf{1}}

\newcommand{\Z}{{\mathcal Z}}
\newcommand{\R}{{\mathcal R}}

\newcommand{\B}{{\mathcal B}}
\newcommand{\C}{{\mathcal C}}
\newcommand{\Ss}{{\mathcal S}}
\newcommand{\D}{{\mathcal D}}

\newcommand\id{\operatorname{id}}
\newcommand\Aut{\operatorname{Aut}}

\newcommand\Tr{\operatorname{Tr}}
\newcommand{\cE}{\mtc{E}}
\newcommand{\cS}{\mtc{S}}

\newcommand\GL{\operatorname{GL}}
\newcommand\cop{\operatorname{cop}}

\newcommand\op{\operatorname{op}}
\newcommand\Hom{\operatorname{Hom}}
\newcommand\Rep{\operatorname{Rep}}

\newcommand\Ind{\operatorname{Ind}}

\newcommand\FPdim{\operatorname{FPdim}}

\newcommand\vect{\operatorname{Vec}}

\newcommand\Irr{\operatorname{Irr}}
\newcommand\rep{\operatorname{Rep}}

\newcommand{\ot}{\otimes}
\newcommand{\mtc}{\mathcal}
\newcommand{\lam}{\lambda}
\newcommand{\lb}{\label}

\newcommand{\al}{\alpha}

\newcommand{\eps}{\epsilon}
\newcommand{\bn}{\begin}

\newcommand{\krn}{\mtr{ker}}

\newcommand{\ts}{\times}

\newcommand{\onh}{On the other hand}

\numberwithin{equation}{section}

\newtheorem{defn}[theorem]{Definition}
\newtheorem{cor}[theorem]{Corollary}
\newtheorem{rem}[theorem]{Remark}
\newtheorem{example}[theorem]{Example}

\newcommand{\fp}{\mtr{FPdim}}
\newcommand{\bp }{\bn{proposition}}
\newcommand{\ep}{\end{proposition}}
\newcommand{\gr}{K_0}

\newcommand{\dw}{\downarrow}
\newcommand{\uw}{\uparrow}

\newcommand{\ch}{\chi}
\newcommand{\mtr}{\mathrm}

\newcommand{\ncm}{\newcommand}\newcommand{\gm}{\gamma}
\numberwithin{equation}{section}
\newcommand{\el}{\end{lemma}}\newcommand{\bl}{\bn{lemma}}
\newcommand{\et}{\end{thm}}\newcommand{\bt}{\bn{thm}}
\newcommand{\ovr}{\overline}
\newcommand{\beqarn}{\begin{eqnarray*}}
\newcommand{\eeqarn}{\end{eqnarray*}}
\newcommand{\beqn}{\bn{equation*}}
\newcommand{\eeqn}{\end{equation*}}
\newcommand{\bpf}{\bn{proof}}
\newcommand{\epf}{\end{proof}}
\ncm{\cX}{\mtc{X}}
\ncm{\wt}{\widetilde}
\newcommand{\ra}{\rightarrow}
\ncm{\ro}{\rho}
\ncm{\sg}{\sigma}
\ncm{\np}{\newpage}
\ncm{\ebl}{\end{thebibliography}}
\ncm{\bbl}{\begin{thebibliography}}
\ncm{\chd}{_{ _{\ch}}}
\ncm{\ald}{_{ _{\al}}}
\ncm{\cP}{\mathcal{P}}
\ncm{\ei}{e_i}
\ncm{\eij}{e_{i,\;j}}
\ncm{\bdef}{\begin{defn}}
\ncm{\edf}{\end{defn}}
\newcommand{\wtan}{With the above notations\; }
\ncm{\bc}{\begin{cor}}
\ncm{\ec}{\end{cor}}
\ncm{\er}{\end{rem}}
\ncm{\br}{\begin{rem}}
\ncm{\bd}{\begin{document}}
\ncm{\ed}{\end{document}}
\ncm{\beq}{\begin{equation}}
\ncm{\eeq}{\end{equation}}
\ncm{\bea}{\begin{eqnarray}}
\ncm{\eea}{\end{eqnarray}}
\ncm{\beanon}{\begin{eqnarray*}}
\ncm{\eeanon}{\end{eqnarray*}}\ncm{\ek}{\eps|_K}\ncm{\diez}{\#}
\ncm{\bwt}{\bowtie}
\ncm{\cC}{\mtc{C}}
\ncm{\cs}{\cS}
\ncm{\cA}{\mathcal{A}}
\ncm{\HKer}{\mtr{HKer}}
\ncm{\LKER}{\mtr{LKER}}
\ncm{\aad}{\mtr{ad}}
\ncm{\Dr}{\mtr{D}}
\ncm{\cD}{\mathcal{D}}
\ncm{\Dc}{\mtc{D}}
\ncm{\Vc}{\mtr{Vec}}
\ncm{\cK}{\mtc{K}}
\ncm{\cM}{\mtc{M}}
\ncm{\ce}{\mtc{E}}
\ncm{\chr}{character }
\ncm{\bw}{\bwt}
\ncm{\cY}{\mtc{Y}}
\ncm{\cd}{\cD}
\ncm{\bx}{\boxtimes}
\ncm{\cc}{\cC}\ncm{\irr}{\Irr}
\ncm{\bgm}{\Gamma}
\ncm{\goesto}{\mapsto}
\ncm{\one}{\mathrm{1}_{\cc}}
\ncm{\bvec}{\mtr{Vec}}
\ncm{\omg}{\omega}
\ncm{\og}{\omega}
\ncm{\dt}{\delta}
\ncm{\cg}{\mtc{G}}
\ncm{\blue}{\textcolor[rgb]{.00, .00, 1.00}}
\ncm{\red}{\textcolor[rgb]{1.00, .00, .00}}
\ncm{\green}{\textcolor[rgb]{0.00, .60, .40}}
\ncm{\dzs}{\#_{\sg}}
\ncm{\md}{\medskip}\ncm{\opl}{\oplus}\ncm{\xra}{\xrightarrow}
\ncm{\Gm}{\Gamma}\ncm{\ind}{\mtr{Ind}}\ncm{\res}{\mtr{Res}}\ncm{\core}{\mtr{core}}\ncm{\lker}{\mtr{LKer}}\ncm{\bb}{\blue}
\renewcommand{\subjclassname}{%}
\textup{2000} Mathematics Subject Classification}
\ncm{\bq}{\beq}\ncm{\eq}{\eeq}
\ncm{\bne}{\begin{enumerate}}
\ncm{\ene}{\end{enumerate}}
\ncm{\bqarn}{\beqarn}
\newcommand{\cz}{\mathcal{Z}}
\newcommand{\ct}{\mathrm{t}}%\rho\mathrm{J}
\newcommand{\bll}{\blue}
\newcommand{\bnex}{\begin{example}}
\newcommand{\enex}{\end{example}}
\ncm{\eqarn}{\eeqarn}
\begin{document}
\ncm{\inv}{^{-1}}
\title[Fusion subcategories of an equivariantization]{On M\"{u}ger's  centralizer for a certain class of braided fusion categories}
\date{\today}
%\newgeometry{top=1.9cm,bottom=1.3cm}
\begin{abstract} 
We  provide a general formula for M\"{u}ger's centralizer of any fusion subcategory of a braided fusion category containing a tannakian subcategory. This entails a description for M\"{u}ger's centralizer of all fusion subcategories of a group theoretical braided fusion category.  
\end{abstract}
\author{S. Burciu}
\address{Inst.\ of Math.\ ``Simion Stoilow" of the Romanian Academy, Research Unit 5,
P.O. Box 1-764, RO-014700, Bucharest, Romania}\email{sebastian.burciu@imar.ro}
\thanks{ The research was done in part during a stay at the Erwin Schr\"odinger Institute, Vienna, in the frame of the Programme ``Modern Trends in Topological Quantum Field Theory'' in February 2014. The author  thanks the ESI and the organizers of the Programme for the support and kind hospitality. This work was partially supported by a grant of the Romanian National Authority for Scientific Research, CNCS-UEFSCDI, project number PN II - ID - PCE-2011-3-0039}
\maketitle
\section{Introduction and main results}\lb{intro}
%\subsection{Conventions and notation} 
It is known that braided fusion categories containing a non-trivial Tannakian subcategory are equivariantizations of $G$-cross braided fusion categories \cite{DGNOI, gnn}.
In this  paper  we present a formula for the M\"uger centralizer in equivariantizations of braided $G$-crossed fusion categories.

Let $\C$ be a fusion category over $k$. The isomorphism class of an object $X$ of $\C$ will be denoted by $[X]$. However, when the context presents no ambiguity, we shall indicate an object and its class by the same letter.
We shall use the notation $\Irr(\C)$ to indicate the set of isomorphism classes of simple objects of $\C$. %\bll{Thus $\Irr (\C)$ is a basis of the Grothendieck ring of $\C$.}
%\subsection{Group actions on fusion categories}

Let  $\ct: \underline G \to \underline \Aut_{\otimes} \C$ be an action of $G$ by tensor autoequivalences on a fusion category $\C$ (see Section \ref{s1} for the precise definition). Such an action induces naturally an action of $G$ on the set $\Irr (\C)$ of isomorphism classes of simple objects of $\C$.

A parameterization for all the fusion subcategories of an equivariantization will be given in \cite{buna2}. For the sake of completeness of the paper we recall this description in Section \ref{fusion-subcategories}. Below we recall the  notion of a fusion datum which classifies all the fusion subcategories, see \cite{buna2}.
%\bll{The definition depends on a global chosen system of isomorphisms (GCSI). I fix a set of isomporphisms for each orbit class.}
%\bab Explain that $\tau_{Y, Z}^{U}$ gives a fixed one dimensional representation of $H$. \eab
\subsection{Equivariantization under group actions}\label{s1}
%\subsection{Equivariantization under a finite tensor group action}
%^{}Along this section $k$ will be an algebraically closed field of characteristic zero.
%\subsubsection{Definition of a tensor $G$-action}
\mdn Let $\C$ be a fusion
category, $G$ a finite group and $\ct: \underline
G \to \underline \Aut_{\ot} \, \C$ be an action of $G$ on $\C$ by
$k$-linear tensor autoequivalences. Thus, for every $g \in G$, we have a
$k$-linear tensor functor $\ct^g: \C \to \C$ and natural isomorphisms of tensors functors
$$\ct^{g,h}_2 : \ct^g \ct^h \to \ct^{gh}, \quad g, h \in G,$$ and
$\ct_0 : \id_\C \to \ct^1$, subject to the following conditions:
\begin{align}\label{ro-2} & (\ct^{ab, c}_2)_X \, (\ct^{a,b}_2)_{\ct^c(X)} =
(\ct^{a, bc}_2)_X \, \ct^a((\ct^{b,c}_2)_X), \\ & (\ct^{a,
1}_2)_X \ct^{a}_{\ct_0(X)} = (\ct^{1, a}_2)_X
(\ct_0)_{\ct^{a}(X)},
\end{align}
for all objects $X \in \C$, and for all $a, b, c \in G$.
%See \cite[Subsection 4.1]{DGNOI}.
By the naturallity of $\ct^{g,h}_2$, $g, h \in G$, we have the following relation:
\begin{equation}\label{nat-ro}\ct^{gh}(f) \, (\ct_2^{g, h})_{Y} = (\ct_2^{g, h})_X \, \ct^g\ct^h(f),
\end{equation} for every morphism $f: Y \to X$ in $\C$. 
We shall assume in what follows that $\ct^1 = \id_\C$ and
$\ct_0$, $\ct^{g, 1}_2$, $\ct^{1, g}_2$ are identities.
%\subsubsection{Definition of the equivariantized fusion category}
\mdn
 Let $\C^G$ denote the corresponding
\emph{equivariantization}. Recall that $\C^G$ is a finite
semisimple $k$-linear category whose objects are $G$-equivariant
objects of $\C$, that is, pairs $(X, \mu)$, where $X$ is an object
of $\C$ and $\mu = (\mu^g)_{g \in G}$, such that $\mu^g:\ct^gX
\to X$ is an isomorphism, for all $g \in G$, satisfying
\begin{equation}\label{deltau} \mu^g \ct^g(\mu^h) = \mu^{gh} (\ct^{g,
h}_2)_X, \quad \forall g, h \in G, \quad \quad
\mu_e{\ct_0}_X=\id_X.\end{equation} A morphism $f: (X, \mu) \to
(X', \mu')$ in $\C^G$ is a morphism $f: X \to X'$ in $\C$ such
that $f\mu^g = {\mu'}^g$, for all $g \in G$.
\mdn We shall also say that an object $X$ of $\C$ is
\emph{$G$-equivariant} if there exists such a collection $\mu =
(\mu^g)_{g \in G}$ so that $(X, \mu) \in \C^G$. Note that $\mu$ is
not necessarily unique.
%\subsubsection{Tensor product in $\cc^{G}$}
\mdn
For all $g \in G$ the autoequivalence $\ct^g$ is endowed with a monoidal structure $(\ct_2^g)_{X, Y}: \ct^g(X \otimes Y) \to \ct^g(X) \otimes \ct^g(Y)$, $X, Y \in \C$ such the following relation holds:
\begin{equation} {(\ct_2^{gh})}_{X, Y} {(\ct_2^{g, h})}_{X \otimes Y} = (({\ct_2^{g, h}}^{-1})_{X} \otimes ({\ct_2^{g, h}}^{-1})_{Y}) \, {(\ct_2^{g})}_{\ct^h(X), \ct^h(Y)} \, {\ct^g((\ct_2^{h})}_{X, Y}),
\end{equation}
for all $g, h \in G$, $X, Y \in \C$.
\mdn Then $\C^G$ is also a fusion category with tensor product
$(X, \mu_X) \otimes (Y, \mu_Y) = (X \otimes Y, (\mu_X \otimes \mu_Y){\ct_2}_{X, Y})$, where for all $g \in G$, ${(\ct_2)}_{X, Y}^g: \ct^g(X \otimes Y) \to \ct^g(X) \otimes \ct^g(Y)$ is given by the monoidal structure on $\ct^g$.
%\subsubsection{Forgetful functor}
\mdn The forgetful functor $F: \C^G \to \C$, $F(X, \mu) = X$, is a dominant functor. The functor $F$ has a left adjoint $L: \C
\to \C^G$, defined by $L(X) = (\oplus_{s \in G}\ct^s(X), \mu_X)$,
where $(\mu_X)_g: \oplus_{t \in G}\ct^g\ct^s(X) \to \oplus_{s \in
G}\ct^s(X)$ is given componentwise by the isomorphisms
$(\ct^{g,s}_2)_X$. %The composition $T_\ct = FL: \C \to \C$ is a faithful $k$-linear monad on $\C$ such that $\C^G$ is equivalent to the category $\C^{T_\ct}$ of $T_\ct$-modules in $\C$. See\cite[Subsection 5.3]{tensor-exact}.
%\subsubsection{Inertia subgroups of objects} 
\subsection{On the chosen set of isomorphisms $\mathrm C$}  Let $Y$ be a simple object of $\C$ and let $G_{Y} \subseteq G$ denote the inertia subgroup of $Y$, that is,
$
G_{Y}=\{g \in G\;|\; \ct^g(Y)\cong Y\}.
$  By the definition of $G_Y$, there exist isomorphisms $c_{Y}^g:\ct^g(Y) \to
Y$, for all $g \in G_Y$. For all $g, h \in G_Y$, the composition
$c_{Y}^g \ct^g(c_{Y}^h)({\ct^{g, h}_{2_Y}})^{-1}(c_{Y}^{gh})^{-1}$  defines
an isomorphism $Y \to Y$. Since $Y$ is a simple object, by Schur's Lemma, there
exists a nonzero scalar $\alpha_Y(g, h) \in k^{*}$ such that
\begin{equation}\label{alfa}  \alpha_Y(g, h)^{-1} \id_Y = c_Y^g \ct^g(c_Y^h)({\ct^{g, h}_{2}})_{Y}^{-1}(c_Y^{gh})^{-1}: Y \to Y. \end{equation}
This defines a map $\alpha_Y: G_Y \times G_Y \to k^*$ which
is a $2$-cocycle on $G_Y$, see \cite{buna}. As explained in \cite{buna} the cohomology class of $\al_{Y}$ does not depend on the chosen set of isomorphisms $c_{Y}^{g}$ with $g \in G_{Y}$. We denote by $\widetilde\alpha_Y$ its cohomologus class. \begin{remark}\label{dependence-alfa} The cocycle $\alpha_Y$ measures the possible obstruction for $(Y, C_{Y})$ to be a $G_Y$-equivariant object, where $u = (c_Y^g)_{g \in G_Y}$.

\mdn
Consider another choice of isomorphisms ${c'}^{g}_{Y}: \ct^g(Y) \to Y$, $g \in G$. Then since $Y$ is a simple object, the composition $c_{Y}^g({c'_{Y}}^g)^{-1}: Y \to Y$ is given by
scalar multiplication by some $f_Y(g) \in k^*$, for all $g \in G$.  Denoting by $\al'_Y$ the 2-cocycle resulting from the choice of $({c'}^g)_{g\in G_Y}$, it easy to see that $\alpha_Y$ and $\al'_Y$ differ by the coboundary of the cochain $f_Y: G_Y \to k^*$. This implies that the cohomology class $[\alpha_Y]$ of $\alpha_Y$ depends only on the isomorphism class of the simple object $Y$. More precisely one has the following relation (see \cite{buna}):
\beq\label{cocycled}
\al'_Y(g,h)=\frac{f_Y(gh)}{f_Y(g)f_Y(h)}\al_Y(g,h)
\eeq
\end{remark} \mdn Fix a set $\{Y_{0}, Y_{1}, \cdots ,Y_{s}\} $ of representative objects for all the isomorphism classes of simple objects of $\cc$. Denote $G_{i}:=G_{Y_{i}}$ and let $c_{Y_{i}}^{g}:\ct^{g}(Y_{i})\ra Y_{i}$ be arbitrary isomorphisms with $c_{Y_{i}}^{1}=\id_{Y_{i}}$ for any $g \in G_{i}$. As above we define a $2$-cocycle $\al_{Y_{i}}:G_{i}\times G_{i}\ra k^{*}$ by the equality
\beq
\al_{Y_{i}}^{-1}(g, g')=c_{Y_{i}}^g \ct^g(c_{Y_{i}}^h)({\ct^{g, h}_{2_Y}})^{-1}(c_{Y_{i}}^{gh})^{-1}.
\eeq
Note that $\al_{Y_{i}}(g,1)=\al_{Y_{i}}(1,g)=1$ for any $g \in G$, i.e $\al_{Y_{i}}$ is a normalized cocycle. Moreover by Remark \ref{unity} one can assume that the values of $\al_{Y_{i}}(g, g')$ are all roots of unity.\mdn 
For any simple object of $Y$  with $Y\cong Y_{i}$ and $g \in G_{Y}$ we denote by $c_{Y}^{g}:T^{g}(Y)\ra Y$ the isomorphism obtained by the composition $f^{-1}c_{Y_{i}}^{g}T^{g}(u)$ where $u:Y\ra Y_{i}$ is an arbitrary chosen isomorphism. 
Note that $c_{Y}^{g}$  does not depend on $u$ and $c_{Y}^{1}=\id_{Y}$ for any simple object $Y$. Moreover $\al_{Y}(g,g')=\al_{Y_{i}}(g,g')$ for all $g,g'\in G_{Y}$.  {\it We denote this common value by $\al_{[Y]}(g,g')$.}
Moreover, since $u$ is unique up to a scalar, the following digram commutes:
{\tiny
\begin{center}
\begin{tikzpicture}
  \matrix (m) [matrix of math nodes, row sep=3em,column sep=4em,minimum width=3em]
  { T^{g}(Y) & Y \\
   T^{g}(Y_{i})  &   Y_{i}\\
    };%5
  \path[-stealth]
  (m-1-1) edge node [above] {$c_{Y}^{g}$} (m-1-2)
    (m-2-1) edge node [above] {$c_{Y_{i}}^{g}$} (m-2-2)
         
               (m-1-1) edge node [left] {$T^{g}(v)$} (m-2-1)
    (m-1-2) edge node [right] {$v$} (m-2-2)
   ;
\end{tikzpicture}
   \end{center}
   }
  \mdn for any $v:Y\ra Y_{i}$ isomorphism.
\mdn For any two simple isomorphic objects $Y \cong Y'$ and any $f:Y\ra Y'$ an isomorphism one has from above that
\beq\label{isocho}
v \circ c^{g}_{Y}=c^{g}_{Y'}\circ \ct^{g}(v).
\eeq
Indeed, since $v$ is unique up to scalar, one can choose $v$ to be composition of two isomorphisms $Y \ra Y_{i} \ra Y'$ where $Y_{i}$ is the corresponding representative for $Y$ and $Y'$. Then the commutativity follows from the definition of $c_{Y}^{g}$ and commutativity of the above diagram.
\mdn Let $g \in G$ and $Y$ be a simple object of $\cc$. Since $\ct^{g}(Y)$ is also a simple object there is  a  scalar$D_{g,Y}(h)\in k^{*}$  such that
\beq\label{D}
c_{\ct^{g}(Y)}^{ghg^{-1}}=D^{-1}_{g, Y}(h)[\ct^{g}(c_{Y}^{h})(\ct_{2}^{g,h})_{Y}^{-1}(\ct_{2}^{ghg^{-1}, g})_{ _{Y}}]
\eeq
for any $h \in G_{Y}$.
From the commutative diagram below it is easy to see that this scalar does not depend on the class of isomorphism of $Y$.
{\tiny
\begin{center}
\begin{tikzpicture}
  \matrix (m) [matrix of math nodes, row sep=6em,column sep=4em,minimum width=4em]
  { \ct^{ghg^{-1}}(\ct^{g}(Y)) & \ct^{gh}(Y) & \ct^{g}(\ct^{h}(Y)) & \ct^{g}(Y)& \ct^{ghg^{-1}}(\ct^{g}(Y))\\
   \ct^{ghg^{-1}}(\ct^{g}(Y_{i}))  & \ct^{gh}(Y_{i}) & \ct^{g}(\ct^{h}(Y_{i})) & \ct^{g}(Y_{i})& \ct^{ghg^{-1}}(\ct^{g}(Y_{i}))\\
    };%5
  \path[-stealth]
  (m-1-1) edge node [above] {$(\ct_{2}^{ghg^{-1}, g})_{ _{Y}}$} (m-1-2)
   (m-1-2) edge node [above] {$(\ct_{2}^{g,h})_{Y}^{-1}$} (m-1-3)
   (m-1-3) edge node [above] {$\ct^{g}(c_{Y}^{h})$} (m-1-4)
   (m-1-4) edge node [above] {$[c_{\ct^{g}(Y)}^{ghg^{-1}}]^{-1}$} (m-1-5)
    (m-2-1) edge node [above] {$(\ct_{2}^{ghg^{-1}, g})_{ _{Y_{i}}}$} (m-2-2)
         (m-2-2) edge node [above] {$(\ct_{2}^{g,h})_{Y_{i}}^{-1}$} (m-2-3)
   (m-2-3) edge node [above] {$\ct^{g}(c_{Y_{i}}^{h})$} (m-2-4)
   (m-2-4) edge node [above] {$[c_{\ct^{g}(Y_{i})}^{ghg^{-1}}]^{-1}$} (m-2-5)
               (m-1-1) edge node [right] {$\ct^{ghg^{-1}}(\ct^{g}(f))$} (m-2-1)
    (m-1-2) edge node [right] {$\ct^{gh}(f)$} (m-2-2)
    (m-1-3) edge node [right] {$\ct^{g}(\ct^{h}(f))$} (m-2-3)
    (m-1-4) edge node [right] {$\ct^{g}(f)$} (m-2-4)
    (m-1-5) edge node [right] {$\ct^{ghg^{-1}}(f)$} (m-2-5)
   ;
\end{tikzpicture}
   \end{center}
}\noindent Note that the last rectangle is commutative by Equation \eqref{isocho}. The other rectangles are commutative by the naturality of the group action. \mdn We denote this common value by $D_{g, [Y]}(h)$.% \bll{A similar diagram to the one above implies that for $g=g' (\text{mod} \;G_{Y})$ one has that $D_{g, [Y]}(h)=D_{g', [Y]}(h)$ for all $h \in G_{Y}$.}
 \mdn Denote by $C_{Y}:=\{c_{Y}^{g}\}_{g \in G_{Y}}$ the set of chosen isomorphisms for any simple object $Y$ of $\cc$. 
One can  choose $C_{Y}$ such that $c_{Y}^{1}=\id_{Y}$. In this case one has that $\al_{Y}(1,h)=\al_{Y}(h,1)=1$. {\it For the rest of the paper we fix a set $\mathrm C:=\cup_{Y \in \mathcal \Irr(\cc)} C_{Y}$ of chosen isomorphisms as above. Moreover note that we can choose $c_{\1}^{g}=\id_{\1}$ for any $g \in G$ since $\ct^{1}=\id_{\cc}$ and $\ct_{0}$ is assumed to be identity.}
\subsection{Hom spaces as projective representations}
Let $(X, \mu) \in \C^G$ and let $Y \in \Irr(\C)$. Then following \cite{buna} the space $\Hom_\C(Y, X)$ becomes a projective representation of $G_Y$ with factor set $\alpha_Y$. Its structure is defined by
\begin{equation}\label{pr}\pi (g) (f) = \mu^g \ct^g(f) (c_Y^g)^{-1} : Y \to X,\end{equation}
for all $f \in \Hom_\C(Y, N)$.
%\subsection{On the projective representation $\tau_{U}^{Y, Z}(t,s)$} \label{newproj}
\mdn
Let $U, Y, Z \in \Irr(\cc)$.
Recall by \cite[Proposition 3.6]{buna} that the vector space $\Hom_{\cc}(U,\ct^m(Y)\ot \ct^n(Z))$ is a projective representation of $T:=G_U\cap G_{\ct^m(Y)} \cap G_{\ct^n(Z)}$ with associated $2$-cocycle $(\al_{[U]}|T)(\;^m\al_{[Y]}|_T)^{-1}(\;^n\al_{[Z]}|_T)^{-1}$. The action of $g \in T$ is given by
\beq\lb{pr'}
g.f=(\ct^{m}(c_Y^{m^{-1}gm})\ot \ct^n(c_Z^{n^{-1}gn}))((\ct^2)^{m, m^{-1}gm}_Y\ot (\ct^2)^{n, n^{-1}gn}_Z)(\ct^g_2)^{\ct^m(Y),\;\ct^n(Z)}\ct^g(f)(c^g_U)^{-1}
\eeq
for all $f \in \Hom_{\cc}(U,\ct^m(Y)\ot \ct^n(Z))$.
\mdn 
Note that the definition of the above projective representation depend on all chosen set of isomorphisms $c_Y, c_Z, c_U$. \mdn As in \cite{buna} denote this projective representation by  $\tau^{Y,Z}_U(m, n)$.  For shortness we write $\tau_U^{Y, Z}:=\tau_U^{Y, Z}(1,1)$. \mdn Note that $\tau_{U}^{Y, Z}(m, n)$ depend on the chosen set of isomorphisms for $U, \ct^{m}(Y), \ct^{n}(Z)$. Clearly changing these system $\mathrm C$ of isomorphisms one obtains an equivalent projective representation. \mdn 
However, note that for the above system $C$ of chosen isomorphisms one has that $\tau_{U}^{Y,Z}(m,n)\simeq \tau_{Y_{i}}^{Y_{j}, Y_{k}}(m,n)$ as $T$-projective representations, where $U\cong Y_{i}$, $Y\cong Y_{j}$ and $Z \cong Y_{k}$. Indeed, choosing isomorphisms $u:U \ra Y_{i}$, $v:Y \ra Y_{j}$ and $w:Z\ra Y_{k}$ it is easy to verify (by the naturality of the $G$-action) that the $k$-linear  map
$
\phi:\Hom_{\cc}(U, \ct^{m}(Y)\otimes \ct^{n}(Z))\ra \Hom_{\cc}(Y_{i}, \ct^{m}(Y_{j})\ot \ct^{n}( Y_{k}))
$
given by $f \mapsto u^{-1}f(\ct^{m}(v) \ot \ct^{n}(w))$ 
is an isomorphism of projective representations. This is equivalent to the commutativity of the diagram below. In this diagram the first and the last rectangles are commutative by Equation \eqref{isocho}. The other rectangles are commutative by the naturality of the action of $G$.
{\Tiny
\begin{center}
\begin{tikzpicture}
  \matrix (m) [matrix of math nodes, row sep=8em,column sep=4em,minimum width=1 em]
  { U &  Y_{i}
%line 2
\\ \ct^{g}(U)  & \ct^{g}(Y_{i})
 %line 3
   \\ \ct^{g}(\ct^{m}(Y) \ot \ct^{n}(Z)) & \ct^{g}(\ct^{m}(Y_{j}) \ot \ct^{n}(Y_{k})) 
      % line 4
   \\ \ct^{g}(\ct^{m}(Y))\ot \ct^{g}(\ct^{n}(Z)) &  \ct^{g}(\ct^{m}(Y_{j}))\ot \ct^{g}(\ct^{n}(Y_{k}))
     % line 5
   \\\ct^{gm}(Y)\ot \ct^{gn}(Z) & \ct^{gm}(Y_{j}) \ot \ct^{gn}(Y_{_{k}}))
   %line 6
   \\ \ct^{m}(\ct^{m\inv gm}(Y))\ot \ct^{n}(\ct^{n\inv gn}(Z)) &  \ct^{m}(\ct^{m\inv gm}(Y_{j}))\ot \ct^{n}(\ct^{n\inv gn}(Y_{k}))
   %line 7
   \\ \ct^{m}(Y)\ot \ct^{n}(Z) & \ct^{m}(Y_{j}) \ot \ct^{n}(Y_{_{k}})
   \\}; \path[-stealth]
  (m-1-1) edge node [above] {$u$} (m-1-2)
    (m-2-1) edge [dashed] node   [above] {$\ct^{g}(u)$} (m-2-2)
       (m-3-1) edge [dashed] node   [above] {$\ct^{g}(\ct^{m}(v)\ot \ct^{n}(w)$} (m-3-2)
     (m-4-1) edge [dashed] node   [above] {$\ct^{g}(\ct^{m}(v))\ot \ct^{g}(\ct^{n}(w))$} (m-4-2)
     (m-5-1) edge [dashed]  node    [below] {$ \ct^{gm}(v) \ot \ct^{gn}(w)$} (m-5-2)
     (m-6-1) edge [dashed] node    [below] {$ \ct^{m}(\ct^{m\inv gm}(v)) \ot \ct^{n}(\ct^{n\inv gn}(w))$} (m-6-2)
      (m-7-1) edge  node    [below] {$ \ct^{m}(v) \ot \ct^{n}(w)$} (m-7-2)
         %columns
               (m-1-1) edge node [left] {$(c^{g}_{U})^{-1}$} (m-2-1)
   (m-1-2) edge node [right] {$(c^{g}_{Y_{i}})^{-1}$} (m-2-2)     
   (m-2-1) edge node [left] {$\ct^{g}(f)$} (m-3-1)
    (m-2-2) edge node [right] {$\ct^{g}(\phi(f))$} (m-3-2)
    
     (m-3-1) edge node [left] {$(\ct^{g})_{2}^{\ct^{m}(Y), \ct^{n}(Z)}$} (m-4-1)
    (m-3-2) edge node [right] {$(\ct^{g})_{2}^{\ct^{m}(Y_{j}), \ct^{n}(Y_{k})}$} (m-4-2)
   
     (m-4-1) edge node [left] {$(\ct_{2})^{g, m}_{Y}\ot (\ct_{2})^{g, n}_{Z}$} (m-5-1)
    (m-4-2) edge node [right] {$(\ct_{2}^{g, m})_{Y_{j}}\ot (\ct_{2}^{g, n})_{Y_{k}}$} (m-5-2)
   (m-5-1) edge node [left] {$(\ct_{2}^{m, m\inv gm})^{-1}_{Y}\ot (\ct_{2}^{n, n\inv gn})^{-1}_{Z}$} (m-6-1)
    (m-5-2) edge node [right] {$(\ct_{2}^{m, m\inv gm})^{-1}_{Y_{j}}\ot (\ct_{2}^{n, n\inv gn})^{-1}_{Y_{k}}$} (m-6-2)
     (m-6-1) edge node [left] {$\ct^{m}(c_{Y}^{m\inv gm})\ot \ct^{n}(c_{Z}^{n\inv gn})$} (m-7-1)
    (m-6-2) edge node [right] {$\ct^{m}(c_{Y_{j}}^{m\inv gm})\ot \ct^{n}(c_{Y_{k}}^{n\inv gn})$} (m-7-2)
   ;
\end{tikzpicture}
  \end{center}
} \noindent {\it For this reason, by fixing the system $C$ as above for the rest of the paper, we use the notation $\ch_{\tau^{[Y],[Z]}_{[U]}(m,n)}$  for the character of these isomorphic representations. By abuse of notation sometimes we also write $\ch_{\tau^{Y, Z}_{U}(m,n)}$ instead of $\ch_{\tau^{[Y],[Z]}_{[U]}(m,n)}$.}
\subsection{Fusion datum} In this subsection we recall the notion of a fusion datum relative to a given set of chosen isomorphisms $C$.
\bn{definition}\label{fusdatum}
A $\mathrm C$-\emph{fusion datum}
with respect to the tensor auto equivalence action $\ct: \underline G \to \underline \Aut_{\otimes} \C$ is a triple $(\Ss, H, \lambda)$ where:
\begin{enumerate}\item $\Ss$ is a $G$-stable fusion subcategory of $\C$.
\item $H$ is a normal subgroup of $G$ acting trivially on $\Irr (\Ss)$.
\item $\al_{Y}|_{H}$ is cohomologus trivial for any simple object $Y$ of $\cs$. Or equivalently the object $Y$ is $H$-equivariant.
\item Every element $h \in H$ acts as a scalar on the projective representation $\tau_{U}^{Y, Z}(t,s)$, for any simple objects $Y,Z, U \in \Irr(\cs)$ and any $t,s\in G$ such that $U$ is a constituent of $\ct^{t}(Y)\ot 
\ct^{s}(Z)$.
\item $\lambda: \Irr (\Ss) \times H \to k^*$ is a map satisfying the following conditions:
\begin{enumerate}\item[(i)] The following holds, for all $[Y] \in \Irr(\Ss)$, $h, h' \in H$:
\beq\lb{hh'}\lambda([Y], hh') = \lambda([Y], h) \, \lambda([Y], h')\al_{[Y]}(h,h').\eeq
\item[(ii)]  For all $h \in H$ the map $\lambda(- , \;h): \Irr(\Ss) \to U(1)$ is a {\it twisted grading} on $\Ss$ with respect to $\tau$, that is,
\beq\lb{yz}
\lam([U],h)=\lam([Y],h)\, \lam([Z],h)(\frac{\ch_{\tau_{[U]}^{[Y],[Z]}}(h)}{\ch_{\tau_{[U]}^{[Y],[Z]}}(1)})
\eeq
for any irreducible constituent $U$ of $Y \ot Z$.% Recall \cite{buna} that $\tau_U^{Y,Z}:=\tau_U^{Y,Z}(1,1)$.
%Moreover, for any constituent $U$ of $Y \otimes Z$, the character of $k_{\al_{[U]} |_H}[H]$ given by $\lam(-,\;Y)
%\dw^{G_Y}_H\lam(-,\;Z)\dw^{G_Z}_H\tau_U^{Y, Z}\dw^T_H$ is a $G_U$-stable character of $H$ with respect to $\al_{[U]}$.% (with multiplicity $\tau^{Y,Z}_U(1)$.)
%\blue{$\tau_{U}^{Y, Z}(t,s)|_H=\tau_{U}^{Y, Z}(t,s)(1)$}
%by a cyclic group $C(h)$, for all $h \in H$,
\item[(iii)] $\lam$ is a $G$-invariant function, i.e:
\beq\label{inv}\lambda([\ct^g(Y)],h)=  \lambda([Y], g\inv hg)D_{g, [Y]}(g^{-1}hg),
\eeq for all $g \in G$, $h \in H$, $Y \in \Irr(\Ss)$.
\end{enumerate}
\end{enumerate}
\end{definition}

%{\bf Do I have to define it as with values in U(1) or it follows from the definion ?}
Note that Condition (2) from above can be expressed as $H \subseteq \bigcap_{Y \in \irr(\cs)}G_Y$, where $G_Y$ is the inertia group of $Y$.
%\bll{Made it as a remark that one can choose it to be mutiplicative; or that will be proven in the other version.}
\mdn Then a  parameterization of fusion subcategories of an equivariantization can be given in terms of the \emph{fusion data} introduced above:
\begin{theorem}\label{cls-by-fusiondata}(\cite{
buna2}) Let $\C$ be a fusion category and let $\ct: \underline
G \to \underline \Aut_{\otimes} \C$ be an action of a finite group $G$ on $\C$ by tensor autoequivalences. Fix a system $\mathrm C$ of chosen isomorphism as above. Then there is a bijective correspondence between the following two sets:
\begin{enumerate}\item[(a)] Fusion subcategories $\D \subseteq \C$.
\item[(b)] $\mathrm C$-fusion data $(\Ss, H,  \lambda)$.
\end{enumerate}
\end{theorem}
We denote by $\cc(\cs, H, \lam)$ the category associated to the $C$-fusion datum $(\cs, H, \lam)$.
%Note that under this correspondence, the Frobenius-Perron dimension of the fusion subcategory $\D( \Ss, H, \lambda)$ associated to the fusion datum $(\Ss, H,  \lambda)$ is  $[G: H] \FPdim \Ss$.%}\mdn \bll{this is a different lemma}
\mdn
For the sake of completeness we will recall the proof of Theorem \ref{cls-by-fusiondata}  in Section \ref{fusion-subcategories}. In order to establish Theorem \ref{cls-by-fusiondata}, we use the explicit description of all simple objects of the equivariantization $\cc^G$ from \cite{buna}. 
The isomorphism classes of simple objects of $\cc^G$ are parameterized by pairs $([Y], \pi)$, where $[Y]$ runs over the orbits of the action of $G$ on the isomorphism classes of simple objects of $\cc$, and $\pi$ is an irreducible $\al_{Y}$-projective representation of the inertia subgroup $G_{[Y]}$ (\cite{buna}). The irreducible projective representation $\pi$ has an associated two cocycle $\al_{Y}\in Z^2(G_{[Y]}, k^*)$ (see Equation \ref{alfa}). The corresponding simple object associated to the pair $([Y],\pi)$ is denoted by $S_{{Y},\pi}$.
%\md
%As a particular case, in Section \ref{example} for  $\cc=\bvec^{\omg\inv}_{G}$ with $\omg \in Z^3(G, k^*)$ we recover the description of all fusion subcategories of the twisted Drinfeld double $\rep(D^{\omega} (G))$ given in \cite{nnw}.
%\subsection{The main results of the paper} 
\subsection{$G$-cross braided fusion categories}
It is known that braided fusion categories containing a non-trivial Tannakian subcategory are equivariantizations of $G$-cross braided fusion categories \cite{gnn, DGNOI}.
\mdn
In the second part of the paper  we consider equivariantizations of braided $G$-crossed fusion categories. For  a braided $G$-crossed category we define a notion of a {\it $G$-centralizer} (see subsection \ref{gcross}). We say that two simple objects $X, Y$ of a braided $G$-crossed category $\cc$  $G$-centralize each other if there are two simple objects of the equivariantization $\cc^{G}$ which centralize each other and have $X$ respectively $Y$ as constituents. For two fusion categories $\cd, \ce$ of $\cc$ we write $\cd \perp_{G} \ce$ if any two objects of $\cd$ and $\ce$ $G$-centralize each other
\mdn  Let $(\cs, H, \lam)$ be a $\mathrm C$-fusion datum and suppose that the support of $\cs$ is a subgroup $K \leq G$. In this part of the paper we adopt the notation $\cd(H,K, \cs, \lam)$ for the fusion subcategory $\cd(\cs,H, \lam)$ associated to the given fusion datum. Recall that if $\cd \subset \cc$ is a fusion subcategory of a graded fusion category $\cc=\oplus_{g \in G}\cc_{g}$ then its support is defined by 
\beq
\mathrm{supp}(\cd)=\{g \in G\;|\; \cd\cap \cc_{g}\neq 0\}.
\eeq 
\mdn In the next theorem we give a description  for the fusion datum associated to the M\"{u}ger centralizer of a given fusion subcategory.
%\mdn Put the theorem for two objects to centralize each other.
%\mdn\bll{Explain how the groups $H$ and $K$ are associated.}\mdn
\begin{theorem}\label{form} Let $\cc$ be a braided $G$-crossed   fusion category and $(\cs, H, \lam)$ a $\mathrm C$-fusion datum with $K:=\mathrm{supp}(\cs)$. Then the following hold:
\bne
\item $[K, H]=1$, where $[x,y]=xyx^{-1}y^{-1}$ is the usual commutator.
\item
One has that:
 \beqn
 \cc(H, K, \cs, \lam)'=\cc(K, H', \mathcal{T}, \tilde{\lam})
 \eeqn
where $\mathcal{T}$ is a fusion subcategory of $\cc$ whose  support is $H'\subseteq H$. Moreover $\mathcal{T}$ is stable under the action of $G$ and it satisfies $\cs \perp_{G}\mathcal{T}$. 
\item With the above notations if $a \in K$ and $b \in H'$ then one has that $\tilde{\lam}:\irr(\mathcal{T})\times K\ra k^{*}$ is given by 
 \beq\label{tilde}
 \tilde{\lam}([Y_{b}], a)=\lam^{-1}([X_{a}], b)\omega
 ([X_{a}], [Y_{b}])
 \eeq
for any object $X_{a}\in \cc_{a}\cap \cs$ and any $Y_{b }\in \cc_{b}\cap \mathcal{T}$.
\item Moreover if $\cc^{G}$ is nondegenerate then $H'=K$.
\ene
\end{theorem}
%\bll{The conjecture is that $\mathcal T=\cs^{\perp_{G}}\cap \cc(H)$ and this is a fusion subcategory of $\cc$.}
%\bll{drop $K$ or explain the notation.}\mdn
%\bll{They get a bicharacter that is twisted in both sides but the formula for $\tilde \lam$ does not depend on other things, it is just the inverse.}
\mdn
The roots of unity $\omega([X_{a}], [Y_{b}])$ are defined in subsection \ref{fmc}. Note that in the nondegenerate case one obtains that
 \beqn
 \cc(H, K, \cs, \lam)'=\cc(K, H, \mathcal{T}, \tilde{\lam})
 \eeqn
 formula that generalizes  \cite[Proposition 3.6]{nnw} and \cite[Proposition 5.2]{nnw} to any nondegenerate equivariantization of a braided $G$-crossed fusion category.
\mdn Recall that \cite[Theorem 5.3 ]{Na} states that every braided group-theoretical category is equivalent to $C(\xi)^G$, for some pointed group crossed category $C(\xi)$ associated to a normalized quasi-abelian 3-cocycle $\psi$ on a finite crossed module $(G, X, \delta)$. Thus it follows from Theorem \ref{form} that we can give a description for the M\"{u}ger centralizer of any fusion subcategory of all braided group theoretical fusion categories. %\bll{This also follows from...}
\subsection{Organization of the paper}The paper is organized as follows. In the first section we recall the description of the simple objects $S_{Y,\pi} \in \cc^{G}$ from \cite{buna}. The fusion rules of $\cc^{G}$ from loc. cit. is also recalled in this section. Next section recalls a proof of  Theorem \ref{cls-by-fusiondata}.  In Proposition \ref{eqf} we show that 
any fusion subcategory of $\cc^{G}$ is also an equivariantization of a fusion subcategory of $\cc$ by a quotient subgroup of $G$. In this section we also give a description for the lattice of fusion subcategories of an equivariantization $\cc^G$. In Section \ref{emc} we recall basic notions on braided $G$-cross fusion categories. In Section \ref{mc} we study equivariantizations of   braided $G$-crossed   fusion categories and describe M\"{u}ger's centralizer for a fusion subcategory. Section \ref{gcp} contains some examples and application of our results.
%\mdn\bll{Section \ref{example} studies centers of Tambara-Yamagami categories as introduced in \cite{ty}. Using results from \cite{gnn} we give a parameterization of all fusion subcategories of these centers and also describe their M\"{u}ger centralizer. }
\section{Simple objects of an equivariantization}\noindent
In this section we recall the description of the simple objects of an equivariantization $\cc^G$ as well as the fusion rules for the fusion category $\cc^G$ which were described in \cite{buna}.

\subsection{The relative adjoint functor} For any simple object $Y$ consider the forgetful functor $F_Y: \C^G \to \C^{G_Y}$.  Let $\R$ be a set of representatives of the left cosets
of $G_Y$ in $G$. Thus $G = \cup_{t \in
\R}tG_Y$.
\mdn
For all $(N, \nu) \in \C^{G_Y}$ let $L_Y^\R(N, \nu) = (\oplus_{t \in \R}\ct^t(N), \mu) \in \cc^{G}$, where, for
all $g \in G$, the equivariant structure $\mu_{L_Y^\R(N, \nu) }^{g}: \oplus_{t \in \R}\ct^g\ct^t(N) \to
\oplus_{t \in \R}\ct^t(N)$ is defined componentwise by the
formula \begin{equation}\label{mu}\mu_{L_Y^\R(N, \nu)}^{g, t} = \ct^s(\nu^h)
(\ct_2^{s, h})^{-1} \ct_2^{g, t}: \ct^g \ct^t(N) \to
\ct^s(N),\end{equation} and the elements $h \in G_Y$, $s \in \R$
are uniquely determined by the relation
\beq\label{indexch}
gt = sh.
\eeq
\mdn
It was noticed in \cite{buna} that if $(N, \nu) \in \C^{G_Y}$ then $L_Y(N, \nu) \in
	\C^G$. Therefore there is a well defined functor $ L_Y^\R: \C^{G_Y} \to \C^G$ which by results of {\it loc.cit} is left adjoint of the forgetful functor $F_Y$. By the uniqueness of the adjoint  functor, it follows that, up to isomorphism,
$L_Y^\R$ does not depend on the particular set choice of the set of representatives $\R$.  For this reason we will use the notation $L_{Y}$ instead of $L_Y^\R$.

\subsubsection{Notations} Let $\pi$ be a $\al$-projective representation of a finite group $G$.
We usually denote by $V_{\pi}$ the vector space associated to the projective representation $\pi$ of $G$. Therefore one can regard $\pi$ as a map $\pi: G\ra \End_k(V_{\pi})$ with $\pi(gg')=\al(g,\;g')\pi(g)\pi(g')$ for all $g,g' \in G$. We also denote by $\ch_{\pi}$ the character of the projective representation $\pi$, thus $\ch_{\pi}(h)=\text{tr}_{V_{\pi}}(L_{h})$. 
\subsubsection{Simple objects of $\cc^{G_{Y}}$}
Recall \cite{buna} that to any $\al_{[Y]}$-projective representation $\pi$ of $G_Y$ one can associate an object $\pi\ot  Y \in \cc^{G_Y}$. As an object of $\cc$ we have that $\pi\ot Y$ coincides to $V_{\pi}\ot Y$ and its equivariant structure $\nu^g$ is defined as follows:
\bq\lb{cy}
\nu^g : \ct^g(V_{\pi} \otimes Y) = V_{\pi} \otimes \ct^g( Y) \xra{\pi(g) \otimes c^g_{Y}}
V_{\pi} \otimes Y.
\eq
for all $g \in G_Y$. Conversely, if $L = (N, \nu)$ is any object of
$\C^{G_Y}$ with $N \simeq \Hom_\C(Y, N) \otimes Y$, then $V =
\Hom_\C(Y, N)$ carries a projective representation $\pi$ of $G_Y$
with the factor set $\widetilde\alpha_Y$. This projective representation is defined as in \eqref{pr}.
\mdn
In other words, there is a bijective correspondence
between isomorphism classes of simple objects $L = (N, \nu)$ of
$\C^{G_Y}$ such that $N \simeq \Hom_\C(Y, N) \otimes Y$ and
irreducible projective representations of the group $G_Y$ with
factor set $\widetilde\alpha_Y$. If the simple object $L$  corresponds to
the projective representation $\pi$,  we have $\FPdim L = \dim \pi
\, \FPdim Y$.\mdn
Note that the simple object $Y\ot V_{\pi}$ depends on the chosen set of isomorphisms $C_{Y}$ at $Y$.
\subsubsection{Parameterization of simple objects}
In  \cite [Theorem 2.12]{buna} it was shown that there is a bijective correspondence between the set isomorphism classes of simple objects $(X, \mu)$ of $\C^G$ and the set of isomorphism classes of pairs $(Y, \pi)$, where $Y$ runs over the orbits of the action of $G$ on $\Irr(\C)$ and $\pi$ is an irreducible $\al_{Y}$-projective representation of the inertia subgroup $G_Y \subseteq G$. \mdn
The simple object corresponding to the pair $(Y,
\pi)$ was denoted by $S^{C_{Y}}_{Y, \pi}$ in \cite{buna}. Then $S^{C_{Y}}_{Y, \pi}:=L_{Y}(Y \ot \pi)$ where $Y \otimes V_{\pi}\in \cc^{G_{Y}}$ was previously defined. Thus we have $S_{Y, \pi} \simeq \oplus_{t \in G/G_Y}\ct^t(V_\pi \otimes Y)$ and its equivariant structure follows from Equation \eqref{mu}.   In particular, one has that $\FPdim S_{Y, \pi} = \dim \pi [G:G_Y]
\FPdim Y$.

\begin{remark}\label{s-object} Note that $S^{C_{Y}}_{Y, \pi}$ as defined above  depends on the set of  chosen isomorphisms $C_{Y}$ defined previously.
Let us describe more explicitly the dependence of the simple object $S^{ C_{Y}}_{Y, \pi}$ on the choice of the isomorphisms $c^g_Y:\ct^g(Y)\ra Y$. 
Suppose we choose another set of isomorphisms $c'_Y= \{{c'}^g_Y\}$. Then, $Y$ being a simple object, implies that for any $g \in G_Y$ one can write that ${c'}^g_Y = f_{Y}^{C, C'}(g) c^g_Y$, for some scalar $f_{Y}^{c, c'}(g) \in k^*$.
It follows easily that $\pi \ot Y = (f_{Y}^{C, C'})^{-1}\pi \ot Y$ as objects in $\cc^{G_Y}$. Hence
\bq\lb{dyt}
S^{ C_{Y}}_{Y,\pi} = L_Y(\pi \ot Y) = L_Y((f_{Y}^{C, C'})^{-1}\pi \ot Y) = S^{C'_{Y}}_{Y, (f_{Y}^{C, C'})^{-1}\pi }. \eq
\end{remark}
\mdn
\br\label{indepofs}
One also can check that if $Y\simeq Y'$ and $c_{Y'}$ is obtained from $c_{Y}$ via an isomorphism $f:Y \ra Y'$ (i.e. satisfying Equation \eqref{isocho}) then $Y\otimes V_{\pi}\cong Y'\otimes V_{\pi}$ as objects of $\cc^{G_{Y}}$. Therefore in this case $S^{c_{Y}}_{Y,\pi}\simeq S^{c_{Y'}}_{Y',\pi}$.
\er
\mdn
{\it For the rest of this paper we write for shortness  $S_{Y,\;\pi}$ for $S_{Y, \pi}^{C_{Y}}$ where $C_{Y}\subset \mathrm C$ is the fixed set of chosen isomorphisms defined in the previous section.}
\subsection{On the choice of the isomorphisms on the orbit}
Let $G$ be a finite group and $\pi$ be a $\al$-projective representation of a subgroup $H\leq G$. For any $g \in G$ define $\;^{g}\pi$ as a projective representation of $\;^{g}H:=gHg^{-1}$ by $\;^{g}\pi(ghg\inv):=\pi(h)$ for $h \in H$. Note that $\;^{g}\pi$ is a $\;^{g}\al$-projective representation of $\;^{g}H$ where $\;^{g}\al(ghg^{-1}, gh'g^{-1})=\al(h,h')$.  \mdn Similarly, for any function $f : G_{Y}\ra k$ define by $\;^{g}f:G_{\ct^{g(Y)}}\ra k$ defined by $\;^{g}f(ghg\inv):=f(h)$ for all $h \in G_{Y}$. \mdn
With the notations from the previous section let $g \in G_{Y}$. Since $\ct^g(Y)$ is a constituent of $F_{Y}(S_{Y, \pi})$, it follows from \cite[Theorem 2.2]{buna} that $S_{Y, \pi} = S_{\ct^g(Y), \delta}$, for some irreducible projective representation $\delta$ of $G_{\ct^g(Y)}$. In this section we discuss the dependence of  $\delta$ upon $\pi$.
%The following Lemma appears in \cite{buna}.
%\bl\cite [ Lemma 2.18]{buna}\label{xh-equivariant} Let $H$ be a subgroup of $G$ and suppose that $(M, \nu) \in \C^H$. Then, for all $x \in G$, $(\ct^x(M), ^x\nu) \in \C^{^xH}$ with equivariant structure $(^x\nu)^{xhx^{-1}}: \ct^{xhx^{-1}}\ct^x(M) \to \ct^x(M)$ defined, for every $h \in H$, as the composition
%\beq
%\ct^{xhx^{-1}}(\ct^x(M))\xra{\ct_2^{xhx^{-1}, x}} \ct^{xh}(M)\xra{(\ct_2^{x,h})^{-1}}\ct^x(\ct^h(M))\xra{\ct^x(\nu^h)} %\ct^x(M).\eeq\el
%\bll{Do I need to go to $\mathcal R$?}
%Given a collection of isomorphisms $c:=\{c^g\}_{g\in G_Y}$, let as in Remark \ref{s-object} $S^{\R, c}_{Y, \pi}$ denote the associated simple object of $\cc^G$.
Recall that by Equation \eqref{D} one has \beq\label{D2}
\ct^{g}(c_{Y}^{h})(\ct_{2}^{g,h})_{Y}^{-1}(\ct_{2}^{ghg^{-1}, g})_{ _{Y}}=D_{g, [Y]}(h)c_{\ct^{g}(Y)}^{ghg^{-1}}
\eeq
for any $h \in G_{Y}$.
\bl \lb{rot} With the above notations one has
\beq
S_{Y, \pi}\cong S_{\ct^g(Y),\; \;^{g}D_{g, [Y]}( -)\;^g\pi},
\eeq for all $g \in G$.
\el
\bpf 
Recall that $S_{Y,\pi} = L_Y(\pi \ot Y)$. As in the proof of \cite[Theorem 2.2]{buna} it is enough to check that $V_{\;^{g}D_{g, [Y]}( -)\;^g\pi} \ot \ct^g(Y)$ is a constituent of the restriction $F_{\ct^g(Y)}(L_Y(\pi \ot Y))$ where $F_{\ct^{g}(Y)}$ is the forgetful functor $F_{\ct^{g}(Y)}:\cc^{G}\ra \cc^{G_{\ct^{g}(Y)}}$. As objects of $\cc$,
$F_{\ct^g(Y)}(L_Y(\pi \ot Y))=\oplus_{s \in G/G_Y}\ct^s(V_{\pi} \ot Y)$. Restricting the equivariant structure $\mu_{S_{Y,\pi}}$ of $S_{Y,\pi}$ to an element $ghg^{-1} \in G_{\ct^g(Y)}$ it follows that $(ghg^{-1})g = gh$, and therefore one can choose $s = g$, with notations as in Equation \eqref{indexch}. Thus on the component $\ct^g(Y)\ot V_{\pi}$ the equivariant structure $\mu_{S_{Y,\pi}}^{ghg^{-1}}$ is given by the composition
\begin{align*}
\ct^{ghg^{-1}}(V_{\pi} \ot \ct^g(Y))& \xra{\id_{V_{\pi}} \ot (\ct^{tht^{-1}, t}_2)_Y} V_{\pi} \ot \ct^{gh}(Y) \xra{} \\ &\xra{\id_{V_{\pi}} \ot (\ct^{t,h}_2)_Y^{-1}} V_{\pi} \ot \ct^g(\ct^h(Y))  \xra{\pi(h) \ot \ct^g(c_Y^h)} V_{\pi} \ot \ct^g(Y).
\end{align*}
This can also be written as
\begin{equation*}
V_{\pi} \ot \ct^{ghg^{-1}}(\ct^g(Y))\xra{\pi(h) \ot {(\ct^g(c_Y^h)(\ct^{g,h}_2)_Y(\ct^{ghg\inv , g}_2)_Y)}}V_{\pi} \ot \ct^g(Y).
\end{equation*}
Formula \eqref{D2} shows that the above morphism coincides with the equivariant structure on $\;^gd_Y(g,-)\;^g\pi \ot \ct^g(Y)$ given by Equation \eqref{cy}.
\epf
\br 
With the above notations one can check by straightforward computation that  \beqn\ct^g(Y\ot V_{\pi})=\ct^g(Y)\;\ot \;\;^gD_{g, [Y]}(-)\;^g\pi \eeqn as equivariant objects of $\cc^{G_{\ct^g(Y)}}$. This  also implies that
\beqarn
S_{Y, \pi}=L_{Y}(\pi \otimes Y)&\simeq& L_{\ct^{g}(Y)}(\ct^{g}(\pi \otimes Y)\simeq \\&\simeq& L_{\ct^{g}(Y)}(\ct^g(Y)\;\ot \;\;^gD_{g, [Y]}(-)\;^g\pi= S_{\ct^{g}(Y), \;^gD_{g, [Y]}(-)\;^g\pi}
\eeqarn
\mdn
It can also shown by straightforward computation that
\bq\lb{dty}
\al_{[Y]}(h,h')=\frac{D_{g, [Y]}(hh')}{D_{g, [Y]}( h) D_{g, [Y]}( h')} \; \al_{[\ct^g(Y)]}(ghg^{-1}, gh'g^{-1}),
\eq
for all $h,h' \in G_Y$, $g \in G$.
\er
\subsection{More on the projective representations $\tau$}
Recall from the previous section that we denoted by $\tau_{U}^{Y,Z}:=\tau_{U}^{Y, Z}(1,1)$.
\bp\label{compst}
With the above notations it follows that
\beqn
\tau_{[U]}^{[\ct^m(Y)], [\ct^n(Z)]}=\tau_{[U]}^{[Y], [Z]}(m, n)\;\;^{m}(D_{m, [Y]})^{-1}\;\;^{n}(D_{n, [Z]})^{-1}
\eeqn
i.e for any $g \in G_{U} \cap mG_{Y}m^{-1}\cap nG_{Z}n^{-1}$ one has:
\beq\label{tauconj}
\tau_{[U]}^{[\ct^m(Y)], [\ct^n(Z)]}(g)=\tau_{[U]}^{[Y], [Z]}(m, n)(g)D^{-1}_{m, [Y]}(m^{-1}gm)D^{-1}_{n, [Z]}(n^{-1}gn)
\eeq
\ep 
\bpf Note that both projective representations are defined on the same vector space $\Hom_{\cc}(U ,\ct^{m}(Y)\ot \ct^{n}(Z))$. Moreover, for any $g \in G_{U} \cap mG_{Y}m^{-1}\cap nG_{Z}n^{-1}$ and $f \in \Hom_{\cc}(U, \ct^{m}(Y)\ot \ct^{n}(Z))$ using Equation 
\eqref{D2} one has the following 
{\small
%(\ct^{t}(c_Y^{t^{-1}gt})\ot \ct^s(c_Z^{s^{-1}gs}))\times \\&&\times ((\ct_2^{t, t^{-1}gt})_Y\ot (\ct_{2}^{s, s^{-1}gs})_Z)
\beqarn
[\tau_{[U]}^{[Y], [Z]}(m,n)](g)f&&= (c^{g}_{\ct^{m}(Y)}\ot c^{g}_{\ct^{n}(Z)})(\ct^g)_2^{\ct^m(Y),\;\ct^n(Z)}\ct^g(f)(c^g_U)^{-1}.\eeqarn
Similarly one has that
\beqarn
[\tau_{[U]}^{[\ct^{m}(Y)], [\ct^{n}(Z)]}](g)f&&= [\ct^{m}(c_{Y}^{m\inv g m}))\ot \ct^{n}(c_{Z}^{n\inv g n})][(\ct_{2}^{m,m^{-1}gm})_{Y}^{-1}\ot (\ct_{2}^{n,n^{-1}gn})_{Z}^{-1}]\\ && [(\ct_{2}^{g, m})_{ _{Y}}\ot (\ct_{2}^{g, n})_{ _{Z}}][(\ct^g)_2^{\ct^m(Y),\;\ct^n(Z)}\ct^g(f)(c^g_U)^{-1}]\\&=& ([\ct^{m}(c_{Y}^{m^{-1}gm})(\ct_{2}^{m,m^{-1}gm})_{Y}^{-1}(\ct_{2}^{g, m})_{ _{Y}}]\ot [\ct^{n}(c_{Z}^{n^{-1}gn})(\ct_{2}^{n,n^{-1}gn})_{Z}^{-1}(\ct_{2}^{g, n})_{ _{Z}}]) .\eeqarn
On the other hand, by Equation \eqref{D} note that 
\beqn
c_{\ct^{m}(Y)}^{g}=c_{\ct^{m}(Y)}^{m(m^{-1}gm)m^{-1}}=D^{-1}_{m, Y}(m^{-1}gm)[\ct^{m}(c_{Y}^{m^{-1}gm})(\ct^{2}_{m,m^{-1}gm})_{Y}^{-1}(\ct^{2}_{g, m})_{ _{Y}}].
\eeqn
Similarly one has that
\beqn
c_{\ct^{g}(Y)}^{n}=c_{\ct^{n}(Y)}^{n(n^{-1}gn)n^{-1}}=D^{-1}_{n, Y}(n^{-1}gn)[\ct^{n}(c_{Y}^{n^{-1}gn})(\ct^{2}_{n,n^{-1}gn})_{Y}^{-1}(\ct^{2}_{g, n})_{ _{Y}}]
\eeqn
Thus
\beqarn  [\tau_{[U]}^{[Y], [Z]}(m,n)](g)f &&= D_{m,[Y]}^{-1}({m^{\inv}}gm)D_{n, [Z]}^{-1}( {n^{\inv}}gn)\\&&([\ct^{m}(c_{Y}^{m^{-1}gm})(\ct_{2}^{m,m^{-1}gm})_{Y}^{-1}(\ct_{2}^{g, m})_{ _{Y}}]\ot [\ct^{n}(c_{Z}^{n^{-1}gn})(\ct_{2}^{n,n^{-1}gn})_{Z}^{-1}(\ct_{2}^{g, n})_{ _{Z}}])\\&&  (\ct^g_2)^{\ct^m(Y),\;\ct^n(Z)}\ct^g(f)(c^g_U)^{-1}\\ &=& D_{m, [Y]}^{-1}({m^{\inv}}gm)D_{n, [Z]}({n^{\inv}}gn)\tau_{[U]}^{[\ct^{m}(Y)], [\ct^{n}(Z)]}(g)f.\eeqarn
}
\epf
%\mdn \bll{The above computations can easily be seen in the diagram below.}
\subsection{On the tensor product of two simple objects}
Let $S_{Y, \pi}$ and $S_{Z,\gm}$ be two simple objects of $\cc^G$. The following Corollary follows from the proof of \cite [Theorem 3.9]{buna}, see also \cite[Corollary 3.6]{buna}.
\begin{corollary}\label{constituent} A simple object $S_{U,\delta}$ is a
constituent of a tensor product of simple objects $S_{Y, \pi}
\otimes S_{Z,\gamma}$ if and only if there exist elements $t, s\in G$ such that
\begin{enumerate}\item[(a)] $\Hom_\C(U, \rho^t(Y) \otimes \rho^s(Z)) \neq 0$ and
\item[(b)] $m_{T}(\delta\vert_T,\;^t\pi\vert_T \otimes ^s\gamma\vert_T \otimes \tau^{[Y], [Z]}_{[U]}(t, s))\neq 0$, where $T = G_U \cap G_{\ct^{t}(Y)} \cap G_{\ct^{s}(Z)}$.
\end{enumerate} 
\end{corollary} \section{Fusion subcategories of an equivariantization}\label{fusion-subcategories}
\mdn In this subsection we shortly recall from \cite{buna2} the description of fusion subcategories of an equivariantization. Recall that for the rest of the paper we have fixed a system $C=\{c_{Y}\}_{Y \in \irr(\cc)}$ of chosen isomorphisms for any $Y \in \irr(\cc)$ satisfying Equation \eqref{isocho} and with $c^{g}_{\1}=\id_{\1}$ for all $g \in G$. Moreover we have $c_{Y}^{1}=\id_{Y}$ for any simple object $Y$ and the values of $\al_{Y}(g,g')$ are  roots of unity for any $g,g' \in G$. By the object $S_{Y, \pi}$ we will understand from now on the simple object $S^{c_{Y}}_{[Y], \pi}$ associated to the chosen set $C$. 
\mdn Let $\pi: G \to \GL(V)$ be a finite dimensional representation of  $G$ on the vector space $V$.
Then the (trivial) object $V \otimes \1 \in \C$ has a $G$-equivariant structure defined by $\pi(g) \otimes \id_\1: \ct^g(V \otimes \1) \to V \otimes \1$. To see that this defines an equivariant structure one needs to use the fact that $c^{g}_{\1}=\id_{\1}$. 
This induces an embedding of fusion categories $\rep G \to \C^G$ that gives rise to an exact sequence of fusion categories
\begin{equation}\label{sec} \rep G \to \C^G \to \C.
\end{equation}
See \cite[Subsection 5.4]{tensor-exact}. Note that if $\pi$ is an irreducible representation of $G = G_\1$ on $V$, then the simple object $(V \otimes \1, (\pi(g) \otimes \id_\1)_g)$ of $\C^G$ is isomorphic to
the simple object $S_{\1, \pi}$ corresponding to the pair $(\1, \pi)$.
\subsection{The fusion datum associated to a fusion subcategory} Let $\D$ be a fusion subcategory of $\cc^G$. One defines the following datum associated to $\cd$. Let $\cs_{\cd}$ be the fusion subcategory of $\cc$ generated by $F(\cd)$ where $F$ is the forgetful functor $F:\cc^G \ra \cc$. Thus $\Irr(  \cS_{\D})$ consists of the set of all $Y\in \Irr(\cc)$ such that there is some $(X,\mu) \in \cc^G$ with $Y$ a subobject of $F(X,\; \mu)$.
 Clearly $\cs_{\cd}$ is a fusion subcategory of $\cc$ stable under the action of $G$. Indeed if $Y$ is a subobject of $F(S)$ with $S$ a simple object of $\cc^G$ then $\ct^g(Y)$ is a subobject of $F(\ct^g(S))\cong F(S)$.
 \mdn
Define also a normal subgroup $H_{\cd}$ of $G$ by  the equality $\cd \cap \rep(G)=\rep(G/H_{\cd}).$ Note that the normal subgroup $H_{\D}$ satisfies also the following equation:
\begin{equation}\label{fh}
 H_{\D}=\bigcap_{\{\gamma\;|\;S_{\1, \gamma}\in \D\}}\krn_G{\gamma}.
\end{equation}

\subsection{The dual of a simple object}\label{involution} Let $Y \in \Irr (\C)$. Then the multiplicity of the unit object of
$\cc$ in the tensor product $Y \otimes Y^*$ is one. Hence
$\tau_Y : = \tau_{\1}^{Y, Y^*}= \Hom_{\cc}(\1, Y\ot
Y^*)$ is a one dimensional (linear) representation of $G =
G_{\1}$. In particular, it follows that the cohomology class of the product $\al_{Y} \al_{Y^*}$ is
trivial on the subgroup $G_Y = G_{Y^*}$.  Recall that the dual $\pi^*$ of the $G_Y$-projective
representation  $\pi$ is defined as $V_{\pi}^*$ with $\pi^*(h)(f)
= f \circ \pi(h)^{-1}$. This is an $\al_{Y}^{-1}$-projective
representation of $G_Y$.  Moreover in these conditions \cite[Proposition 3.12]{buna} implies that 
\beq\label{dualf}
S_{Y, \pi}^{*}=S_{Y^{*}, \pi^{*}\tau_{Y}^{-1}}.
\eeq
\bl\label{1const}
One has that
$S_{\1, \delta}$ is a constituent of $S_{Y,\pi}\otimes (S_{Y, \gm})^{*}$ if and only if $\delta$ is a constituent of $(\pi\gm^{*})\uw^{G}_{G_{Y}}$.\el
\bpf
Applying Corollary \ref{constituent} for $r=s=1$ it follows $S_{1, \delta}$ is a constituent of $S_{Y,\pi}\otimes (S_{Y, \gm})^{*}$ if and only if $\delta$ is a constituent of $\pi\gm^{*}\tau_{Y}^{-1}\tau_{Y}=\pi\gm^{*}$. The statement of the lemma follows now by Frobenius reciprocity.
\epf
 %\bll{ Recall also that for a given chosen system $c_{Y}$ of isomorphisms one can define $c^{*}_{Y^{*}}$ as a chosen set of isomorphisms for $Y^{*}$ by $Y^{*}\xra{(c_{Y}^{g})^{*}} T^{g}(Y^{*})\simeq T^{g}(Y^{*})$ for any $g \in G_{Y}$.}
%\bp\cite [ Proposition 3.4]{buna}\label{dual} The dual object of $S_{Y,\pi}\in \cc^G$ is determined by $$S_{Y,\pi}^* \simeq S_{Y^*,\;\pi^*}.$$
%$$S_{Y,\pi}^* \simeq S_{Y^*,\;\pi^*\ot\tau_Y^{-1}},$$ where $\tau_Y = \Hom_{\cc}(\1, \;Y\ot Y^*)$.\ep
%\begin{corollary}\cite [ Corollary 3.13]{buna}\label{onecst} Let $Y, Z \in \Irr(\C)/G$ and let $\pi$, $\gm$ be projective representations of $G_Y$ and $G_Z$ with factor sets determined by \eqref{alfa}.
%$\alpha_Y$, $\alpha_Z$, respectively.Let also $\delta$ be an irreducible representation of $G$. Then $S_{\1,\delta}$ is a constituent of $S_{Y, \pi}\ot S^*_{Z, \gm}$ if and only if $Z = Y^*$ and $\delta$ is a constituent of $(\pi \ot\gm^*)\uw^G_{G_Y}$. \qed
%\end{corollary}
\bn{lemma}\label{zentre}
Suppose that $\D$ is a fusion subcategory of
$\cc^G$. Then $H_{\D}\subset \mtr{core}_G(G_Y)$ for all $Y \in
\Irr(\cs_{\cd})$.  If $S_{Y,\;\gamma} \in \D$ and $S_{Y,\; \gamma'}\in
\D$ then for all $h \in H$ one has
\begin{eqnarray*}
\frac{\ch_{\gm}(h)}{\ch_{\gm}(1)} &=& \frac{\ch_{\gm'}(h)}{\ch_{\gm'}(1)}
=\lam_{\cd}(Y,h)
\end{eqnarray*}
for some  root of unity $\lam_{\cd}(Y,h)  \in k^{*}$.
\end{lemma}
\bpf 
If $S_{Y ,\;\gamma}\in \D$ and $S_{Y,\; \gamma'}\in \D$ then
$S_{Y ,\gamma}\ot (S_{Y ,\gamma'})^*\in \D$. But  $(S_{Y ,\gamma'})^*=S^{}_{Y^* ,\gamma'^*\tau_Y^{-1}}$ by Equation \eqref{dualf}. Then from Lemma \ref{1const} it follows that for all $\delta $, a constituent of $(\gm\gm'^*\tau_{Y}\tau_{Y}^{-1})\uw_{G_Y}^G=(\gm\gm'^*)\uw_{G_Y}^G$,  one has that $S_{1,\; \delta}\in \D$. Thus $H_{\cd}$ acts as a scalar on $\delta$ and therefore on $(\gm\gm'^*)\uw_{G_Y}^G$. By Equation \eqref{fh} we have
$H_{\D}\subset\krn_G((\gm\gm'^*)\uw_{G_Y}^G)$. Then Lemma \ref{krnind}
implies the first statement.

\mdn The second statement follows as well
from this lemma since for all $h \in H_{\cd}$ one has $|\ch_{\gamma}( h)|\leq \ch_{\gamma}( 1)$ and $\ch_{\gamma'^*}(h)=\overline{\ch_{\gamma'}(h)}$. Thus if $(\gm\gm'^*)(h)=(\gm\gm'^*)(1)$ then by Lemma \ref{a1}  there is a root of unity $\lam_{\cd}(Y, h)$ such that \beqn
\frac{\ch_{\gamma}( h)}{\ch_{\gamma}( 1)} = \frac{\ch_{\gm'}(h)}{\ch_{\gm'}(1)}=\lam_{\cd}(Y, h).
\eeqn
\epf
\noindent For any fusion subcategory $\cd \subset \cc^{G}$ define a function $\lambda_{\cd} :\irr(\cs_{\cd})\times H \ra k^*$ as following:
\bq
\lam_{\cd}([Y], h)={\frac{\ch_{\pi}(h)}{\ch_{\pi}(1)}}
\eq
if $S_{Y, \pi}\in \cd$. Note that $\lam_{\cd}$ is well defined since if $Y\simeq Y'$ then $S_{Y',\pi}\simeq S_{Y, \pi}\in \cd$ by Remark \ref{indepofs}.
\bl \label{ews}  Let $Y \in \cs_{\cd}$ be a simple object. With the above notations note that $S_{Y, \delta} \in \cd$ if and only if $ \frac{\ch_{\delta}(h)}{\ch_{\delta}(1)}=\lam_{\cd}([Y], h)$ for all $h \in H_{\cd}$.\el
\bpf If $S_{Y, \delta} \in \cd$ then clearly $ \frac{\ch_{\delta}(h)}{\ch_{\delta}(1)}=\lam_{\cd}([Y], h)$ by defintion. \mdn  Conversely, since $ Y \in \cs_{\cd}$ then there is a $\al_{Y}$-projective representation $\pi$ of $G_{Y}$ such that $S_{Y,\pi}\in \cd$. Thus $\frac{\ch_{\pi}(h)}{\ch_{\pi}(1)}=\lam_{\cd}([Y], h)$ for all $h \in H_{\cd}$. \mdn Note that as above one has that $S_{1, (\gamma\delta^{*})\uw^{G}_{G_{Y}}}$ is a constituent of $S_{Y, \gamma}\otimes S_{Y, \delta}^{*}$. But if $ \frac{\ch_{\delta}(h)}{\ch_{\delta}(1)}=\lam_{\cd}([Y], h)$ then $H_{\cd}$ acts as a identity on the representation $ \gamma\delta^{*}$. Therefore by Lemma \ref{krnind} $H_{\cd}$ acts as a identity also on $(\gamma\delta^{*}) \uw^{G}_{G_{Y}}$.  Thus $S_{1, (\gamma\delta^{*})\uw^{G}_{G_{Y}}} \in \cd$ and therefore $S_{Y, \delta}\in \cd$ since
\beqn
\Hom_{\cc^{G}}(S_{Y, \delta}, S_{1, (\gamma\delta^{*})\uw^{G}_{G_{Y}}}^{*} \ot S_{Y, \pi})\simeq \Hom_{\cc^{G}}(S_{1, (\gamma\delta^{*})\uw^{G}_{G_{Y}}}, S_{Y,\pi}\otimes S_{Y, \delta}^{*})\neq 0.
\eeqn
\epf

\bl \label{bicharacter} 
Let $\D$ be a fusion subcategory of $\cc^G$. Then $(\cs_{\cd}, H_{\cd}, \lam_{\cd})$ is a fusion datum as in Definition \ref{fusdatum}.
\el
\bpf By Lemma \ref{rot} it is clear that $\cs_{\cd}$ is $G$-stable. Moreover $H_{D}$ is a normal subgroup of $G$ and it acts trivially on $\Irr(\cs_{\cd})$ since $H_{\cd} \subset G_{Y}$ for all $Y \in \irr(\cs_{\cd})$. On the other hand since $H$ acts as a scalar on each projective representation $\pi$ with $S_{Y, \pi}\in \cd$ it follows that $\al_{Y}\vert_{H}$ is cohomologus trivial for any simple object $Y$ of $\cs_{\cd}$. This proves the first three conditions from the definition of a fusion datum.
\mdn
%{\bf Step 3.1: Twisted multiplicativity on the second term.}
Suppose that $S_{Y, \gm}\in \D$ and let $h,l \in H_{\cd}$. Since by Lemma \ref{zentre} $k_{\al_{[Y]}|_{H_{\cd}}}H_{\cd}$ acts as a scalar on $\gamma$  it follows that
\beqn (\frac{\ch_{\gamma}(hl)}{\ch_{\gamma}(1)})=\al_{[Y]}(h,l)(\frac{\ch_{\gamma}(h)}{\ch_{\gamma}(1)})(\frac{\ch_{\gamma}(l)}{\ch_{\gamma}(1)}).\eeqn Then
\beqarn
\lam_{\D}([Y], hl) &= & \frac{\ch_{\gamma}( hl)}{\ch_{\gamma}( 1)}=\frac{\ch_{\gamma}( h)}{\ch_{\gamma}( 1)}
\frac{\ch_{\gamma}( l)}{\ch_{\gamma}( 1)}\al_{[Y]}(h,l)\\ & = & \lam_{\D}([Y],h)\lam_{\D}([Y],l)\al_{[Y]}(h,l).
\eeqarn
for any $Y \in \irr(\cs_{\cd})$ and any $h,l \in H_{\cd}$. This proves Equation \eqref{hh'}.
\mdn
Now suppose that $U \in \Irr(\cs_{\cd})$ such that $\Hom_{\cc}(U, \ct^{t}(Y)\ot \ct^{s}(Z))\neq 0$. Then from the definition of $\cs_{\cd}$ there are $\gm,\; \gm', \; \delta$  irreducible projective representations of  respectively $G_Y,\;G_Z,\;G_{U}$ such that $S_{Y,\gm}, \; S_{Z,\gm'}\in \D$ and $\Hom_{\cc^G}(S_{U,\delta},\; S_{Y,\gm}\ot S_{Z,\gm'})\neq 0$. For example, by \cite[Lemma 2.15]{buna} one can take $\delta:=\Hom_{\cc}(U, S_{Y,\gm}\ot S_{Z,\gm'})$. From Corollary \ref{constituent} this implies that $\delta$ is a constituent of $(\;^t\gamma \ot \;^s\gamma' \ot \tau_U^{Y, Z}(s,t))\uw^{G_{U}}_{T}$.  Since by Lemma \ref{zentre} $H_{\cd}$ acts by the same scalar on each such constituent $\delta$ it follows that $H_{\cd}$ acts as a scalar on the whole induced representation $(\;^t\gamma \ot \;^s\gamma' \ot \tau_U^{Y, Z}(s,t))\uw^{G_{U}}_{T}$. On the other hand by Lemma \ref{a2} it follows that 
%On the other hand note that since $H_{\cd} \subset \core_G(G_U)\cap \core_G(G_Y)\cap \core_G(G_Z)$ this implies that
%\bq
%h(\bar{x}\ot_{k_{\al_{[U]}|_T}T}m)=\bar{x}\ot_{k_{\al_{[U]}|_T}T}hm
%\eq
$H_{\cd}$ acts as a scalar on $\;^t\pi \ot \;^s\dt \ot \tau_U^{Y, Z}(s,t)$.  Since $H_{\cd}$ acts also as a scalar on $\;^t\pi $ and $ \;^s\dt$ it follows that $H_{\cd}$ acts as a scalar on $\tau_U^{Y, Z}(t,s)$ which proves the fourth condition of a fusion datum.
\mdn
Then one has that
\beqn
\frac{\ch_\delta(h)}{\ch_\delta(1)}=\frac{\ch_{(\gm|_{G_U\cap G_Y\cap G_Z}^{G_Y}
\otimes \gm'|_{T}^{G_Z}\ot \tau_U^{Y,Z})}(h)}
{\ch_{(\gm|_{G_U\cap G_Y\cap G_Z}^{G_Y}\otimes \gm'|_{T}^{G_Z}\ot \tau_U^{Y,Z})}(1)}
=\frac{\ch_\gamma(h)}{\ch_\gm(1)}\frac{\ch_{\gamma'}(h)}{\ch_{\gm'}(1)}\frac{\ch_{\tau_U^{Y,Z}}(h)}{\ch_{\tau_U^{Y,Z}}(1)}
\eeqn
for all $h \in H_{\cd}$.
This shows that
\beqarn
\lam_{\D}([Y],h)\lam_{\D}([Z],h) & = & (\frac{\ch_\gamma(h)}{\ch_\gamma(1)}\frac{\ch_{\gamma'}(h)}{\ch_{\gamma'}(1)})
=\frac{\ch_\delta(h)}{\ch_\delta(1)}(\frac{\ch_{\tau_U^{Y,Z}}(h)}{\ch_{\tau_U^{Y,Z}}(1)})\inv\\ & = & \lam_{\D}([U],h)(\frac{\ch_{\tau_U^{Y,Z}}(h)}{\ch_{\tau_U^{Y,Z}}(1)})^{-1}
\eeqarn
which proves Equation \eqref{yz}.
\mdn
As noticed in Lemma \ref{rot} for any $g \in G$ one has that 
\beqn 
S_{Y, \pi}\cong S_{\ct^g(Y),\; \;^{g}D_{g, [Y]}( -)\;^g\pi},
\eeqn 
Therefore if $S_{Y, \pi}\cong S_{\ct^g(Y),\; \;^{g}D_{g, [Y]}( -)\;^g\pi}, \in \cd$ then one has that
\beqarn
\lam_{\cd}([\ct^g(Y)], h)&=&\frac{\;^{g}D_{g, [Y]}( h)\;^g\pi(h)}{\;^{g}D_{g, [Y]}( 1)\;^g\pi( 1)}=\frac{D_{g, [Y]}( g^{-1}hg)\ch_{\pi}( g^{-1}hg)}{\ch_{\pi}( 1)}\\&=&\lam_{\cd}([Y], g\inv h g)D_{g, [Y]}( g^{-1}hg)
\eeqarn which proves Equation \eqref{inv}.
\epf
\subsection{Construction of a fusion subcategory from a fusion datum}
Conversely suppose we have a fusion datum $(\cs, H,\lam)$ as in Definition \ref{fusdatum}. 
%\bab We fix this system $c$ and all objects will be constructed relative to this system. \eab 
\mdn To this datum we associate the fusion subcategory $\cc(\cs, H, \lam)$ generated by all the simple objects $S_{Y,\pi}\in \Irr(\cc)$ that satisfies
\beq\label{defi}
\frac{\ch_{\pi}(h)}{\ch_{\pi}(1)}=\lam([Y], h)
\eeq
for all $h \in H$.
%\mdn{This simple object depends on the chosen system $c_{Y}$ which should be $H$-trivial.} 
\begin{theorem}
 $\cc(\cs, H, \lam)$ as defined above is a fusion subcategory of $\cc^G$.
\end{theorem}

\bpf %Using the same arguments as in \cite{nat-faith} 
It is enough to show only that $\cc(\cs, H, \lam)$ is closed under tensor products.
Suppose that two simple objects $S_{Y, \pi}, S_{Z,\gm} \in \cc(\cs, H, \lam)$. We have to show that all the simple constituents $S_{U,\delta}$ of the tensor product
$S_{Y, \pi}\ot S_{Z,\gm}$ are also in $\cc(\cs, H, \lam)$. Using Corollary \ref{constituent} it follows that there are $t, s$ such that $U$ is a constituent of $\ct^{t}(Y)\ot \ct^{s}(Z)$ and
 $\Hom_T(\delta\vert_T,\;^t\pi\vert_T \otimes ^s\gamma\vert_T \otimes \tau_U^{Y, Z}(t,s))\neq 0$. Since $H$ acts as a scalar on the representation $^t\pi\vert_T \otimes ^s\gamma\vert_T \otimes \tau_U^{Y, Z}(s,t)$ it follows that
\bq \frac{\ch_{\dt}(h)}{\ch_{\dt}(1)}=\frac{\ch_{\;^t\pi}(h)}{\ch_{\;^t\pi}(1)}\frac{\ch_{\;^s\gm}(h)}{\ch_{\;^s\gm}(1)}\frac{\ch_{\tau_U^{Y, Z}(s,t)}(h)}{\ch_{\tau_U^{Y, Z}(s,t)}(1)}=\frac{\ch_{\pi}(t^{-1}ht)}{\ch_{\pi}(1)}\frac{\ch_{\gm}(s\inv hs)}{\ch_{\gm}(1)}\frac{\ch_{\tau_U^{Y, Z}(t,s)}(h)}{\ch_{\tau_U^{Y, Z}(s,t)}(1)}\eq
\mdn
Note that in this case by applying Equation \eqref{tauconj} one has:
\beqarn
\lam([U], h)&=&\lam([\ct^{t}(Y)], h)\lam([\ct^{s}(Z)], h)\frac{\ch_{\tau_{[U]}^{[\ct^{t}(Y)], [\ct^{s}(Z)]}}(h)}{\ch_{\tau_{[U]}^{[\ct^{t}(Y)], [\ct^{s}(Z)]}}(1)}\\&=& \lam([Y], t\inv ht)D_{t, [Y]}(t^{-1}ht)\lam([Z], t\inv ht)D_{s, [Z]}(s^{-1}hs)\frac{\ch_{\tau_{[U]}^{[\ct^{t}(Y)], [\ct^{s}(Z)]}}(h)}{\ch_{\tau_{[U]}^{[\ct^{t}(Y)], [\ct^{s}(Z)]}}(1)}\\&=& \lam([Y], t\inv ht)\lam([Z], t\inv ht)[D_{t, [Y]}(t^{-1}ht)D_{s, [Z]}(s^{-1}hs)\frac{\ch_{\tau_{[U]}^{[\ct^{t}(Y)], [\ct^{s}(Z)]}}(h)}{\ch_{\tau_{[U]}^{[\ct^{t}(Y)], [\ct^{s}(Z)]}}(1)}]\\&=& \frac{\ch_{\;^t\pi}( h)}{\ch_{\;^t\pi}( 1)}\frac{\ch_{\;^s\gm}(h)}{\ch_{\;^{s}\gm}(1)}\frac{\ch_{\tau_{[U]}^{[(Y)], [(Z)]}(t,s)}(h)}{\ch_{\tau_{[U]}^{[(Y)], [(Z)]}(t,s)}(1)}\\&=& \frac{\ch_{\dt}(h)}{\ch_{\dt}(1)}\eeqarn
We also used Equation \eqref{yz},  Equation \eqref{inv} and Proposition \ref{compst}
in the last equality. Thus $\lam([U],h)=\frac{\dt(h)}{\dt(1)}$ which implies that $S_{U,  \dt}\in \cc(\cs, H, \lam)$.\epf
%\noindent
%\bll{It is needed to see in the fusion rules form the published paper which system is used for $S_{U, \delta}$. \mdn For the $G$-invariance the argument that works for any $Y$ is that $S^{c}_{\ct^{l}(\ct^{g}(Y_{i})), \pi}\simeq S^{c}_{\ct^{gl}(Y_{i}), \;\pi}=S^{c}_{Y_{i}, \;^{gl}\pi}=S^{c}_{Y_{i}, \;^{g}(\;^{l}\pi)}=S^{c}_{\ct^{g}(Y_{i}), \;^{l}\pi}$.}
\subsection{Proof of Theorem \ref{cls-by-fusiondata}}
The proof of Theorem \ref{cls-by-fusiondata} follows from the following Propositions.
\begin{proposition} Let $(\Ss, H, \lambda)$ be a fusion datum and let $\cc(\Ss, H, \lambda) \subseteq \C^{G}$ be the associated fusion subcategory associated as above.
Then the fusion datum $(\Ss_\D, H_\D, \lambda_\D)$ corresponding to $\D:=\cc(\Ss, H, \lambda)$ coincides with $(\Ss, H, \lambda)$. \end{proposition}

\begin{proof}Straightforward. \end{proof}

\begin{proposition} Let $\D \subseteq \C$ be a fusion subcategory and let $(\Ss_\D, H_\D, \lambda_\D)$ be the fusion datum associated to $\D$.
 Then $\cc(\Ss_\D, H_\D, \lambda_\D) = \D$ \end{proposition}

\begin{proof}It is clear from Remark \ref{ews} that  $\cc(\Ss_\D, H_\D, \lambda_\D) \subseteq \D$. On the other hand the opposite inclusion follows from definitions  of $\lambda_\D$ and of the fusion datum $(\Ss_\D, H_\D, \lambda_\D)$.
\epf
\begin{corollary}\label{fpdd} With the above notations one has that $ \FPdim \cc(\Ss, H, \lambda)=\frac{ |G|\FPdim \cs}{|H|}.$
\end{corollary}
\bpf One has that
\beqarn
 \FPdim \cc(\Ss, H, \lambda)&=&\sum_{S_{Y, \pi} \in \cc(\Ss, H, \lambda)} \FPdim^{2}S_{Y, \pi}=\\&=&\sum_{S_{Y, \pi} \in \cc(\Ss, H, \lambda)}\frac{|G|^{2}}{|G_{Y}|^{2}}  \FPdim^{2} Y \dim^{2} \pi \\ &=& \sum_{Y \in \irr(\cs)/G} (\frac{|G|^{2}}{|G_{Y}|^{2}}  \FPdim^{2} Y)(\sum_{\pi \in \irr_{\al_{Y}}(G_{Y}), \;\pi|_{H}=\lam(Y, -)}\dim^{2}\pi)
\eeqarn
Then  \cite[Lemma 5.8]{nnw} shows that
\beqn
\sum_{\pi \in \irr_{\al_{Y}}(G_{Y}),\; \pi|_{H}=\lam(Y, -)}\dim^{2}\pi=\frac{|G_{Y}|}{|H|}
\eeqn and the equality $ \FPdim \cc(\Ss, H, \lambda)=\frac{ |G|\FPdim \cs}{|H|}$ follows from the equality $  \FPdim \cs= \sum_{Y \in \irr(\cs)/G} \frac{|G|}{|G_{Y}|}  \FPdim^{2} Y$.  
\end{proof}
\br\label{trivh}Note that if $H=1$ then necessarily $\lam:\irr(\cd)\ra k^{*}$ is trivial since $\lam ([Y], 1)=1$  for all $Y \in \irr( \cs)$. 
\er
\bp\label{eqf}
Any fusion subcategory of $\cc^{G}$ is also an equivariantization of a fusion subcategory of $\cc$ by a quotient subgroup of $G$.
\ep
\bpf
One has the exact sequence 
\beq
\rep G\ra \cc^{G} \ra \cc
\eeq
and its subsequence 
\beq
\rep (G/H) \ra \cc(H, \cs, \lam) \ra \cs
\eeq
which is {\it central} since the previous sequence is central.
It follows by \cite[Theorem 3.6]{tensor-exact2} that the second sequence is also an equivariantization, i.e. $\cc(H, \cs, \lam)=\cs^{G/H}$ for some action of $G/H$ on $\cs$.\epf
\subsection{On the lattice of fusion subcategories}
\bp\label{incl}
With the above notations one has that $\cc(\cs,H,\lam)\subset \cc(\cs',H',\lam')$ if and only if $\cs\subset \cs'$, $H \supset H'$ and $\lam|_{H'\times \Irr(S)}=\lam'|_{H'\times \Irr(S)}$.
\ep
\bpf
If $\cc(\cs,H,\lam)\subset \cc(\cs',H',\lam')$ then clearly $$\cs= F(\cc(\cs,H,\lam))\subset F(\cc(\cs',H',\lam'))=\cs'.$$ Moreover $\rep(G/H)=\rep (G) \cap \cc(\cs,H,\lam)\subset\rep(G)\cap  \cc(\cs',H',\lam')= \rep(G/H')$ which implies $H \supset H'$. Suppose that $S_{X, \pi}\in \cc(\cs, H, \lam)$. It follows that $\al_{X}|_{H'}=1$. Then definition of $\lam$ and $\lam'$ implies that
\beq
\lam([X],h)=\frac{\ch_{\pi(h)}}{\ch_{\pi(1)}}=\lam'([X],h)
\eeq
Thus the other equality $\lam|_{H'\times \Irr(S)}=\lam'|_{H'\times \Irr(S)}$ also holds. \mdn The converse follows directly from the definition of $\cc(\cs, H,\lam)$.
\epf

\mdn
Suppose there are given two fusion datum $(\cs, H, \lam)$ and $(\cs', H', \lam')$. Define 
\beq
\phi_{\lam, \lam'}:\irr(\cs \cap \cs') \times {H \cap H'}\ra k^{*}
\eeq
by $\phi_{\lam, \lam'}([Y], h) = \lam([Y], h)\lam'^{-1}([Y],h)$ for all $h \in H$.

Let also $\ker \phi_{\lam, \lam'}$ be the abelian subcategory of $\cs \cap \cs'$ consisting of those objects $[Y]\in \irr(\cs \cap \cs')$ such that $\phi_{\lam, \lam'}([X], h)=1$ for all $h \in H \cap H'$.
Equation \eqref{yz} implies that $\ker \phi_{\lam, \lam'}$ is a fusion subcategory of $\cc$. Define also
\beq
\psi_{\lam, \lam'}:\irr(\ker \phi_{\lam, \lam'})\times HH'\ra k^{*}
\eeq
defined on generators by
\beqn
\psi_{\lam, \lam'}([X], hh')=\lam([X], h)\lam'([X], h')
\eeqn
It is not difficult to see that $\psi_{\lam, \lam'}$ is well defined. Indeed, suppose that $hh'=ll'$ with $h,l\in H$ and $h',l'\in H'$. Then $l\inv h=l'{h'}\inv \in H \cap H'$. Thus for any $X \in \Irr\ker  \phi_{\lam, \lam'}$ one has that:
\beqn
\lam([X], l)\inv \lam([X],  h)=\lam([X], l\inv h)=\lam([X], l'{h'}\inv )=\lam([X], l' )\lam([X], {h'})\inv 
\eeqn
which shows that
\beqn
\lam([X], h)\lam([X], h')=\lam([X], l)\lam([X], l').
\eeqn
i.e. $\psi_{\lam, \lam'}([X], -)$ is well defined.
%\blue{It is easy to see that $\psi_{\lam, \lam'}$ is well defined. Indeed suppose that $hh'=h_{1}h_{1}'$ with $h, h_{1} \in H$ and $h' , h'_{1} \in H'$. Then by Equation \eqref{hh'} one has that $\lam ([X], h_{1}^{-1}h)=\lam([X], h)\lam([X], h_{1}^{-1})\al_{X}( h_{1}^{-1}, h)$. On the other hand \beqn\lam' (X, h_{1}'h^{'-1})=\lam'([X], h_{1}')\lam'([X], h^{'-1})\al_{X}( h_{1}', h^{'-1}).\eeqn}
\bp\label{int}
\wtan one has that $\cc(H,\cs,\lam)\cap \;\cc(H', \cs',\lam')=\cc(\ker \phi_{\lam, \lam'}, HH',\psi_{\lam, \lam'})$.\ep

\bpf Suppose that  $\cc(\cs,H,\lam)\cap \;\cc(\cs',H',\lam')=\cc(\cs'', H'', \lam'')$.
If $S_{X , \pi} \in \cc(\cs,H,\lam)\cap \cc(\cs',H',\lam')$ then  one has that 
 $\lam([X], h)=\frac{\ch_{\pi}(h)}{\ch_{\pi}(1)}=\lam'([X, ], h)$ which shows that $X \in \irr(\ker \phi_{\lam, \lam'})$. Thus $\cs''\subseteq \ker \phi_{\lam, \lam'}$.\mdn
Moreover note that $\Rep(G)\cap (\cc(\cs,H,\lam)\cap \cc(\cs',H',\lam'))=\Rep(G/H)\cap \Rep(G/H')=\Rep(G/HH')$ which implies that $H''=HH'$.  Then by its definition one has  $\lam''([X], hh')=\frac{\ch_{\pi}(hh')}{\ch_{\pi}(1)}$. Using Lemma \ref{zentre} it follows that $\lam''([X], hh')=\al_{[X]}(h, h')\frac{\ch_{\pi}(h)}{\ch_{\pi}(1)}\frac{\ch_{\pi}(h')}{\ch_{\pi}(1)}=\al_{[X]}(h, h')\lam([X], h)\lam'([X], h')=\psi_{\lam, \lam'}([X],h)$  which also shows that  $\psi_{\lam, \lam'}$ is well defined.
\epf

\bp \wtan one has that $\cc(\cs,H,\lam)\vee \cc(\cs',H',\lam')=\cc(\cs \vee \cs', H_{\vee},\lam_{\vee})$ where 
\beqn
H_{\vee}=\bigcap_{X \in \irr(\cs\cap \cs')}\ker_{G}[\lam([X], - )\uw^{G}_{H}\overline{\lam'([X^{*}], -)}\uw^{G}_{H'}]\eeqn
and 
\beqn\lam_{\vee}([U],h):=\lam([Y],h)\lam'([Z], h)\eeqn for all $Y \in \irr(\cs)$, $Z\in \irr(\cs')$ and all simple objects $U$ with  $\Hom_{\cc}(U,Y\ot Z)\neq 0$.
\ep
\bpf
Clearly $F(\cc(\cs,H,\lam)\vee \cc(\cs',H',\lam'))=\cs \vee \cs'$. Moreover $\Rep(G/H\cap H')=\Rep(G/H)\vee \Rep(G/H')\subset \cc(\cs,H,\lam)\vee \cc(\cs',H',\lam')$. On the other hand if $S_{1, M} \in \cc_{1}\vee \cc_{2}$ then $S_{1, M}$ appears as constituent of $S_{X, \pi} \ot S_{X^{*}, \delta}$. It follows that $M$ is a constituent of $(\pi\delta^{*})\uw^{G}_{G_{X}}$. This implies that
\beqn
H_{\vee}=\bigcap_{X \in \cs_{1}\cap \cs_{2}}\ker_{G}(\lam([X], - )\uw^{G}_{H}\overline{\lam'([X^{*}], -)}\uw^{G}_{H'}).\eeqn
Formula for $\lam_{\vee}$ follows from definition.
\epf

\br
Note that in the previous Proposition one has $H_{\vee}\subseteq H \bigcap H'$.
\er
\begin{corollary}\label{veecom}
Suppose that $\cc$ is a  braided $G$-crossed fusion category with $K_{0}(\cc)$ commutative. Then $H_{\vee}=H_{1}\cap H_{2}$.
\end{corollary}
%\bll{And also for the fusion $\cs_{i}$ the same relation. One passes to the pf dimension.}
%\bll{For $G$-crossed fusion categories both relations are satisfied.}
\bpf
It follows from the fact that in a braided fusion category one has that $\fp\; \cc_{1} \vee \cc_{2}=\frac{\fp \cc_{1}\; \fp \cc_{2}}{\fp\;\cc_{1} \cap \cc_{2}}$, see \cite{DGNOI}.
\epf
%\bb{For the trivial grading $\lam(h)=1$ for all $h \in H$ write $\cc(\cs, H)$ for $\cc(\cs, H,1)$. When is this possible?}
%\bb{\brIf $\cs$ contains $\rep(G)$ then $\cs$ is of the type $\cs(\cd)$ for some fusion subcategory of $\cc$.
%\er}
\section{Equivariantization of braided $G$-crossed fusion categories}\label{emc}

In this section we prove Theorem \ref{form}. It entails a description for M\"{u}ger's centralizer of all fusion subcategories of a group theoretical braided fusion category. %\bll{Note that this can also be obtained from Naidu's result.}
\mdn
\subsection{Modular fusion categories} Recall that a {\it braided tensor category} $\cC$ is a tensor category equipped for any two objects $X,Y \in \cc$ with natural isomorphisms $c_{X,\;Y} : X \otimes Y \ra  Y \otimes X$  satisfying the two braiding hexagon axioms, see for example \cite{BaKi, JS}.\mdn By \cite[Proposition 2.1]{JS} it follows that for any object $A\in \cc$ one has that $c_{A, I}=l_{A}r_{A}^{-1}$ and $c_{I, A}=l_{A}^{-1}r_{A}$ where $l_{A}:A\ot I\ra A$ and $r_{A}:I\otimes A \ra A$ are the canonical isomorphisms associated to the unit $I$ of $\cc$. %\subsubsection*{Ribbon fusion categories}
\mdn
A {\it twist} on a braided fusion category $\cc$ is a natural automorphism $\tet:\id_{\cc}\ra \id_{\cc}$ satisfying $\tet_{1}=\id_{1}$
and
\beq
\tet_{X\ot Y}=(\tet_{X}\ot \tet_{Y})c_{YX}c_{XY}.
\eeq
A braided fusion category is called {\it premodular} or {\it ribbon} if it has a twist satisfying $\tet_{X^{*}} =\tet_{X}^{*}$ for all $X \in \cc$, see \cite{br2}. By \cite{DGNOI} a premodular category is a braided fusion category endowed with a spherical structure. Recall that the entries of the $S$-matrix,  $S=\{s_{X,Y}\}$ of a premodular category are defined as the quantum trace $s_{X,Y}:=tr_{q}(c_{YX}c_{XY})$, see \cite{BaKi}. A premodular category $\cc$ is called {\it  modular} if the above $S$-matrix is nondegenerate, see \cite{tuv0}.

\subsection*{Centralizers in braided fusion categories}
Let $\cc$ be a braided fusion category and $\cd$ be a fusion subcategory of $\cc$. M\"uger introduced the notion of {\it centralizer of $\cd$ } as the fusion subcategory $\cd'$ of $\cc$ generated by all simple objects $X$ of $\cc$ satisfying  
\beq\lb{centr} c_{X,\; Y}c_{Y,\;X}=\mtr{id}_{X\otimes Y}\eeq
for all objects $Y \in\cd$.% (see \cite{proclond}).
\mdn
In the case of a ribbon category $\cc$ the condition of Equation \eqref{centr} is equivalent to 
\beq\lb{centr2} s_{X,\;Y}=\fp(X)\fp(Y)\eeq for all objects $Y \in \cd$.
Note that  in  general \beq\label{is}|s_{X,\;Y}|\leq \fp(X)\fp(Y)\eeq  by \cite[Proposition 2.5]{proclond}. In the situation of equality in Equation \eqref{centr} we say that the objects $X$ and $Y$ centralize each other.
\mdn We also say that $X, Y$ {\it projectively centralize} each other if and only if there is $\omega \in k^{*}$ a root of unity such that  $s_{X,\;Y}=\omega\fp(X)\fp(Y)$.
\mdn
Let $\cc$ be a braided fusion category. We say that a Tannakian subcategory
$\ce\subset \cd$ is {\it Lagrangian} if $\ce'=\ce$.\mdn
Recall that the {\it M\"{u}ger center} of a braided fusion category is the M\"{u}ger centralizer $\cz_{2}(\cc)=\cc'$ of the whole category $\cc$. A braided fusion category is called {\it nondegenerate } if its M\"{u}ger center is trivial.
%\subsubsection*{\bf Properties of the centralizer}
If $\cC$ is a braided fusion category then by \cite[Theorem 3.23, Theorem 3.24]{DGNOI} one has that \beq\label{fpdc}\fp(\cd)\fp(\cd')=\fp (\cC)\fp(\cd\cap \cz_{2}(\cc)).\eeq
Moreover if $\cd$ is nondegenerate then there is an equivalence of braided fusion categories $\cc \simeq \cd \boxtimes \cd'$. In this case $\cc$ is nondegenerate if and only if $\cd'$ is. On the other hand by \cite[Corollary 3.11]{DGNOI} one has that $\cd''=\cd\vee \cz_{2}(\cc)$. In particular if $\cc$ is nondegenerate then $\cd''=\cd$.
Note that for modular fusion categories the above properties of the centralizers were proven in \cite{proclond, subf2}.
\mdn
Let $\cD$ and $\cE$ be fusion subcategories of a modular fusion category $\cC$. Then one also has 
\beq\lb{centrofint} (\cD\vee\cE)'=\cD'\cap \cE' \text{\;and \;} (\cD\cap \cE)'=\cD'\vee \cE',\eeq see \cite{proclond}.
%\subsubsection{Properties of the $S$-matrix}
%$$s_{X,\;Y}=\teta^{-1}_X\teta^{-1}_Y\sum_{W\in \Lam_{\cc}}N^W_{XY}\teta_Wd(W)$$
%$$s_{X,\;Y} s_{X,\;Z} =d(X)\sum_{W\in \Lam_{\cc}}N^W_{YZ} s_{X,\;W},$$ for all $X,Y,Z\in \Lam_{\cc}$.
A braided fusion category is called {\it slightly nondegenerate } if its M\"{u}ger center is braided equivalent to the category SuperVec of super vector spaces.
\mdn
Recall that a braided fusion category $\cc$ is called {\it symmetric} if $\cc'=\cc$. Equivalently, $c_{X,Y}c_{Y,X}=\id_{X\otimes Y}$ for all objects $X, Y\in \cc$. In this case the braiding $c$ is also called symmetric.
\subsection{Braided $G$-crossed fusion categories}
In this section we recall a description (due to Kirillov. Jr. \cite{orbif} and M\"{u}ger \cite{gtm}) of the structure of braided fusion categories containing a subcategory $\rep(G)$ for some finite group $G$ as a braided fusion subcategory. They are described as equivariantizations of {\it braided $G$-crossed categories.} which were introduced by Turaev \cite{tuv1}.
\mdn
Recall that a braided $G$-crossed fusion category consists of the following data:
\mdn
$(1)$ 
$G$ a finite group,
\mdn $(2)$
$\cC$ a fusion category with a (not necessarily faithful) $G$-grading 
 \beqn
 C = \oplus_{g\in G}
\;\cC_g,
 \eeqn
\mdn $(3)$ an action of $G$ $\ct: \underline G \to \underline \Aut_{\otimes} \C$  on $\cc$  by tensor autoequivalences, satisfying 
 $\ct_g(\cC_h) \subset \cC_{ghg^{-1}}$.
 %, withtensor structure $\gamma$ and unit-preserving structure $i$,
\mdn $(4)$ a family of natural isomorphisms,
$c_{X, Y  }: X \ot Y \xra{\cong} \ct^g(Y ) \ot X$ for any  $X \in \cc_g$, $Y \in \cc$,  called {\it $G$-braidings},
satisfying the following compatibility conditions:

{\Small
\beq\label{c1}
 \;\;\;\;\;\;\;\;((\ct_{2}^{g,h})^{-1}_{Y} \ot  \id_{T_{g}(X)})\circ((\ct_{2}^{ghg^{-1},g})_{Y}
\ot \id_{T_{g}(X)}) \circ c(\ct^{g}(X), \ct^{g}(Y )) \circ (\ct^{g}_2)_{X, Y }
= \eeq \beqn = (\ct^{g}_2)_{\ct^{h}(Y ),X} \circ \ct^{g}(c(X, Y )),
\eeqn
}
for all $g, h \in G$ and objects $X \in \cc_h, Y \in \cc$.
{\Small \beq\label{c22}
\;\;\;\;\; \al^{-1}_{\ct^g(\ct^{h}(Z)),X,Y} \circ((\ct^{2}_{g,h})^{-1}_{Z} \ot \id_{X\ot Y }) \circ c(X \ot Y,Z) \circ \al^{-1}_{
X,Y,Z}= \eeq \beqn =(c(X, \ct^{h}(Z)) \ot \id_Y )  \circ \al^{-1}_{
X,\ct^{h}(Z),Y} \circ (\id_X \ot c(Y,Z)),
\eeqn}
for all $g, h \in G$ and objects $X \in \cc_g, Y \in \cc_h, Z \in \cc$.
{\Small \beq\label{c3}
\;\;\;\;\;\;\;\;\;\;\al_{\ct^g(Y ),\ct^g(Z),X}\circ ((\ct^{g}_{2})_{Y,Z} \ot \id_X) \circ c(X, Y \ot Z) \circ \al_{X,Y,Z}\eeq \beqn =(\id_{\ct^g(Y )} \ot c(X,Z)) \circ \al_{\ct^g(Y ),X,Z }\circ (c(X, Y ) \ot \id_Z),
\eeqn
}
It is well known that in this case $c_{\1, X}=\id_{X}=c_{X, \;1}$ for any object $X$ of $\cc$. \mdn
For any simple object $X$ of $\cc$ we define the degree  $\deg(X) \in G$  of $X$, to be the element that satisfies $X \in \cc_{\deg(X)}.$

 \br Note that the trivial component $\cc_1$ of a braided $G$-crossed fusion category
$\cc$ is a braided fusion category with the action of $G$ on it given by {\it braided autoequivalences.}
This can be seen by taking $X, Y \in \cC_1$ in the second and the third compatibility condition from above.
\er\noindent
It is well known that in this case $\cc^{G}$ is a braided fusion category with the braiding
$
\tilde{c}_{S, T}:S\ot T\ra T\ot S\; 
$ given by
{\Small \beq\label{br}
S\ot T=\bigoplus_{g \in K_{S}}S_{g} \ot T \xra{\oplus_{g \in K_{S}}c_{ _{S_{g}, T}}} \bigoplus_{g \in K_{S}}\ct^{g}(T) \ot S \xra{\oplus_{g \in K_{S}}(\mu^{g}_{T} \ot 1)} \bigoplus_{g \in K_{S}}T \ot S_{g}=T\ot S.
\eeq}
for all $S, T \in \cc^{G}$.
\section{M\"{u}ger centralizer for the equivariantization of  a braided $G$-crossed fusion categories}\label{mc}
%\subsection{General settings}
Suppose that
$
\cc=\oplus_{g \in G}\;\cc_{g}
$
is a braided $G$-crossed fusion category with $\cd=\cc_{1}$. Let $G_{1}:=\{g \in G\;|\;\cc_{g}\neq 0\}$. Then clearly  $G_{1}$ is a normal subgroup of $G$.
\mdn For any object $S=\oplus_{g\in G}S_{g} \in \cc$ define its support $K_{S}$ as $K_{S}:=\{g \in G\;|\;S_{g} \neq 0\}$. Let $\cc^{G}$ be the equivariantized fusion category and $F:\cc^{G}\ra \cc$ be the forgetful functor. By abuse of notations sometimes we identify $F(X)$ with $X$ for any object  $X\in \cc^{G}$. For this reason, for any $S \in \cc^{G}$ we also sometimes denote by $K_{S}$ the support of the object $F(S)\in \cc$.

\br\label{inclgrps} Let $M=\oplus_{g \in K_{M}}M_{g}$ be an object of $\cc^{G}$. Since $\mu^{g}_{T}:\ct^{g}(M)\ra M$ is an isomorphism in $\cc$ it follows that  $\ct^{g}(M_{h})\simeq M_{ghg\inv}$ by the morphism $\mu^{g}_{M}|_{M_{h}}$.
\er 
%\mdn \bll{For an action of $G$ on $\cc$ recall also the inertia subgroup $G_{Y}=\{g \in G\;|\ct^g(Y) \cong Y\}$ for any object $Y \in \cc$.
%\md Moreover one has $\rep(G) \subset \cc^{G}$ as in Equation \eqref{sec}.}
\bl\label{commuting}
Suppose that $\cc$ is a braided $G$-crossed fusion category as above. Let $S=\bigoplus_{g \in K_{S}}S_{g}$ and $T=\bigoplus_{g \in K_{T}}T_{g}$ be two objects of $\cc^{G}$ with $S_{g}, T_{g}\in \cc_{g}$. Then $S$ and $T$ projectively centralize each other in $\cc^{G}$ if and only if the following conditions are satisfied:
\bne 
%\item \beq\label{incls} K_{S}\subset \bigcap_{h \in K_{T}}G_{T_{h}}\text{\;and\;}  K_{T}\subset \bigcap_{g \in K_{S}}G_{S_{g}} \eeq 
\item
$K_{S}$ commutes elementwise with $K_{T}$.

\item There is a root of unity $\omega \in k^{*}$ such that:
\beq\label{compo}
(\mu^{h}_{S}|_{\ct^{h}(S_{g})}\ot \id_{T_{h}})c_{T_{h},S_{g}}(\mu^{g}_{T}|_{\ct^{g}(T_{h})}\ot \id_{S_{g}})c_{S_{g}, T_{h}}=\omega\id_{S_{g}\ot T_{h}}
\eeq
for all $g \in K_{S}$ and $h \in K_{T}$.
\ene
\el

\bpf
%Suppose that $\cc=\oplus_{g \in G}\cc_{g}$ is a   braided $G$-crossed   fusion category with the crossed braiding $c_{X, Y}:X\ot T \xra{c_{X, Y}} T^{g}(Y) \ot S$ with $X \in \cc_{g}$.  Then 
From Equation \eqref{br} the braiding $\tilde{c}_{S,T}$ between $S$ and $T$ can be written on the components as following:
\beq
\widetilde c^{\;g, h}_{S, T}:S_g\ot T_h \xra{c_{S_g, T_h}} 
\ct^g(T_h)\ot S_g \xra{\mu^g_{T}|_{{\ct^g(T_h)}}\ot \id_{S_g}} T_{ghg^{-1}} \ot S_g
\eeq
\noindent
It follows that the braiding $\tilde{c}_{S, T}$ sends the homogenous subobject $S_{g}\ot T_{h}$ of $S \ot T$ into the homogenous subobject $T_{ghg^{-1}}\ot S_{g}$ of $T \ot S$. Applying now $\widetilde{c}_{T, S}$ one has that $\widetilde c^{ghg^{-1}, g}_{T, S}$ is given by
{\Tiny \beq 
T_{ghg^{-1}}\ot S_g \xra{c_{_{ T_{ghg^{-1}}, S_g}}}\ct^{ghg^{-1}}(S_g)\ot T_{ghg^{-1}} \xra{{\mu^{ghg^{-1}}_{S}|_{{\ct^{ghg^{-1}}(S_g)}}}\ot \id_{T_{ghg^{-1}}}} S_{ghgh^{-1}g^{-1}}\ot T_{ghg^{-1}}
\eeq}
Thus  the homogenous subobject $T_{ghg^{-1}}\ot S_{g}$ of $T \ot S$ is sent by $\widetilde{c}_{T, S}$ to $S_{(ghg^{-1})h(ghg^{-1})^{-1}}\ot T_{ghg^{-1}}$. Since the composition of these two maps should be a scalar multiple of the identity it follows that $ghg^{-1}=h$ i.e $g$ commutes to $h$. This shows that $K_{S}$ and $K_{T}$ commutes elementwise.
\mdn
Note that imposing that the composition $\widetilde c_{S, T} \widetilde c_{T, S}$ of the two braidings being $\omega\id_{S \ot T}$ it implies the other equation.
\epf
\br Suppose that $S=\oplus_{g\in K_{S}}S_{g}$ projectively centralizes $T=\oplus_{h\in K_{T}}T_{h}$. It follows from previous lemma that $gh=hg$ for any $g\in K_{S}$ and any $h\in K_{T}$.  \onh\; Remark \ref{inclgrps} implies that
$\mu^{g}_{T}|_{\ct^{g}(T_{h})}:\ct^{g}(T_{h})\ra T_{h}$ is an isomorphism in $\cc$ which implies  that $g \in G_{T_h}$.  Thus $K_{S}\subset \cap_{h \in K_{T}}G_{T_{h}}$. Similarly $K_{T}\subset \cap_{g \in K_{S}}G_{S_{g}}$ for all $g \in K_{S}$.\er

\br
The proof of the previous lemma also shows that $S$ and $T$ centralize each other if and only if $K_{S}$ commutes elementwise with $K_{T}$ and 
\beq\label{compo2}
(\mu^{v}_{S}|_{\ct^{v}(S_{u})}\ot \id_{T_{v}})c_{T_{v},S_{u}}(\mu^{u}_{T}|_{\ct^{u}(T_{v})}\ot \id_{S_{u}})c_{S_{u}, T_{v}}=\id_{S_{u}\ot T_{v}}
\eeq
for all $u \in K_{S}$ and $v \in K_{T}$.
\er
%\bab \noindent Can one do a similar thing to compute te dyslectic modules?\eab
\subsection{Two subgroups associated to a fusion subcategory of $\cc^{G}$}
For a fusion subcategory $\cd \subset \cc^{G}$ define $H_{\cd}$, the normal subgroup of $G$ such that $\rep(G)\cap \cd=\rep(G/H_{\cd})$. Let also $K_{\cd}$ be the subgroup generated by the support of all simple objects of $\cd$. 
\bl\label{c2} Suppose that $\cc$ is a braided $G$-crossed fusion category as above and let $X$ be an object of $\cc^{G}$. Then for any $V \in \rep(G)$ one has that $S_{\1, V}$ centralizes $X$ if and only if $K_{X} \subset \ker_{G}(V) $.
\el
\bpf Let $\rho:G \ra \End_{k}(V)$ be  the representation of $G$ associated to $V$.
Recall that $S_{\1, V}:=V\ot \1$ is a on object of $\cc^{G}$ with the equivariant structure $\ct^{g}(V \ot \1) \xra{\rho(g) \ot 1} V \ot \1$. On the other hand since $\ct^{1}=\id_{\cc}$ note that the braiding 
\beqn
\tilde{c}_{S_{\1,V},\;X}:S_{\1,V}\ot X \xra{c_{V \ot \1, X} }\ct^{1}(X)\ot S_{\1,V} \ra X \ot S_{\1,V}
\eeqn
is the identity morphism on $X$, via the identification $l_{X}r^{-1}_{X}$.
\mdn
Also, since $\ct^{g}(\1)=\1$ the braiding $\tilde{c}_{X, (V\ot \1)}:X\ot (V\ot \1) \ra \ct^{g}(V\ot \1)\ot X=(V\ot \1) \ot X$
is given on the homogenous component $X_{g}\ot V$ by multiplication with $g$ on $V$.
Thus the double braiding 
\beq
c_{ _{S_{\1, V}, X}}c_{ _{X, S_{\1, V}}}:X\ot S_{\1, V} \ra S_{\1, V}\ot X
\eeq
regarded as an endomorphism of $X \ot V$ is given on the component $X_{g }\ot V$ just by multiplication by $g$ on $V$. This implies that the two above objects centralize each other if and only if $K_{X}$ acts as identity on $V$, i.e. $K_{X} \subset \ker_{G}(V) $.
\epf
\begin{corollary}\label{c2} Suppose that $Y\in \cc_{g}$ is a simple object of $\cc$. Then one has that $S_{\1, \;V}$ centralizes $S_{Y, \pi}$ if and only if $g \in \ker_{G}(V) $.
\end{corollary}
\bpf
Note that $G_{Y }\subseteq C_{G}(g)$ and therefore the support of $S_{Y, \;\pi}$ coincides to $K_{g}$, the conjugacy class of $G$.
\epf

\bn{corollary}\label{czc1}
One has that $S_{1, V}\in \cz_{2}(\cc^{G})$ if and only if $G_{1}\subseteq \ker_{G}(V)$. In particular $S_{1, G/G_{1}}\in \cz_{2}(\cc^{G})$ and $H_{\cz_{2}(\cc^{G})}=G_{1}$.
\end{corollary}

\begin{theorem}\label{switch} Let $\cc$ be a  braided $G$-crossed fusion category and let $\cd$ be a fusion subcategory of $\cc^{G}$.
Then with the above notations one has that $K_{\cd'}$ and $K_{\cd}$ commute elementwise and the following equalities hold $H_{\cd'}=K_{\cd}$. Moreover  $K_{\cd'}=H_{\cd''}\subseteq H_{\cd}\cap G_{1}$.\end{theorem}

\bpf
Lemma \ref{commuting} implies that $K_{\cd}$ and $K_{\cd'} $ commutes elementwise. On the other hand by Lemma \ref{c2} one has that $S_{\1, V} \in \rep(G) \cap \cd'$ if and only $K_{\cd}\subset \ker_{G}(V)$. This implies that $\rep(G) \cap \cd' =\rep(G/K_{\cd})$ and therefore $H_{\cd'}=K_{\cd}$. Then $K_{\cd'}=H_{\cd''}\subseteq H_{\cd}\cap G_{1}$.
\epf
\bp\label{nondeg}
Suppose that $\cc$ is a braided $G$-crossed fusion category with a faithful grading by $G$ and $\cc_{1}$ a  nondegenerate braided fusion category. Then for any fusion subcategory $\cd \subseteq \cc^{G}$ one has that $H_{\cd}$ commutes elementwise to $K_{\cd}$ and $H_{\cd}=K_{\cd'}$.
\ep

\bpf Applying \cite[Proposition 4.56]{DGNOI} it follows that $\cc^{G}$ is braided nondegenerate.
Therefore $\cd''=\cd$ and applying the results of the previous theorem for $\cd'$ the conclusion follows.
\epf\noindent

\br
The above proposition generalizes \cite[Proposition 3.6]{nnw} and \cite[Proposition 5.2]{nnw} to equivariantizations of any braided $G$-crossed fusion category.
\er
\subsection{On the $d$-morphisms}\label{dd} For any element $g \in G$ denote by $K_{g}$ the conjugacy class of $g$ in $G$.
\mdn%\subsubsection{} 
Suppose that are given two simple objects $X\in \cc_{g}$ and $Y \in \cc_{h}$ such that $K_{g}$ commutes elementwise to $K_{h}$. For all $ a,b \in G$ consider the following isomorphism  $ d_{X, Y}^{a,b}:\ct^{a}(X) \ot \ct^{b}(Y) \ra \ct^{b}(Y) \ot \ct^{a}(X)$ in $\cc$ given by
\beqn
d_{X, Y}^{a,b}:=(\ct^{b}(\ct_{2}^{aga\inv, b})_{Y}({\ct_{2}^{b, b^{-1}aga\inv b})_{Y}}^{-1}(c^{b^{-1}aga\inv b}_{Y})\ot \id_{\ct^{a}(X)})c_{ \ct^{a}(X),\ct^{b}(Y)}
\eeqn
\mdn Recall that we fixed a set of representatives for the isomorphism classes of simple objects of $\cc$.
Suppose  that $X, Y$ are simple objects  as above such that $X \simeq Y_{i}$ and $Y\simeq Y_{j}$. Note that if $u:X \ra Y_{i}$ and $v:Y\ra Y_{j}
$ are isomorphisms then
\beq
(\ct^{b}(v)\ot \ct^{a}(u))d_{Y_{i}, Y_{j}}^{a,b}=d_{X, Y}^{a,b}(\ct^{a}(u)\ot \ct^{b}(v)).
\eeq
Indeed, the above relation is equivalent to the commutativity of the  diagram below. Note that the last rectangle is commutative by Equation \eqref{isocho}. Moreover, the other rectangles are commutative by the natural properties of the group action and Equation \eqref{ro-2}. 
{\tiny
\begin{center}
\begin{tikzpicture}
  \matrix (m) [matrix of math nodes, row sep=5em,column sep=3 em,minimum width=0.5 em]
  { \ct^{a}(X)\ot \ct^{b}(Y) & \ct^{a}(Y_{i})\ot \ct^{b}(Y_{j}))   
%line 2
\\  \ct^{aga\inv}( \ct^{b}(Y)) \ot \ct^{a}(X)  & \ct^{aga\inv}( \ct^{b}(Y_{j})) \ot \ct^{a}(Y_{i})
 %line 3
   \\ \ct^{aga^{-1}b}(Y) \ot \ct^{ab}(X) & \ct^{aga^{-1}b}(Y_{j}) \ot \ct^{ab}(Y_{i})
      % line 4
   \\ \ct^{b}(\ct^{(b^{-1}a)g(b^{-1}a)^{-1}}(Y))\ot \ct^{a}(X)) & \ct^{b}(\ct^{(b^{-1}a)g(b^{-1}a)^{-1}}(Y_{j}))\ot \ct^{a}(Y_{i}))   % line 5
   \\\ct^{b}(Y) \ot \ct^{a}(X)) & \ct^{b}(Y_{j}) \ot \ct^{a}(Y_{i})). \\}; \path[-stealth]
  (m-1-1) edge node [above] {$\ct^{a}(u)\ot \ct^{b}(v)$} (m-1-2)

    (m-2-1) edge [dashed] node   [above] {$$} (m-2-2)
       (m-3-1) edge [dashed] node   [above] {$$} (m-3-2)
     (m-4-1) edge [dashed] node   [above] {$$} (m-4-2)
    
     (m-5-1) edge  node   [below] {$ \ct^{b}(v) \ot \ct^{a}(u)$} (m-5-2)
    
         %columns
               (m-1-1) edge node [left] {$c(\ct^{a}(X), \ct^{b}(Y))$} (m-2-1)
    
   (m-1-2) edge node [right] {$c(\ct^{a}(Y_{i}), \ct^{b}(Y_{j}))$} (m-2-2)     
   (m-2-1) edge node [left] {$(\ct_{2}^{aga^{-1},b})_{Y}\ot 1$} (m-3-1)
    (m-2-2) edge node [right] {$(\ct_{2}^{aga^{-1},b})_{Y_{j}}\ot 1$} (m-3-2)
    
     (m-3-1) edge node [left] {$(\ct_{2}^{b, b^{-1}aga^{-1}b})^{-1}_{Y}\ot 1$} (m-4-1)
    (m-3-2) edge node [right] {$(\ct_{2}^{b, b^{-1}aga^{-1}b})^{-1}_{Y_{j}}\ot 1$} (m-4-2)
   
     (m-4-1) edge node [left] {$\ct^{b}(c_{Y}^{(b^{-1}a)g(b^{-1}a)^{-1}})\ot 1$} (m-5-1)
    (m-4-2) edge node [right] {$\ct^{b}(c_{Y_{j}}^{(b^{-1}a)g(b^{-1}a)^{-1}})\ot 1$} (m-5-2)
  
   ;
\end{tikzpicture}
  \end{center}
}
\mdn
Similarly one has that 
\beq
(\ct^{b}(u)\ot \ct^{a}(v))d_{{Y_j}, Y_{i}}^{b,a}=d_{Y, X}^{b,a}(\ct^{a}(v)\ot \ct^{b}(u)).
\eeq
Combining the above two equations one has that
\beq\label{isochocomp}
(\ct^{a}(u)\ot \ct^{b}(v))d^{a,b}_{Y_{i}, Y_{j}}d^{b,a}_{Y_{j}, Y_{i}}=d^{a,b}_{X,Y}d^{b,a}_{Y,X}(\ct^{a}(u)\ot \ct^{b}(v))
\eeq
\bl \label{cdcomp1} With the above notations one has that 
\beq
\ct^{a} (d_{X, Y}^{b,c}) (\ct_{2}^{a})^{T^{c}(Y), T^{b}(X)}= (\ct_{2}^{a})^{T^{b}(X), T^{c}(Y)}d_{X, Y}^{ab, ac}
\eeq
\el
\bpf Note that the identity is equivalent to the commutativity of the diagram below. In this diagram the first column represents $\ct^{a}(d_{X,Y}^{b,c})$ while the last column represents $d_{X, Y}^{ab, ac}$.
\mdn 
The above rectangle is commutative by the first condition \eqref{c1} in the definition of a crossed braided fusion category. All the other squares are commutative by the properties of the group action by tensor equivalences. 
{\Tiny
\begin{center}
\begin{tikzpicture}
  \matrix (m) [matrix of math nodes, row sep=4em,column sep=2 em,minimum width=0.5 em]
  { \ct^{a}(\ct^{b}(X)\ot \ct^{c}(Y)) & \ct^{a}(\ct^{b}(X))\ot \ct^{a}(\ct^{c}(Y)) & \ct^{ab}(X) \ot \ct^{ac}(Y)   
%line 2
\\  \ct^{a}( \ct^{bgb^{-1}}(\ct^{c}(Y)) \ot \ct^{b}(X) ) & \ct^{a}(\ct^{bgb^{-1}}((\ct^{c}(Y))\ot \ct^{a}(\ct^{b}(X)) & \ct^{abgb^{-1}a^{-1}}(\ct^{ac}(Y) )\ot \ct^{ab}(X)  
 %line 3
   \\ \ct^{a}(\ct^{bgb^{-1}c}(Y) \ot \ct^{b}(X)) & \ct^{a}( \ct^{bgb^{-1}c}(Y))\ot \ct^{a}(\ct^{b}(X)) & \ct^{abgb^{-1}c}(Y)  \ot \ct^{ab}(X)
   % line 4
   \\ \ct^{a}(\ct^{c}(\ct^{(b^{-1}c)^{-1}g(b^{-1}c)}(Y))\ot \ct^{b}(X)) & \ct^{a}( \ct^{c}(\ct^{(b^{-1}c)^{-1}g(b^{-1}c)}(Y)))\ot \ct^{a}(\ct^{b}(X)) & \ct^{ac}(\ct^{c^{-1}bgb^{-1}c}(Y))  \ot \ct^{ab}(X)
   % line 5
   \\\ct^{a}(\ct^{c}(Y) \ot \ct^{b}(X)) & \ct^{a}(\ct^{c}(Y))\ot\ct^{a}(\ct^{b}(X)) & \ct^{ac}(Y)  \ot \ct^{ab}(X)\\}; \path[-stealth]
  (m-1-1) edge node [above] {$(\ct^{a})_{2}^{\ct^{b}(X),\ct^{c}(Y)}$} (m-1-2)
   (m-1-2) edge node [above] {$( \ct_{2}^{a,b})_{X}\ot( \ct_{2}^{a,c})_{Y}$} (m-1-3)

    (m-2-1) edge [dashed] node   [above] {$$} (m-2-2)
   (m-2-2) edge [dashed] node [above] {$$} (m-2-3)
    (m-3-1) edge [dashed] node   [above] {$$} (m-3-2)
    (m-3-2) edge [dashed] node   [above] {$$} (m-3-3)
     (m-4-1) edge [dashed] node   [above] {$$} (m-4-2)
    (m-4-2) edge [dashed] node   [above] {$$} (m-4-3)
     (m-5-1) edge  node   [below] {$(\ct^{a})_{2}^{\ct^{c}(Y), \ct^{b}(X)}$} (m-5-2)
    (m-5-2) edge node   [below] {$( \ct_{2}^{a,c})_{Y}
    \ot ( \ct_{2}^{a,b})_{X} $} (m-5-3)
         %columns
               (m-1-1) edge node [left] {$\ct^{a}(c(\ct^{b}(X)\ot \ct^{c}(Y))$} (m-2-1)
    
    (m-1-3) edge node [left] {$c(\ct^{ab}(X), \ct^{ac}(Y))$} (m-2-3)
     (m-2-1) edge node [left] {$\ct^{a}((\ct_{2}^{bgb^{-1},c})_{Y}\ot 1)$} (m-3-1)
    (m-2-2) edge [dashed] node [right] {$$} (m-3-2)
    (m-2-3) edge node [left] {$(\ct^{abgb^{-1}a^{-1}, ac}_{2})_{Y}\ot 1$} (m-3-3)
     (m-3-1) edge node [left] {$\ct^{a}((\ct_{2}^{c, c^{-1}bgb^{-1}c})_{Y}\ot 1)$} (m-4-1)
    (m-3-2) edge [dashed] node [right] {$$} (m-4-2)
    (m-3-3) edge node [left] {$(\ct^{ac, c^{-1}bgb^{-1}c}_{2})^{-1}_{Y} \ot 1$} (m-4-3)
     (m-4-1) edge node [left] {$\ct^{a}(\ct^{c}(c_{Y}^{(b^{-1}c)^{-1}g(b^{-1}c)})\ot 1)$} (m-5-1)
    (m-4-2) edge [dashed] node [right] {$$} (m-5-2)
    (m-4-3) edge node [left] {$\ct^{ac}(c_{Y}^{c^{-1}bgb^{-1}c})\ot 1$} (m-5-3)
   ;
\end{tikzpicture}
  \end{center}
}
\epf %\bll{Include the diagram in the tex. }
\bn{corollary} With the above notations one has that
\beq
(\ct_{2}^{a})^{T^{b}(X), T^{c}(Y)}\ct^{a}(d_{X, Y}^{c,b}d_{X, Y}^{b, c})=(d_{X, Y}^{ac,ab}d_{X, Y}^{ab, ac})(\ct_{2}^{a})^{T^{b}(X), T^{c}(Y)}
\eeq
\end{corollary}
\bpf
Applying Lemma \ref{cdcomp1} for $d_{X, Y}^{b,c}$ and its reverse $d^{c,b}_{Y, X}$ and then composing the two of them one obtains the above equation:
\epf
\br\label{obsd}
Note that the above corollary implies that if $d_{X, Y}^{b,c}d_{Y, X}^{c, b}=\omega\id_{\ct^{b}(X)\ot \ct^{c}(Y)}$ then  $d_{X, Y}^{ab,ac}d_{Y, X}^{ac, ab}=\omega\id_{\ct^{ab}(X)\ot \ct^{ac}(Y)}$. 
\er
\bl \label{segonde}With the above notation one has:
\beq
\ct^{a}(d^{b, c}_{ _{X, Y}})((\ct_{2}^{ca^{-1},a}\ct_{2}^{a,c})_{Y}\ot (\ct_{2}^{ba^{-1},a}\ct_{2}^{a,b})_{X})=((\ct_{2}^{ba^{-1},a}\ct_{2}^{a,b})_{X}\ot (\ct_{2}^{ca^{-1},a}\ct_{2}^{a,c})_{Y})d^{b, c}_{ _{\ct^{a}(X), \ct^{a}(Y)}}
\eeq
\el
\bpf The proof is similar to the proof of Lemma \ref{cdcomp1}.
It follows by a straightforward computation from Equation \eqref{c22}.
\epf
\noindent As before, by composition this implies that 
\beqn
\ct^{a}[d^{b, c}_{ _{X, Y}}d^{c, b}_{ _{Y, X}}][(\ct_{2}^{ba^{-1},a}\ct_{2}^{a,b})_{X}\ot (\ct_{2}^{ca^{-1},a}\ct_{2}^{a,c})_{Y}]=\eeqn \beqn= [(\ct_{2}^{ba^{-1},a}\ct_{2}^{a,b})_{X}\ot (\ct_{2}^{ca^{-1},a}\ct_{2}^{a,c})_{Y}][d^{b, c}_{ _{\ct^{a}(X), \ct^{a}(Y)}}d^{c, b}_{ _{\ct^{a}(Y), \ct^{a}(X)}}]
\eeqn
Thus if $d^{b, c}_{ _{X, Y}}d^{c, b}_{ _{Y, X}}=\omega \id_{\ct^{b}(X)\ot \ct^{c}(Y)}$ for some scalar $\omega$ then

$d^{b, c}_{ _{\ct^{a}(X), \ct^{a}(Y)}}d^{c, b}_{ _{\ct^{a}(Y), \ct^{a}(X)}}=\omega \id_{\ct^{b}(\ct^{a}(X))\ot \ct^{c}(\ct^{a}(Y))}$.

\subsection{M\"{u}ger's centralizer: necessary and sufficient conditions}

Next proposition generalizes \cite[Lemma 5.1]{nnw}. We keep the fixed system $\mathrm C=\cup_{Y}C_{Y}$ of  isomorphisms as above. Recall  that we have chosen the systems of chosen isomorphisms $C_{Y_{i}}$ such that the values of $\al_{Y_{i}}$ are all roots of unity.  It follows that the values of $\al_{[Y]}$ are all roots of unity. \bp\label{centrm} Suppose that $X \in \cc_{g}$ and $Y \in \cc_{h}$ are two simple objects of a braided $G$-crossed fusion category. Then two simple objects $S_{X, \pi}$ and $S_{Y, \delta}$ centralize each other in $\cc^{G}$ if and only if the following conditions hold:
\bne
\item  $K_{g}$ commutes elementwise with $K_{h}$.
\item $K_{g}\subseteq G_{Y}$ and $K_{h}\subseteq G_{X}$.
\item For all $ m,n \in G$ there are some roots of unity $\omega^{m,n}_{X, Y}\in k^{*}$ such that the following two relations are satisfied:
\beq\label{compd}
d_{X, Y}^{m,n}d_{Y, X}^{n,m}=\omega^{m,n}_{X, Y}\id_{\ct^{m}(X)\ot \ct^{n}(Y)}%c_{\ct^{m}(X), \ct^{n}(Y)}(\ct_{2}^{nhn\inv, m})_{X}({\ct_{2}^{m, nhn\inv})_{X}}^{-1}\ot \id_{\ct^{n}(Y)})=\omega_{m, n}\id_{\ct^{n}(X)\ot \ct^{m}(Y)}
\eeq
and	
\beq\label{compi}
\pi(m^{-1}nhn^{-1}m)\ot \delta(n^{-1}mgm^{-1}n)=(\omega^{m,n}_{X, Y})^{-1}\id_{V_{\pi}\ot V_{\delta}}.
\eeq
%\beqn\frac{\pi(nhn^{-1})}{\pi(1)}\ot \frac{\delta(mgm^{-1})}{\delta(1)}=\omega_{m, n}^{-1}\id_{V_{\pi}\ot V_{\delta}}.\eeqn
\ene
\ep
\newcommand{\mcr}{\mathcal{R}}

\bpf 
Note that $G_{X}\subseteq C_{G}(g)\subseteq G$ where $g=\deg(X)$. Fix a set of representatives $\mcr_{g}$ for the left cosets $G/C_{G}(g)$ and $\mcr_{X}$ a set of representatives for the left cosets of $C_{G(g)}/G_{X}$. Recall that 
\beqn 
S_{X, \pi}=\bigoplus_{m \in G/G_{X}}\ct^{m}(X)\ot V_{\pi}=\bigoplus_{r \in \mcr_{g},\;s\in \mcr_{X}}\ct^{rs}(X)\otimes V_{\pi}.
\eeqn 
Similarly 
\beqn 
S_{Y, \delta}= \bigoplus_{n \in G/G_{Y}}\ct^{n}(Y)\ot V_{\delta}=\bigoplus_{a \in \mcr_{h},\;\ b \in \mcr_{Y}}\ct^{ab}(Y)\ot V_{\delta}.
\eeqn 
Then one can apply  Lemma \ref{commuting} with $S=S_{X , \pi}$ and $T=S_{Y, \delta}$. In this situation one has $K_{S}=K_{g}$, $K_{T}=K_{h}$ with
$S_{rgr^{-1}}=\bigoplus_{s \in \mcr_{X}} \ct^{rs}(X)\ot V_{\pi}$ and $T_{aha^{-1}}=\bigoplus_{b \in \mcr_{Y}}\ct^{ab}(Y) \ot V_{\delta}$. 
\mdn
Then Equation \eqref{compo2} for $u=rgr^{-1}$ and $v = aha^{-1}$ can  be written as follows
\beqarn
\bigoplus_{s \in \mcr_{X},b \in \mcr_{Y}}(d_{X, Y}^{rs,ab}d_{Y, X}^{ab,rs}\ot (\pi((ab)g(ab)^{-1})\ot \delta((rs)h(rs)^{-1})))=\\=\bigoplus_{s \in \mcr_{X},b \in \mcr_{Y}}\id_{\ct^{rs}(X) \ot \ct^{ab}(Y)}\ot \id_{V_{\pi}\ot V_{\delta}}.
\eeqarn Let $m:=rs$ and $n:=ab$.
The last identity holds in the $k$-vector space \beqarn \End_{\cc}(\ct^{m}(X)\otimes  \ct^{n}(Y))\ot_{k}(\End_{k}(V_{\pi})\otimes \End_{k}(V_{\delta}))\simeq \\ \simeq \End_{\cc}(\ct^{m}(X)\otimes \ct^{n}(Y))\ot_{k}\End_{k}(V_{\pi}\otimes V_{\delta}).\eeqarn
\mdn Note that one has $\End_{\cc}(X\ot V)\simeq \End_{\cc}(X)\ot \End_{k}(V)$. This implies that there is some nonzero scalar $\omega_{m, n}\in k^{*}$ such that $d_{X, Y}^{m,n}d_{Y, X}^{n,m}=\omega_{m,n}\id_{\ct^{m}(X)\ot \ct^{n}(Y)}$ and $\pi(m^{-1}nhn^{-1}m)\ot \delta(n^{-1}mgm^{-1}n)=\omega_{m, n}^{-1}\id_{V_{\pi}\ot V_{\delta}}$. \mdn
On the other hand, by the choice of the system $C_{X}, C_{Y}$ it follows that if  $\pi(g)=\al\id_{V_{\pi}}$ then $\al$ must be a root of unity.  This implies that all $\omega_{m, n}$ are roots of unity.
\epf 
%\bll{Here I have to explain that omega does not depend on the isomorphism class.}
\begin{corollary} \label{cp}With the above notations the objects $S_{X,\pi}$ and $S_{Y, \delta}$ centralize each other if and only if for any $p\in G$ there is a root of unity $\omega_{p}(X, Y)\in k^{*}$  such that the following two relations hold:
\beq\label{cp1}
d^{1,\;p}_{X,Y}d^{p,\;1}_{Y, X}=\omega_{p}(X, Y)\id_{X\ot \ct^{p}(Y)}
\eeq
\beq\label{reform}
\frac{\pi(php\inv)}{\pi(1)}\ot \frac{\delta(p\inv gp)}{\delta(1)}=\omega_{p}(X, Y)^{-1}\id_{V_{\pi}\ot V_{\delta}}
\eeq
\end{corollary}
\bpf
Remark \ref{obsd} implies that if $X\perp_{G}Y$ then $\omega^{m,n}_{X, Y}=\omega^{1, m^{-1}n}_{X, Y}=\omega^{n^{-1}m,1}_{X, Y}$. The corollary follows by denoting $\omega_{p}(X, Y)
:=\omega^{1, \;p}_{X, Y}$ for any $p\in G$.c%\bll{this is correct}.
\epf
\br
Suppose that $G$ is an abelian group.  Let $X \in \cc_{g}$ and $Y \in \cc_{h}$ and suppose that the two simple objects $S_{X , \pi}$ and $S_{Y, \delta}$ centralize each other in $\cc^{G}$. Then for any $m\in G$ the two identities above can be written as follows:
\beq\label{omegas}
d_{X, Y}^{1,m}d_{Y, X}^{m,1}=\omega\id_{X\ot \ct^{m}(Y)}
\eeq
and
\beq
\pi(h) \ot \delta(g)=\omega^{-1} \id_{V_{\pi}\ot V_{\delta}}
\eeq
Thus in this situation $S_{X, \pi}$ centralizes $S_{Y, \pi}$ if and only if $\deg(X)\in G_{Y}$, $\deg(Y) \in G_{X}$ and there is $\omega \in k^{*}$ root of unity such that the above two identities are satisfied.
\er
\subsection{Definition of the $G$-centralizer}\label{gcross}
Let $\cc$ be a braided $G$-crossed fusion category and two simple objects $X \in \cc_{g}$ and $Y \in \cc_{h} $. We say $X$ and $Y$ {\it $G$-centralize} each other  (and we write $X \perp_{G}Y$)  if there are projective representations $\pi \in \irr(k_{\al_{X}}[G_{X}])$ and $\delta \in \irr(k_{\al_{Y}}[G_{Y}])$ such that the simple objects $S_{X, \pi}$ and $S_{Y, \delta}$ centralize  in $\cc^{G}$. Equation \eqref{isochocomp} shows that if $X\perp_{G} Y$ and $X \simeq Y_{i}$, $Y \simeq Y_{j}$ then $Y_{i} \perp_{G} Y_{j}$. Moreover by the same relation we have that
\beq
\omega^{m,n}_{X, Y}=\omega^{m,n}_{Y_{i}, Y_{j}}
\eeq
for all $m,n \in G$. {\it For this reason we denote the above common value by $\omega_{[X],[Y]}^{m,n}$. Thus there is also a well defined quantity $\omega_{p}([X],[Y]):=\omega^{1,p}_{X, Y}$.} For the rest of this paper we denote $\omega([X],[Y]):=\omega_{1}([X],[Y])$.
  
\mdn
By Corollary \ref{cp} it follows that $X\perp_{G}Y$ if and only if there are roots of unity $\omega_{p}([X], [Y])\in k^{*}
$ and projective representations $\pi$, $\delta$ such that the above two Equations \eqref{cp1} and \eqref{reform} are satisfied. 
\mdn
Note also that by its definition $\perp_{G}$ is a symmetric relation, i.e. $X \perp_{G} Y$ if and only if $Y \perp_{G} X$. Moreover if $X \perp_{G}Y$ then $\ct^{m}(X)\perp_{G} \ct^{n}(Y)$ for all $m,n \in G$.\mdn 

\br 
Note that if $X, Y \in \cc_{1}$ then 
\beqn
d_{X,Y}^{m, n}=c_{\ct^{m}(X),\; \ct^{n}(Y)}
\eeqn
coincides to the braiding of $\cc_{1}$ and does not depend  on the chosen system of isomorphisms. Note that $d_{X,Y}^{m, n}$ depends only on $c_{Y}^{1}$ and $c_{X}^{1}$ which we supposed to be the identity.
 It follows that $X \perp_{G} Y$ if and only if $\ct^{m}(X)$ projectively centralize $\ct^{n}(Y)$ in $\cc_{1}$ and there are projective representations $\pi$ and $\delta$ such that for any $p \in G$ and Equations \eqref{cp1} and \eqref{reform}  are satisfied.
\er
%\subsubsection{Extend the function $d$ on the whole grothendieck ring.}
\subsection{A formula for M\"{u}ger centralizer}\label{fmc}
%\bll{recall here the definition of the fusion datum}\mdn
%\bll{Note that if $X$ and $Y$ are from $\cc_{1}$  then the $d_{X,Y}=c_{X,Y}$ and therefore it coincides with the M\"{u}ger centralizer.}
Recall that for a fusion datum $(H, \cs, \lam)$ as in Section \ref{intro}
we denote by $\cc(H, K, \cs, \lam) $ the  corresponding fusion subcategory $\cc(H, \cs, \lam) $ where $K$ is the support of $\cs$.
\mdn
{\bf Proof of Theorem \ref{form}:}  Since $H$ leaves invariant the objects of $\cs$ it follows that $H$ commutes elementwise with $K$ which proves i). Indeed, if $a \in K$, $b \in H$ and $X_{a} \in \cc_{a} \cap \cs$ then $X_{a}\cong \ct^{b}(X_{a}) \in \cc_{bab^{-1}}$ which shows that $ab=ba$. Clearly $\mathcal T \perp_{G} \cs$ by definition. The rest of the second item follows directly from  Theorem \ref{switch}. 
\mdn
For the next item we use Proposition \ref{centrm}. Suppose that $S_{X_{a}, \pi} \in \cc(H,K,\cs, \lam)$ and $S_{Y_{b}, \delta}\in \cc(K, H, \mathcal{ T}, \tilde \lam)$. Then by Equation \eqref{defi} one has  $\lam([X_{a}],b)=\frac{\ch_{\pi}(b)}{\ch_{\pi}(1)}$ and $\tilde \lam([Y_{b}],a)=\frac{\ch_{\delta}(a)}{\ch_{\delta}(1)}$. \mdn Since $S_{X_{a}, \pi}$ centralize $S_{Y_{b}, \delta}$ it  follows by Equation \eqref{reform} that for any $p \in G$ one has
\beqarn
\tilde \lam ([Y_{b}], pap^{-1})&=&\frac{\ch_{\delta}(p^{-1}ap)}{\ch_{\delta}(1)}=(
\frac{\ch_{\pi}(pbp\inv)}{\ch_{\pi}(1)})^{-1}\omega_{p}^{-1}([X_{a}], [Y_{b}])\\ &=& \lam^{-1}([X_{a}], pbp^{-1})\omega_{p}^{-1}([X_{a}], [Y_{b}])\eeqarn%For $m=n=1$ one obtains Equation \eqref{tilde}.
For $p=1$ one obtains the desired equality. \qed
%\bll{It seems that these omega's might depend only on the homogenous degrees of $X$ and $Y$, i.e., $a$ and $b$.}
\br Suppose that $S_{X , \pi} \in   \cc(H, K, \cs, \lam)$ with $ X \in \cc_{a}$. Since $ \lam([X],1)=1$ it follows from Equation \eqref{tilde} of the previous theorem that 
\beq\label{f1tilde}
\tilde \lam([Y], a)=\omega([X], [Y])
\eeq
for any $Y\in \cc_{1}$ such that $S_{Y, \pi}\in \cc(H, K, \cs, \lam)'$ and any $a \in K$.
\er
%\bp One has that $\cc(H,K,\lam, \cs, \lam)'=\cc(K, H_{\cd''}, \mathcal T, \tilde \lam)$.\ep
\br\label{q1} Note that $\mathcal T$ is contained in the largest fusion subcategory of $\cc$ whose objects $G$-centralize any object of $\mathcal{S}$. \er
\section{Braided $G$-crossed fusion categories from crossed $G$-modules}\label{gcp}

In this section we apply our results to equivariantizations of crossed pointed fusion categories obtaining a formula for the centralizer in any braided group theoretical category. 
\mdn Recall that a  fusion category $\C$  is \emph{pointed} if all simple objects of $\C$ are invertible. Then there is an equivalence of fusion categories $\C \simeq \C(X, \omega)$, where $X $ is the group of isomorphism classes of invertible objects in $\C$ and   $\omega: X \times X \times X \to k^*$ is an invertible normalized 3-cocycle. Moreover, $\C(X, \omega) = \vect_{X}^{\omega}$ is the category of finite dimensional $X$-graded vector spaces with associativity constraint induced by $\omega$.

\subsection{Group actions on $\C(X, \omega)$ and equivariantizations}(see \cite{buna})
\noindent
Let $\C = \C(X, \omega)$ and let $G$ be a finite group. An action $\ct: \underline G \to \underline \Aut_{\otimes} \C$ of $G$ on $\C$ is determined by an action by group automorphisms of $G$ on $X$,  that we shall indicate by $x \mapsto {}^gx$, $x \in X$, $g \in G$, and two maps $\tau: G \times X \times X \to k^*$, and $\sigma: G \times G \times X \to k^*$, satisfying
\begin{align*}\frac{\omega(x, y, z)}{\omega({}^gx, {}^gy, {}^gz)} & = \frac{\tau(g; xy, z) \, \tau(g; x, y)}{\tau(g; y, z) \, \tau(g; x, yz)} \\
1 & = \frac{\sigma(h, l; x) \, \sigma(g, hl; x)}{\sigma(gh, l; x) \, \sigma(g, h; {}^lx)} \\
\frac{\tau(gh; x, y)}{\tau(g; {}^hx, {}^hy) \, \tau(h; x, y)} & = \frac{\sigma(g, h; x) \, \sigma(g, h; y)}{\sigma(g, h; xy)},
\end{align*} for all $x, y, z \in X$, $g, h, l \in G$.
\mdn
We shall also assume that $\tau$ and $\sigma$ satisfy the additional normalization conditions $\tau(g; x, y) = \sigma(g, h; x) = 1$, whenever some of the arguments $g$, $h$, $x$ or $y$ is an identity.

\mdn The action $\ct: \underline G \to \underline \Aut_{\otimes} \C$ determined by this data is defined by letting $\ct^g(x) = {}^gx$, for all $g\in G$, $x\in X$, and $\ct^g = \id$ on arrows, together with the following constraints:
\begin{equation}\label{formulas-pointed}(\ct_2^{g, h})_x = \sigma(g, h; x)^{-1} \id_{^{gh}x}, \quad (\ct^g)_2^{x, y} = \tau(g; x, y)^{-1}\id_{xy}, \quad \ct^g_0 = \id_e,
\end{equation} for all $g, h \in G$, $x, y \in X$. See also \cite[Section 7]{tyinv}.
\mdn
Recall that a fusion category is called \emph{group-theoretical}
if it is Morita equivalent to a pointed fusion category. By \cite[Theorem 7.2]{nnw}, a \emph{braided}
fusion category is group-theoretical if and only  if it is an
equivariantization of a pointed fusion category. Moreover, it  was shown in \cite[Theorem 5.3]{Na}  that
every braided group-theoretical fusion category is equivalent to
an equivariantization $\C(\xi)^G$ of a crossed pointed fusion
category $\C(\xi)$ associated to a quasi-abelian 3-cocycle $\xi$
on a finite crossed module $(G, X, \partial)$, under a canonical
action of $G$ on $\C(\xi)$.
\mdn Recall that a finite crossed module is a triple $(G, X, \partial)$ consisting of a finite group $G$ acting by automorphisms on a finite group $X$, and a group homomorphism $\partial: X \to G$ such that $${}^{\partial(x)}y = xyx^{-1}, \quad \partial({}^gx) = g\partial(x)g^{-1}, \qquad g \in G,\; x, y \in X,$$ where $x \mapsto ^gx$, $x \in X$, $g \in G$, denotes the action of $g$ on $X$.
\mdn As a fusion category $\C(\xi) = \C(X, \omega)$, and the action of $G$ on $\C(\xi)$ is determined by the action of $G$ on $X$ and formulas \eqref{formulas-pointed}, with respect to $\sigma(g, h;x) : = \gamma(g, h; x)^{-1}$, $\tau( g\;x, y) : = \mu(g; x, y)^{-1}$, $x, y \in X$, $g, h \in G$. See \cite[Subsection 4.1]{Na}. Following loc. cit. we will use the notation $\gm_{g, h}(x):=\gm(g,h;x)$.
\mdn
Recall that a quasi-abelian 3-cocycle $\xi$ on $(G, X, \partial)$ is a quadruple $\xi = (\omega, \gamma, \mu, c)$, where $\omega: X \times X \times X \to k^*$ is a $3$-cocycle, $\gamma: G \times G \times X \to k^*$, $\mu: G \times X \times X \to k^*$ and $c: X \times X \to k^*$ are maps satisfying the compatibility conditions in \cite[Definition 3.4]{Na}.
%\subsubsection{The category $\cc(\omega, \gm,\mu, c)$}
%To a quasi abelian $3$-cocycle it is associated in \cite{Na} a braided $G$-crossed pointed fusion category as following:
%\subsubsection{The action of $G$ on $\cc(\omega, \gm,\mu, c)$} Let $\cc={\bvec^{\omg}}_X$.
%}With the above notations $G$ acts on $\bvec^{\omg}_X$ given by the action $\ct^{g}(k_x)=k_{\;^gx}$ and the following data:
%\beq
%(\ct^{g,h}_{2})_{k_{x}}:\ct^g(\ct^h(k_x))\ra \ct^{gh}(k_x),\;\; v \mapsto \gamma_{g,h}(x)v\eeq
%and \beq (\ct^{g}_{2})^{k_{x}, k_{y}}:\ct^g(k_x\ot k_y)\ra %\ct^g(k_x)\ot \ct^g(k_y),\;\;v\ot w\mapsto \mu_g(x,y)(v \ot w).\eeq
\subsection{Simple objects of the equivariantized category $\cc(\omega)^G$}
%\subsubsection{Chosen isomorphisms}
Fix the system $C\cong \mtr{1}$ of chosen isomorphisms given by $c^g_x:\;k_{^gx}\ra k_{x}=\id_{k_{x}}$ as identity, for any $g \in G_x$, the stabilizer group of $x$. Fix also a set $\Gm$ of representatives for the orbits of the action.
\mdn Using Equation \eqref{alfa} it follows for any $a \in \Gm$ one has that
\bq\lb{cna}
\al_x(g,h)=\gamma_{g,h}(x)
\eq
for any $x \in X$ and any $g,h \in G$.\mdn
%\subsubsection{The difference $d_a(t,h)$}
%\bb{
Following Equation \eqref{dyt} it follows that
\bq
D_{g, y}(h)=\frac{\gm_{ghg^{-1},g}(y)}{\gm_{g,h}(y)}
\eq
\mdn 
Note that the projective representation $\tau^{y,z}_{yz}$ is one dimensional given by 
\beqn
g.1_{k}=\tau(g;y,z)^{-1}1_{k}.
\eeqn
%\blue{There is no scalar involved!!}
%\subsubsection{On the simple modules $S_{x,\pi}$ of $\cc^{G}$} 
\mdn Let $\pi$ be a $\al_{x}$-projective representation of the stabilizer $G_{x}$. By definition $S_{x, \pi}$ is given by:
\bq
S_{x, \pi}=\oplus_{t\in G/G_x}\ct^t(x)\ot V_{\pi}
\eq
The equivariant structure of $S_{x,\pi}$ is given on the components by
\bq\lb{eqvr}
\ct^g(\ct^t(a))\ot v\xra{\gm_{g, t}(x)\ot 1}\ct^{gt}(x)\ot v \xra{\gm_{t', h}^{-1}(x)\ot \pi(h)} \ct^{t'}(x)\ot hv
\eq
where $t'$ is chosen such that $gt=t'h$ with $h\in G_{x}
$. It can be easily seen that this coincides to the simple module  denoted by $(a, \pi)$ in \cite{Na}.
\subsection{$G$ grading on $\cc$}
For each $g \in G$, let $\cc_{g}$ denote the full abelian
subcategory consisting of objects of $\vvec_{X}$ supported on $\delta^{-1}(g) \subset X$, i.e., objects of $\cc_{g}$
are defined to be finite-dimensional $\delta^{-1}(g) $-graded vector spaces.\mdn
%\subsection{The cross $G$-braiding in $\cc(\omega, \gm,\mu, c)$ }
Recall that the G-braiding $\tilde c$ is defined by
\beq
c(x, y)\id_{k_{(\;^{g}x)y}} =k_{x}\ot k_{y}\ra k_{\;^{g}x}\ot k_{y},
\eeq
where $g=\delta(x)$.
\subsection{Fusion subcategories of $\cc^{G}$}
By Definition \ref{fusdatum} all fusion subcategories of $\cc^{G}$ are parameterized by triples $(H, Y, \lam)$ where $H$ is a normal subgroup of $G$, $Y$ a subgroup of $X$ which is stable under the action of $G$ and fixed by $H$.  Moreover  $\lam: H\times Y\ra k^{*}$ is a twisted bicharacter satisfying the following relations:
\beqn
\lam(y, hh')=\lam(y,h)\lam(y',h)\gm_{h, h'}(y)
\eeqn
\beqn
\lam(yz, h)=\lam(y,h)\lam(z, h)\tau(h;y,z)
\eeqn
\beqn
\lam(\;^{g}y, h)=\lam(y,g^{-1}hg)D_{g,y}(h)
\eeqn
for any $y,z\in Y$, any $g \in G$ and any $h \in H$.  
\subsection{The double distance in $X$ and the $G$-centralizer.}
Suppose that we have two objects $k_{x}, k_{y}$ for which their conjugacy classes in $G$ of their degrees commute element by element. Then the distance from subsection \ref{dd} is given by  
\beq
d^{m,n}_{x,y}:k_{\;^{m}x}\otimes k_{\;^{n}}y\ra k_{\;^{n}}y\otimes k_{\;^{m}x}
\eeq
is given by multiplication with the scalar 
\beq
\gm_{mgm^{-1},n}(y)\gm^{-1}_{n,n^{-1}mgm^{-1}n}(y)c_{\;^{m}x, \;^{n}y}c_{\;^{m}x,\;^{n}y}
\eeq
where $g=\deg(x)$. 
In particular $d^{1,n}_{x, y}=\gm_{g,n}(y)\gm^{-1}(n, n^{-1}gn)$
\mdn
Thus for any two elements $x, y\in X$ such that the conjugacy classes in $G$ of $\delta(x)$ and $\delta(y)$ commute each other the elements $k_{x}$ and $k_{y}$ of $\cc$ $G$ centralize each other.
\mdn
Thus by Equation \eqref{omegas} one has that 
\beq
\omega(x,y)=d^{1,n}_{x, y}d^{n,1}_{y, x}=
\gm_{g,n}(y)\gm^{-1}(n, n^{-1}gn)c_{x,\;^{n}y}\gm_{nhn^{-1},1}(x)\gm^{-1}_{1,h}(x)c_{\;^{n}y, \;x}
\eeq
where $h=\deg(y)$.
\subsection{M\"uger centralizer formula for equivariantizations of cross pointed fusion categories}
Applying Theorem \ref{form} one obtains the following:
\begin{theorem}
Suppose that $\delta$ is surjective and $(\omega, \gamma, \mu, c)$ is nondegenerate in the sense
of \cite[Definition 3.10]{Na}. Then for any fusion subcategory one has the following
\beq
\cc(H,Y,\lam)'=\cc(\delta(Y), T, \tilde\lam)
\eeq
where $T $ is a subgroup of $X$ commuting with $Y$ and  with $ \delta(T)=H$. Moreover in this case one has that: \beqn\tilde\lam(\delta(y), t)=\lam^{-1}(h, y)\omega(y,t)\eeqn
for any $h\in H$ and $y \in Y$ and $t \in T$ with $\delta(t)=h$.
\end{theorem}\bpf Following \cite[Proposition 5.6]{Na} we deduce that $\cc(\xi)^{G}$ is nondegenerate and we can apply Theorem \ref{form}. 
\epf\mdn On the other hand from the non-degeneracy condition it follows by Equation \eqref{fpdc} that
\beq
\fp(\cc(H,Y, \lam)')=\frac{\fp(\cc^{G})}{\fp(\cc(H,Y, \lam))}=\frac{|G||X|}{|G||Y|/|H|}=\frac{|X||H|}{|Y|}
\eeq
\mdn
One has that $\cc(H,Y,\lam)'=\cc(\delta(Y), \mathcal{T}, \lam')$ where the support $\delta(\mathcal{T})$ of $\mathcal{T}$ coincides to $H$. Thus $\mathcal{T}\subseteq \delta^{-1}(H)$.
\mdn On the other hand using Equation \eqref{fpdd} one has that
\beq
\fp(\cc(\delta(Y), \mathcal{T}, \tilde\lam))=\frac{|G||T|}{|\delta(Y)|}
\eeq

Note that $|\delta(Y)|=\frac{|Y|}{|\ker(\delta)\cap Y|}$ and $|\ker (\delta)||G|=|X|$. Then the last two equation give that
\beq
|T|=|H|\frac{|\ker(\delta)|}{|\ker(\delta) \cap Y|}.
\eeq
Thus choosing a subgroup $Y$ such that $Y\cap \ker \delta\neq 1$ it follows that $\mathcal T  $ is not the largest subcategory mentioned above in Remark \ref{q1}.\mdn 
As explained in the last paragraph of \cite{Na} the twisted Drinfeld double $D^{\omega}(G)$ of a finite group can be recovered as an equivariantization of a crossed pointed fusion category. Thus the last theorem generalizes the results obtained in \cite[Theorem 5.11]{nnw}.
\section{Appendix}\noindent
In this Appendix we give a brief account of the results on projective representations used throughout this paper.
See for instance \cite{karpilovsky}. Most of the facts of this Appendix are also contained in the Appendix of \cite{buna}.
\mdn Let $G$ be a finite group and let
$\alpha: G \times G \to k^*$ be a (normalized)
2-cocycle on $G$, that is,
$$\alpha(g, h)\alpha(gh, t) = \alpha(g, ht) \alpha(h, t), \quad \alpha(g, e) = 1 = \alpha(e, g),\qquad \forall g, h, t \in G.$$
%for all $g, h, t \in G$.

\mdn A \emph{$\al$-projective representation} $\pi$ of $G$ with
\emph{associated cocycle} $\alpha$
 on a vector space $V$ is a map $\pi: G\ra
\GL(V)$, such that $$\pi(e) = \id_V, \quad
\pi(gh)=\al(g,\;h)\pi(g)\pi(h),\qquad \forall g, h\in
G.$$
%for all $g, h \in G$.
In other words, $\pi$ is a representation of the twisted group
algebra $k_{\alpha} G$ on the vector space $V$. We shall
also use the notation  $V_{\pi} = V$ to indicate such a projective
representation.
\mdn Two projective representations $\pi$ and $\pi'$ of $G$
are called \emph{(projectively) equivalent}  if there is a linear
isomorphism $\phi: V_\pi \to V_{\pi'}$ and a map $f: G \to k^*$
such that $\phi \pi(g) = f(g)\pi'(g) \phi$, for all $g \in G$. In
this case we shall use the notation $\pi' \simeq \pi$.

\noindent
If $\pi' \simeq \pi$, then the associated cocycles
$\alpha$ and $\alpha'$ are related by
$$\alpha(g, h) = \alpha'(g, h)f(g)f(h)f(gh)^{-1}, \quad g, h \in
G,$$  that is, $ \alpha$ and $ \alpha'$ are
cohomologus cocycles, and thus they belong to the same cohomology class $\alpha \in H^2(G, k^*)$. %We shall also call $\pi$ an
%$\alpha$-projective representation. 
Note that the map $f: G \to
k^*$ induces an algebra isomorphism $\tilde f:
k_{\alpha} G \to k_{\alpha'} G$ in the form
$\tilde f(g) = f(g)g$, for all $g \in G$. Thus $\pi$ and $\pi'$
are equivalent projective representations if and only if $V_{\pi}
\simeq {\tilde f}^*(V_{\pi'})$ as $k_{\widetilde\alpha}G$-modules.
\br\label{unity}(see \cite{karpilovsky})
Given a two cocycle $\alpha: G \times G \to k^*$ one can always find a cohomologus cocycle $\al'$ such that the values $\al'(g,g')$ are all roots of unity for all $g,g' \in g$.
\er
\mdn Let $\pi$ and $\pi'$ be projective representations of $G$ with associated cocycles $\alpha$ and $\alpha'$,
respectively. The \emph{tensor product} $\pi \otimes \pi'$ is the $\alpha \al'$-projective representation of $G$ on
the vector space $V_{\pi} \otimes V_{\pi'}$ defined by $(\pi
\otimes \pi')(g) (u \otimes v) = \pi(g)u \otimes \pi'(g)v$. In
particular, if $\pi$ is a representation of $G$, then $\pi \otimes
\pi'$ is again a $\al'$-projective representation with associated cocycle $\al'$.
$\al$.
\mdn If $\pi_1$ and $\pi_1'$ are projective representations
projectively equivalent to $\pi$ and $\pi'$, respectively, then
the tensor products $\pi_1\otimes \pi'_1$ and $\pi\otimes \pi'$
are projectively equivalent. Further, suppose that $\pi'$ is a
one-dimensional representation, that is, a linear character of
$G$. Then  $\pi$ and  $\pi \otimes \pi'$ are projectively
equivalent via the canonical isomorphism $\phi: V_{\pi} \to
V_{\pi} \otimes k$, $v \mapsto v \otimes 1$, and the map $f: G \to
k^*$ given by $f(g) = \pi'(g)^{-1}$, for all $g \in G$.
\mdn A nonzero projective representation $\pi: G \to \GL(V)$
of $G$ is called \emph{irreducible} if $0$ and $V$ are the only
subspaces of $V$ which are invariant under $\pi(g)$, for all $g
\in G$. Hence, $\pi$ is irreducible if and only if it is not
projectively equivalent to a projective representation
$\rho$ of the form $$\rho(g) = \left(\begin{matrix}\pi_1(g) & * \\
0 & \pi_2(g)
\end{matrix}\right), \quad g \in G,$$ where $\pi_1$ and $\pi_2$
are nonzero projective representations or, equivalently, if $V$ is
a simple $k_{\widetilde\alpha} G$-module, where $\widetilde
\alpha$ is the factor set of $\pi$ \cite[Theorem
3.2.5]{karpilovsky}.

\mdn Let $\pi: G \to \GL(V)$ be a projective representation
of $G$ with factor set $\widetilde\alpha$. Since the group algebra
$k_{\widetilde\alpha} G$ is semisimple, then $V = V_\pi$ is
\emph{completely reducible}, that is, $V_\pi \simeq V_{\pi_1}
\oplus \dots \oplus V_{\pi_n}$, where $V_{\pi_i}$ is a simple
$k_{\widetilde\alpha} G$-module, for all $i = 1, \dots, n$. If
$\pi'$ is an irreducible projective representation with factor set
$\widetilde\alpha'$, then $\pi'$ is called a \emph{constituent} of
$\pi$ if $\pi'$ is projectively equivalent to $\pi_i$ for some
$1\leq i \leq n$. In this case, the \emph{multiplicity} (or
\emph{intertwining number}) of $\pi'$ in $\pi$ is defined as
$$m_G(\pi', \pi) : = \dim \Hom_{k_{\widetilde\alpha} G}(V_{\pi_i}, V_{\pi}).$$
Observe that if $\pi'$ is a constituent of $\pi$, then the
cocycles $\widetilde\alpha'$ and  $\widetilde\alpha$ belong to the
same class  in $H^2(G, k^*)$. Letting $\widetilde\alpha' df =
\widetilde\alpha$, with $f: G \to k^*$, we have that $m_G(\pi',
\pi) : = \dim \Hom_{k_{\widetilde\alpha} G}({\tilde
f}^*(V_{\pi'}), V_{\pi})$, where $\tilde f: k_{\widetilde\alpha} G
\to k_{\widetilde\alpha'} G$ is the isomorphism associated  to
$f$.
\mdn The \emph{character} of a projective representation
$\pi: G \to \GL(V)$ is defined as the map $\chi = \chi_V: G \to
k^*$ given by $\chi(g) = \Tr(\pi(g))$, for all $g \in G$. Let
$\widetilde\alpha$ be the factor set of $\pi$. If $\pi'$ is an
irreducible projective representation of $G$ with factor set
$\widetilde\alpha$ and character $\chi'$, then the multiplicity of
$\pi'$ in $\pi$ can be computed by the formula
\begin{align*}m_G(\pi', \pi)  = \langle \chi', \chi \rangle : & = \frac{1}{|G|} \sum_{g \in G}  \frac{1}{\alpha(g^{-1}, g)} \chi'(g) \chi(g^{-1})  \\
& = \frac{1}{|G|} \sum_{g \in G^0}
\frac{1}{\alpha(g^{-1}, g)} \chi'(g) \chi(g^{-1}),
\end{align*} $G^0 \subseteq G$ is the subset of $\widetilde\alpha$-regular elements of $G$. See \cite[Chapter 5]{karpilovsky}.
\mdn Let $\alpha: G \times G \to k^*$ be a
2-cocycle and let $H \subseteq G$ be a subgroup. Consider a
projective representation $W$ of $H$ with 2-cocycle 
$\alpha\vert_H$. The \emph{induced projective
representation} of $G$ is defined as $\Ind_H^GW =
k_{\alpha}G \otimes_{k_{\alpha} H}W$. This is
a projective representation of $G$ with factor set
$\alpha$. By Frobenius reciprocity, we have  natural
isomorphisms
\begin{equation}\label{frob} \Hom_{k_{\alpha} G}(\Ind_H^GW, V) \simeq \Hom_{k_{\alpha} H}(W, V\vert_H),\end{equation}
for every projective representation $V$ of $G$ with factor set
$\alpha$, where $V\vert_H$ denotes the restricted
projective representation of $H$.
\subsection{On the center of a projective representation}
Let $M$ be a $\al$-projective $G$-representation with associated cocycle $\al$. Define \beq
Z_{\al}(M):=\{g \in G\;|\;gm=\psi(g)m \; \text{for some} \; \psi(g)\in k^{*}\; \text{ for all } \; m \in M\}.\eeq
\noindent Clearly $Z_{\al}(M)$ is a subgroup of $G$ and $\psi_M:=\psi$  defines a one dimensional representation of $k_{\al\vert_{Z_{\al}(M)}}[Z_{\al}(M)]$.
\bl\label{a1}
Let $\pi$ be a $\al$-projective representation of $G$ with  $\al(g,g')$ all roots of unity. Then for $h \in G$ one has that $|\ch_{\pi}(h)|\leq \ch_{\pi}(1)$ and one has equality if and only if $h \in Z_{\al}(\pi)$.
\el
\bpf
Since all the values $\al(h^{i}, h^{j})$ are all roots of unit it follows that there is $n>0$ such that $h^{n}=1$ in $k_{\al}[H]$. Therefore $\ch_{\pi}(h)$ is a sum of roots of unity and the statement follows.
\epf
\bl\label{a2} Let $H\subset T \subset G $ be a tower of finite groups and $\pi$ be $\al$-projective representation of $G$.  If $H\subset Z_{ \al\vert_{H}}(\pi\uw^{G}_{T})$ then $H\subset Z_{\al}(\pi)$.
\el
\bpf
Using Equation \eqref{frob} it can be easily checked that $\pi$ is a constituent of $\pi\uw^{G}_{T}\dw^{G}_{T}$.
\epf
\mdn
It is straightforward to check the following lemma.
\bn{lemma}\label{krnind} Let $M$ be a subgroup of $G$ and $\gm$ be the
character of a representation of $M$. Then $\krn_G(\gm\uw^G_{M})=\mtr{core}_G(\krn_M(\gamma)).$
\end{lemma}
\section{Acknowledgments} The author wishes to thank Sonia Natale for the fruitful conversations they had  during the stay at Erwin Schr\"oedinger Institute. Especially the author is deeply indebted to Sonia Natale for fixing an error in the proof of the key Proposition \ref{centrm} of on an earlier version of this manuscript.
\bibliographystyle{amsalpha}

\ed
\bl\lb{scalar}
Let $G$ be a finite group acting transitively on a finite set $X$. Suppose that $\al$ is a $2$-cocycle of $G$ with values in $k^X$ and $M=\oplus_{x \in X}M_x$ is a twisted representation of $G$ with cocycle $\al$.  Then  for any fixed $y \in X$ the vector space $M_y$ is a twisted representation of {\bf $T:=\mtr{St}_G(y)$} with cocycle $\al|_T$.

Then
\bne
\item
$Z_{\al}(M)\subset \cap_{x \in X}\mtr{St}_G(x)$,
\item $Z_{\al}(M)\cap T=Z_{\al|_T}(M_y)$ and $\psi_M|_{Z_{\al|_T}(M_y)}=\psi_{M_y}$.
\ene
Moreover, for any $z \in Z_{\al|_T}(M_y)$ and all $g \in G$ one has that
$$\psi_{M_y}(g^{-1}zg)=\al(z,g)\al(g, g^{-1}zg)^{-1}\psi_{M_y}(z).$$
\bpf

Note that if $z \in Z_{\al|_T}(M_y)$ then $$\psi_M(z)(\bar{1}\ot _{k_{\al}[T]} m)=z(\bar{1}\ot _{k_{\al}[T]} m)=\bar{1}\ot  _{k_{\al}[T]}zm=\psi_{M_y}(z)(\bar{1}\ot _{k_{\al}[T]} m)$$
which shows that $\psi_M(z)=\psi_{M_y}(z)$.
On the other hand for all $m \in M_y$ one has that
\begin{eqnarray*}
\psi_{M}(z)(\overline{g}\ot _{k_{\al}[T]} m) & = &\overline{z}(\overline{g}\ot_{k_{\al}[T]}  m)
\\ & = & \al(z,g)(\overline{zg}\ot_{k_{\al}[T]} m)\\ & = & \al(z,g)(\overline{g(g^{-1}zg)}\ot_{k_{\al}[T]} m)
\\ & = & \al(z,g)\al^{-1}(g,g^{-1}zg) (\overline{g}\overline{g^{-1}zg}\ot_{k_{\al}[T]} m)
\\ & = & \al(z,g)\al^{-1}(g,g^{-1}zg) (\overline{g}\ot_{k_{\al}[T]} ( \overline{g^{-1}zg})m)
\\ & = & \al(z,g)\al^{-1}(g,g^{-1}zg)\psi_{M_y}(g^{-1}zg)(\overline{g}\ot _{k_{\al}[T]} m)
\end{eqnarray*}
which shows that $$\psi_{M_y}(g^{-1}zg)=\al(z,g)\al(g, g^{-1}zg)^{-1}\psi_{M_y}(z).$$
\epf

\bn{defn}\lb{stable}
Let $H$ be a normal subgroup of $G$ and $\al$ be a $2$-cocylcle of $G$. A linear character $\psi$ of $k_{\al|H}[H]$ is called $G$-stable with respect to $\al$ if it verifies the above property
\beqn
\psi(g^{-1}hg)=\al(h,g)\al(g, g^{-1}hg)^{-1}\psi(h)
\eeqn
for all $g \in G$ and all $h \in H$.
\end{defn}
%Moreover
%$M \cong k_{\al}G\ot_{k_{\al}H}M_y$ as projective representations of $G$ with cocycle $\al$.
\el

%textcolor[rgb]{1.00,0.00,0.00}{
\subsection{A cocycle identity}
Suppose that $\al \in H^2(G,k^*)$ is a two cocycle on a finite group $G$. Thus $\al$ verifies the identity
\beq\lb{cocy}
\al(a,b)\al(ab,c)=\al(b,c)\al(a, bc)
\eeq
for all $a,b,c \in G$.
\bl\lb{idnt}
Suppose that $\al \in H^2(G,k^*)$ is a two cocycle on $G$. Then
$$
c_{\al}(x,h)=\frac{\al(x,h)\al(xh, x^{-1})} {\al(x, x^{-1})}=\frac{ \al(h, x^{-1})}{\al(x^{-1},xhx^{-1})}
$$
\el
\bpf
Applying formula \ref{cocy} for $a=x$, $b=h$ and $c=x^{-1}$ one has that $$\al(x,h)\al(xh, x^{-1})=\al(h, x^{-1})\al(x,hx^{-1})$$
On the other hand applying the same formula for $a=x^{-1}$, $b=x$ and $c=hx^{-1}$ one has that
$$
\al(x^{-1},x)=\al(x,hx^{-1})\al(x^{-1}, xhx^{-1})
$$
Thus
$$
\frac{ \al(h, x^{-1})}{\al(x^{-1},xhx^{-1})}=\frac{ \al(h, x^{-1})\al(x,hx^{-1})}{\al(x^{-1},x)}=\frac{\al(x,h)\al(xh, x^{-1})}{\al(x^{-1},x)}
$$
Note also that for $a=c=x$ and $b=x^{-1}$ formula \ref{cocy} gives that $\al(x,x^{-1})=\al(x^{-1},x)$ for all $x \in G$.
\epf

\bll{\br\noindent
Note that $\lam_{\cd}(Y, h)$ depends on the set $C_{Y}$. Changing $C_{Y}$ into $C'_{Y}$ one has by Equation \eqref{dyt} that
\bq
\lam_{\cd}^{C'_{Y}}(Y, h)=f^{C_{Y},C'_{Y}}_Y(h)\lam^{C_{Y}}_{\cd}(Y, h)
\eq
for all $h \in H$.
\er}
\bibitem{gnn}  \textsc{S. Gelaki} and {\sc D. Nikshych} and {\sc D. Naidu}, \emph{Centers of graded fusion categories}, Alg. Num. Th.  \textbf{ 3}, 959-990, (2009).
\bibitem{Ostr} {\sc V. Ostrik}, \emph{Module categories, weak Hopf algebras and modular invariants,} Transf. Groups. \textbf{8}, 177-206, (2003).
\bibitem{mns} {\sc  Maier, J.}, {\sc Nikolaus, T.}, {\sc Schweigert, C}., \emph{Equivariant Modular Categories
via Dijkgraaf-Witten Theory}, arXiv:1103.2963  (2011).

\bibitem{ty} {\sc D. Tambara}, \emph{Tensor categories with fusion rules of self-duality for finite abelian groups} J. Algebra 209 (1998), no. 2, 692-707.

\bibitem{agaitsgory} {\sc S. Arkhipov} and {\sc D. Gaitsgory}, \emph{Another realization of the category of modules over the small quantum group}, Adv. Math. \textbf{173} (2003), 114--143.

\bibitem{clifth} {\sc S. Burciu}, \emph{Clifford theory for cocentral extensions}, Israel J.\ Math. \textbf{181}, 111--123 (2011).

\bibitem{bnda} {\sc S. Burciu}, \emph{Normal Hopf subalgebras of Drinfeld doubles}, preprint  (2013).

\bibitem{dpr}  {\sc R. Dijkgraaf} and {\sc V. Pasquier} and {\sc Ph. Roche, } \emph{Quasi-Quantum Groups Related to Orbifold Models,} Proceedings of the ``International Colloquium on Modern Quantum Field Theory'', Tata Institute of Fundamental Research, 375-383 (1990).

\bibitem{ENO2}  {\sc P. Etingof}, {\sc D. Nikshych} and {\sc V. Ostrik}, \emph{Weakly group-theoretical and solvable fusion categories}, Adv. Math. \textbf{226}, 176--205 (2011).

\bibitem{agaitsgory} {\sc S. Arkhipov} and {\sc D. Gaitsgory}, \emph{Another realization of the category of modules over the small quantum group}, Adv. Math. \textbf{173} (2003), 114--143.
\bibitem{goff} {\sc C. Goff}, \emph{Fusion rules for abelian extensions of Hopf algebras}, to appear, Alg. Number Theory \textbf{}, (2011).
\bibitem{AD}{N. Andruskiewitsch and Devoto}, Extensions of Hopf algebras, Alg. i Analiz (1995).
\bibitem{ext-ty} {\sc S. Natale}, \emph{Hopf algebra extensions of group algebras and Tambara-Yamagami categories}, Algebr. Represent. Theory \textbf{13} (6),  673--691 (2010).
\bibitem{nat-faith} {\sc S. Natale}, \emph{Faithful simple objects, orders and gradings of fusion categories }, to appear Algebr. Geom. Topol. Preprint arXiv:1110.1686, (2011). 

\bibitem{ws} {\sc S. Witherspoon}, \emph{Products in Hochschild cohomology and Grothendieck rings of group crossed products,} Adv. Math. \textbf{185}, 136ï¿½158, (2004).

%\ed

\section{Applications}\label{appl}
\subsection{Invertible objects in $\cc^{G}$}
\begin{theorem}
One has that $\innv(\cc^{G})=Eq\innv(\cc)\times \innv(\rep(G))$ where $Eq\innv(\cc)$ is the set of all invertible objects in $\cc$ that are equivariant.
\end{theorem}
\bpf
Suppose that $S_{X, \pi}$ is an invertible object of $\cc^{G}$. Recall that
\beq
\dim S_{X, \pi}=[G:G_{X}]\fp(X)\dim \pi
\eeq

Conversely, for any pair $(X, \pi)\in Eq\innv(\cc)\times \innv(\rep(G))$ define the element $S_{X,\pi}:X\otimes V_{\pi}$.
\epf
Thus $G_{X}=G$ and $\fp(X)=1$ i.e $X$ is an invertible object of $\cc$. Moreover $\pi$ is a one dimensional representation of $G$ with cocycle $\al_{X}$. It follows that in fact $\pi$ is a one dimensional representation. Thus $(X, \pi)\in Eq\innv(\cc)\times \innv(\rep(G))$.
\begin{corollary}
One has that $\innv(\cc^{G})=\cc([G,G], K_{0}
, \cs, 1)$.
where $K_{0}\subset \cz(G)$ is the support of all equivariant $G$-objects of $\cc$.
\end{corollary}
\bpf
One has that $\innv(\cc^{G})\cap \rep(G)=\rep(G/[G,G])$.
\epf

\subsection{Nondegenerate fusion subcategories}
\bp One has that $Z_{2}(\cc(H,K, \cs, \lam))=\cc(\ker, HK, \tilde\cs, \tilde\lam)$.
\ep
\bpf
One has that $Z_{2}(\cc(H,K, \cs, \lam))=\cc(H,K, \cs, \lam)\cap (\cc(K, H'', \mathcal{T}, \tilde \lam)=\cc(HK, K_{int}, \tilde \cs))$
\epf
\blue{\bp Suppose that $G=\mathbb{Z}_{p^{n}}$ and $\cc$ is faithfully graded by $G$ with $\cc_{1}$ nondegenerate and $\mathrm{K}_{0}(\cc)$ a commutative ring. If $\cc(H, K, \cs, \lam)$ is nondegenerate then $\{K,H\}=\{G, 1\}$ and . \ep}
\bll{for $2^{n}$ it should be prime.}
\bpf
By definition one has that $\cc(H, K, \cs, \lam)$ is nondegenerate if and only if $\cc(H, K, \cs, \lam) \cap \cc(H, K, \cs, \lam)'=\mathrm{Vec}$. Suppose that
$\cc(H, K, \cs, \lam)'=\cc(K, H, \mathcal{T}, \tilde \lam)$. It follows from Proposition \ref{int} that $HK=G$. Since $G$ is cyclic it follows that $H=G$ or $K$=G. \onh \; Corollary \ref{veecom} implies that $H \cap K=1$. Thus $\{K,H\}=\{G, 1\}$. On the other hand $H\subseteq \cap_{Y \in \irr(\cs)G_{Y}}$ implies that $G_{Y}=G$ for any $Y \in \cs$. Thus $\cs$ is {\it  fixed} by the action of $G$.
\mdn
One may suppose that $H=G$ and $K=1$. Thus $\cs \subset \cc_{1}$ and $\fp(\mathcal{T})=\frac{\fp(\cc)}{\fp(\cs)}$\epf
\bn{defn}
We say that a twisted bicharacter $\lam:\irr(\cs)\ra \hat H$ is {\it nondegenerate} if  and only if for any $S \in \cs$ the relation $\lam([S],-)=1$ implies that $S\simeq 1_{\cc}$.
\end{defn}
\bp
A fusion subcategory $\cc(H, K, \cs, \lam)$ is non degenerate if and only if $\lam$ is non degenerate.
\ep

\bpf
$\cd$ is nondegenerate iff $\cd'\cap\cd=Vec=\cc(G, 1, Vec, 1)$. One has that $\cc(H, K, \cs, \lam)\cap \cc(K, H, \mathcal{T}, \tilde \lam)=\cc(HK, K_{int}, \ker \phi_{\lam, \tilde\lam})$. Thus $HK=G$ and $K_{int}=1$, i.e $\lam \tilde\lam$ is nondegenerate.
\epf
\subsection{M\"{u}ger's centre of $\cc^{G}$}
\br\label{suppone} One has that
$\cc(1, G_{1}, \cc_{1}, 1)\simeq \cc_{1}^{G}$ is the largest subcategory of $\cc^{G}$ supported on $\cc_{1}$. Moreover the forgetful functor restricted to this subcategory is a braided tensor functor.
\er
\subsection{Slightly degenerate fusion subcategories}

\bl
One has that $\cc^{G}$ is slightly nondegenerate if and only if $G_{1}=G$ and $\cc_{1}$ is slightly nondegenerate.
\el
\bl
If $\cc(H,K, \cs, \lam)'=sVec$ then the grading is faithful and $K=G=\mathbb Z_{2}$ and $\cz_{2}(\cc_{1})\supseteq sVec$.
\el
\bpf
One has that $\cc(H,K,\cs, \lam)'=\cc(K, H, \mathcal{T},\tilde\lam)$. Thus $\rep(G/K)\subset \mathrm{sVec}$ which implies that $K=G$. But $K\subseteq G_{1}$, thus $K=G_{1}=G$. One the other hand $2=\fp(\cc(G, H, \mathcal{T},\tilde\lam))=\fp(\mathcal T)$ implies that $H\leq 2$ On the other hand Corrolary \ref{czc1} implies that $H=1$. Thus $\mathcal T\subset \cc_{1}$. The above remark implies that $\mathcal T$ is $\mathrm{sVec}$ and centralizes any object of $\cc_{1}$, thus 
\epf
\subsection{M\"uger center of $\cc^{G}$}
By definition one has that $\lam_{\cd}([X],1)=1$ for any fusion subcategory $\cd\subseteq \cc^{G}$. \mdn Thus if $H_{\cd}=1$ then $\lam=1$. Note that one has:
$\vvec=\cc(G, 1, \vvec, 1)$ and $\cc^{G}=\cc(1, G_{1}, \vvec, 1)$.
The equality  $\;\;\vvec'=\cc^{G}$ can be written as $\cc(G, 1, \vvec, 1)=\cc(1, G_{1}, \cc)$ which shows that the subgroup $H' \neq H$ if $G_{1}\neq G$.
\bp
One has that $(\cc^{G})'=\cc(1, G_{1}, \cc, 1)'=\cc(G_{1} , 1, \mathcal T_{1},1)$ for some fusion subcategory $\mathcal{T}_{1}\subseteq \cc_{1}$ where $G_{1}$ is the grading group.
\ep\bl
If $\cc$ is a braided category containing a tannakian subcategory then $\cc$ cannot be slightly degenerate.
\el
\bpf Let $\cd=\cc_{G}$.
One has that $\cc=\cc(1, G_{1}, \cc,1)$. Thus $\cc'=\cc(G_{1}, 1, \mathcal T,\lam)$ for some $\mathcal T\subset \cd_{1}$. If $\cc'=sVec$. Then $G_{1}=G$  and $\fp(\mathcal T)=2$.
\epf
\bp Let $K$ be a normal subgroup of $G$. With the above notations one has that 
\beqn
\cc(1, K, \cc(K), 1)'=\rep(G/K)
\eeqn
where $\cc(K)=\oplus_{k \in K}\cc_{k}$.
\ep

\bpf
This is clear form the formula for the centralizer.
\epf
\bp
One has that $\cc(1,1, \cs, 1)'=\cc(1,1,\cs', 1)$. In particular $\cc(1,1,\cc_{1}, 1)'=\cc(1, 1, \cz_{2}(\cc_{1}), 1)$.
\ep

\subsection{Symmetric and Lagrangian's subcategories}
We say that a fusion subcategory $\cs$ of $\cc$ is $G$ symmetric if $X\perp _{G}Y$ for any objects $X,Y \in \cc$.
\bp
$\cc(H,K,\cs,\lam)$ is symmetric if and only if $K\subseteq H$, $\cs$ is $G$-symmetric and there is a chosen system $c$ such that $\lam$ verifies the following property:
\beq
\lam_{c}([X_{k}],l)\lam_{c}([Y_{l}], k)=\omega_{c}([X_{k}], [Y_{l}])
\eeq
for all $l, k \in K$, $X_{k}\in \cc_{k}$ and $Y_{l}\in \cc_{k}$.
\ep
\bpf
One has that $\cc(H,K,\cs,\lam)$ is symmetric if and only if $\cc(H,K,\cs,\lam)\subseteq \cc(H,K,\cs,\lam)'$. Let as before $\cc(H,K,\cs,\lam)'=\cc(K, H, \tilde\cs, \tilde \lam)$. Thus $\cc(H,K,\cs,\lam)\subseteq \cc(K,H, \tilde \cs, \tilde\lam)$. Proposition \ref{incl} implies that $K \subset H$ and $\cs\subset \tilde \cs$, i.e $\cs$ is a $G$-symmetric fusion category. Moreover, since $\lam|_{\irr(\cs)\times K}=\tilde\lam|_{\irr(\cs)\times K}$ it follows that $\lam$ verifies the above condition.
\mdn
The converse is immediate from the criterion of two objects to centralize each other.
\epf
\bll{Define alternating $\lam_{c}$. Prove that it is alternating with respect to a system then it is alternating to any other system.}
\bp
A fusion subcategory is Lagrangian if and only if $H=K$, $\cs$ is $G$-Lagrangian and $\lam$ is an alternating function.
\ep
A fusion subcategory of $\cc$ is called Lagrangian if and only if $\cs^{\perp_{G}}=\cs$. 
\bpf
\epf
%\subsubsection{Fusion subcategories as equivariantizations}

\bp Suppose that $K_{0}(\cc)$ is a commutative ring. Then if $\cc^{G}$ not prime then $G$ decomposes  nontrivially as $G=H\times K$.\ep
\bpf
Suppose that $\cc^{G}$ is not prime. Then there is a 
\beqn
\cc(H, K, \cs, \lam)\boxtimes \cc(K, H, \tilde \cs, \tilde \lam)=\cc^{G}
\eeqn
Then $\cc(H, K, \cs, \lam)\cap \cc(K, H, \tilde \cs, \tilde \lam)=\mathrm{Vec}$ which implies that $KH=G$. Moreover $\cc(H, K, \cs, \lam)\vee \cc(K, H, \tilde \cs, \tilde \lam)=\cc(H\cap K, )=\cc(1, G, \cc, )$ which implies that $H \cap K=1$.
\epf
\bll{In $\cc_{1}$ two objects $G$-cross centralize each other if and only if they M\"{u}ger centralize since $c_{Y}^{1}=\id_{Y}$ for any $Y\in \cc$.}
\bc 
Suppose that $\cc^{G}$ is a nondegenerate braided fusion category.
If $\cc^{G}$ is prime then $\cc_{1}$ is prime.
\ec
Note that $\cc(1,1,\cd, 1)\subseteq \rep(G)$.
\bpf
Suppose that $\cc_{1}=\cd\boxtimes \cd'$. Then using Theorem it follows 
\beqn
\cc(1,1, \cd, 1)'=\cc(1,1,\cd 1).
\eeqn
\epf
\bc
Criterion for $\cc^{G}$ to be group theoretical.
\ec

\bll{Dyslectic modules with respect to the commutative algebra $Fun(G/H)$ are those supported on  $\oplus_{h \in H}\cc_{h}$.} \mdn

\bp
Suppose that $\rep(L)$ is another Tannakian  fusion subcategory of $\cd=\cc^{G}$. Suppose that
$\rep(L)=\cc(H_{2}, K_{2}, \cs_{2}, \lam_{2})$.
\ep
\bpf
Since $\rep(L)\subseteq \rep(L)'$ one has that $K_{2}\subseteq H_{2}$ and $\lam|_{\irr(\cs)\times K_{2}}=1$.
\epf
\bll{In the pointed case $\tilde \cs$ is not the maximal one, it is the maximal one supported on $H$.}
\mdn
\bll{See also the correspondence from 4.30 dgno}
\subsection{Applications to the center of a fusion category}
\bll{\bp
If $G=H\times K$ then $\cz(\cc)$ si decomposable
\beq
\cz(\cc)\cong \cc(H, K)\times \cc(K, H)
\eeq
\ep}
\bpf\bll{Define $\cc(K, H)$ as all those that $S_{Y,\pi}$ such that $G_{Y}\supset H$ is $H$ equivariant and $\pi|_{H}=\eps$. It is easy to compute the dimension of this category:
\beq
\sum_{Y\;|G_{Y}\supset H}||S_{Y, \pi}|^{2}
\eeq}
\epf

\subsection{What do I want to happen:}
\bt
Two objects $X$ and $Y$ $G$-cross centralize each other if and only if 
\beq
N(g)\subset G_{Y}
\eeq
\beq
N(h)\subset G_{X}
\eeq
\beq
d^{p,1}_{x,y}d^{1,p}_{y,x}=\omega_{p}\id
\eeq
\bll{(or maybe it is enough to be a scalar of the identity)}
\mdn
and 
\beq 
\al_{X}\times \al_{Y}\text{ \;\;is cohomologus ttivial on}\;\; N(h)\times N(g).
\eeq 
\et

\bt
If $X$ and $Y$ are such that they $G$- centralize then for any representation $\pi$ with cocycle $\al_{X}$ there is at least one $\delta$ such that $S_{X, \pi}$ cross centralize $S_{Y,\delta}$.
\et\mdn
\red{It is an elementary problem to decide if it is enough for the cocycles to be cohomologus trivial in order for the second theorem to be true.}
\mdn
\red{As a consequence of these two theorem it follows that $\mathcal T$ is the largest that cross centralize $\mathcal S$ and it is supported on $H$.
}
\subsection{If the grading group $G$ is simple non abelian}
Suppose that the grading group is simple nonabelian. Let $\cc(H,\cs, \lam)$ a fusion subcategory of $\cc$. Since $H, K$ can be onlu $1$ and $G$ itself we have the following three possibilities:

$\cc(1, 1, \cs, 1)$ where $\cs$ is a fusion subcategory of $\cc_{1}$ which is invariant by $G$.

$\cc(G, 1, \cs, \lam)$  where $\cs$ is a fusion subcategory of $\cc_{1}$ which is invariant by $G$.

$\cc(1, G, \cs, \lam)$  where $\cs$ is a fusion subcategory of $\cc$ whose support is the entire group $G$.
\bc
If the group $G$ is simple then any fusion subcategory is an equivariantization of $G$.
\ec
%\mdn \bll{For sure it is true at the level of tensor products.}
%\blue{Note that the equivariant structure sends $T^g(Y_h)$ into $Y_{ghg^{-1}}$.}
\bb{\bne \item Find an example where the grading group is not too small.
\item In the case of the center $\cz(\cc)$ give the crossing structure of $S(a, M)$ in terms of the crossing structure of $M$ regarded as object of $\cz_{\cc_{1}}(\cc_{a})$ and the equivariantized structure of $M$.
\item Determine other Tannakian subcategories by rewriting in the other system.
\ene}
\subsection{The nondegenerate case}
\bab
Use Sonia's explanations, $\cc^{G}$ is non-degenerate if and only if $\cc_{1}$ nondegenerate and in this case 
\beq
\cz(\cc)\simeq \cc^{G}\boxtimes \cc_{1}^{rev}.
\eeq
\eab
\section{Formulae}

Suppose that $K\subset G_{1}$. Then
 \beq
\rep(G/K)=\ca(K, 1, \vvec, 1),\;\;\rep(G/K)'=\ca(1, K, \cc(K),1)
\eeq
More general
\beq
\rep(G/K)'=\ca(1, G_{1}\cap K, \cc(K\cap G_{1}),1)\ca
\eeq
\bll{see in the their proof how the intersection appear}\mdn
Suppose that $X\perp_{G}Y$.
\mdn \bne 
\item look at the twisted drinfeld double what happens; 
\item put it as a particular example; 
\item take the formulae from the paper with sonia on the fusion rules.
\ene

\mdn
\bll{Suppose that $X \perp Y$. Recall the scalar $\omega([X], [Y])$ the scalar from the composition homomorphisms of Equation \eqref{perp}. Denote by $\cs^{\perp_{G}}$ the abelian subcategory consisting of all objects that cross centralize any object of $\cs$.
}
%\subsection{Some invariants}
%\input invarianti braided.tex
\subsection{Changing systems to become more transparent} 
\br \label{oneomega}Note that one can change the systems of chosen isomorphisms of $X$ and $Y$ by the roots of unity arising from each $\pi(php^{-1})$ and $\delta(p\inv gp)$ such that $\omega_{p}([X], [Y])$ disappear. Since the systems of chosen isomorphisms are also {\bll {\it good}} it follows that $\pi(php\inv)=\omega_{p}^{(1)}$ and $\delta(p\inv gp)=\omega_{p}^{(2)}$ for some roots of unity $\omega_{p}^{(i)}$ with the property $\omega_{p}^{(1)}\omega_{p}^{(2)}=\omega_{p}^{-1}$. Changing $c_{Y}^{pgp^{-1}}$ to $\omega_{p}^{1}c_{Y}^{pgp^{-1}}$ and $c_{X}^{p^{-1}gp}$ to $\omega_{p}^{2}c_{Y}^{pgp^{-1}}$ note that in these two new systems one has that
\beq
d_{X,Y}^{1,\;p}d_{Y, X}^{p,\;1}=id_{X\ot \ct^{p}(Y)}
\eeq 
Then Corrolary \ref{cp} imposes that in this case
\beq
\pi(pgp^{-1})\ot \delta(p^{-1}hp)=\id_{V_{\pi}\ot V_{\delta}}
\eeq
\er 
\noindent For the rest of the paper we denote $d_{X,Y}:=d_{X, Y}^{1,1}$.  
\bab
\bl
With the above notations note that $d_{\ct^{p}(Y), X}=d_{Y, X}^{p,\;1}$. Moreover, if the system \bll{$c_{Y}$ is conjugate closed } then $d_{X, \ct^{p}(Y)}=d_{X, Y}^{1,\; p}$.% and  where for $\ct^{p}(Y)$ it is chosen the conjugate system $\;^{p}c_{Y}$.
\el
\bpf It is also immediate that $d_{\ct^{p}(Y), X}=d_{Y, X}^{p,1}$.  Straightforward computation using the definition of $\;^{p}c_{Y}$.
\epf
%\bll{Rewrite everything so only $d_{X, Y}$ and $\omega_{X, Y}$ appears. But in this way we loose the form from their paper.}
\begin{corollary} \label{c22}Suppose we have chosen two good systems $c_{X}$ and $c_{Y}$ such that $c_{Y}$ is conjugate closed. With the above notations the objects $S_{X,\pi}$ and $S_{Y, \delta}$ centralism each other if and only if for any $p\in G$ there is a roots of unity $\omega_{p}\in k^{*}$  such that the following two relations hold:
\beq
d_{X, \ct^{p}(Y)}d_{\ct^{p}(Y), X}=\omega_{p}\id_{X\ot \ct^{p}(Y)}
\eeq
and
\beq\label{reform'}
\pi(php\inv)\ot \delta(p\inv gp)=\omega_{p}^{-1}\id_{V_{\pi}\ot V_{\delta}}
\eeq
\end{corollary}
%\mdn\bll{if $p \in G_{X}\cap G_{Y}$ then the value of $\omega_{p}$ clearly determined by $\omega_{1}$.}\mdn
\br Note that in the above situation $\al_{X}\times \al_{Y}$ is cohomologus trivial on the subgroup $N(h)\times N(g)\subset G_{X}\times G_{Y}$. Here $N(g)$ is the smallest normal subgroup of $G$ and clearly $N(g)\subseteq G_{Y}$\er
\bab
\begin{remark} In the case that $G$ is abelian  there is just one condition that has to be satisfied by the value at $g$. Let $N_{Y}(g)$ be the smallest normal subgroup $N(g)\subseteq G_{Y}$ and $\delta$ has a fixed value on $N(g)$. Then $\delta$ is a constituent of the induced representation from $N_{Y}(g)$ to  $G_{Y}$. \end{remark}\eab
\mdn
\bab
The above relations with root of unity hold only for good chosen systems $c_{Y}$ and $c_{X}$.\eab
\br
Suppose that $G$ is an abelian group.  Let $X \in \cc_{g}$ and $Y \in \cc_{h}$ and suppose that the two simple objects $S_{X , \pi}$ and $S_{Y, \delta}$ centralize each other in $\cc^{G}$. Then for any $m\in G$ the two identities above can be written as follows:
\beq
d_{X, \ct^{m}(Y)}d_{\ct^{m}(Y), X}=\omega\id_{X\ot \ct^{m}(Y)}
\eeq
and
\beq
\pi(h) \ot \delta(g)=\omega^{-1} \id_{V_{\pi}\ot V_{\delta}}
\eeq
Thus in this situation $S_{X, \pi}$ centralizes $S_{Y, \pi}$ if and only if $\deg(X)\in G_{Y}$, $\deg(Y) \in G_{X}$ and there is $\omega \in k^{*}$ root of unity such that the above two identities are satisfied.
\er
\eab

\bll{Since $\al_{Y}|_H$ is cohomologus trivial one can choose a set of isomorphisms $c_{Y}$ such that $\al_{Y}|_{H}=1$. \bab Moreover by Lemma \ref{cjt} one may suppose that $c_{\ct^{g}(Y)}=\;^{g}c_{Y}$ for any $g \in G$. Same argument with chosen representatives of the orbits.\eab}\mdn
\bl\label{cjt} Suppose that $Y \in \cs$ and that the chosen set of isomorphisms $c_{Y}$ is given such that $\al_{Y}|_{H}=1$. Then also $\al_{\ct^{g}(Y)}|_{H}=1$ for the chosen set of isomorphisms $\;^{g}c_{Y}$ for $\ct^{g}(Y)$.
\el 
\bpf
The result follows from the commutativity of the following diagram which express $\al_{\ct^{g}(Y)}(ghg^{-1}, glg^{-1})$. \bll{The big diagram from the yellow notebook.}
\epf
\subsubsection{On the dependence inside the isomorphism class}
\mdn For any other object $Y'$ isomorphic to $Y$ we choose the following isomorphism $c_{Y'}^{h}=f^{-1}c_{Y}^{h}\ct^{h}(f)$ for all $h \in G_{Y}$ where $f:Y'\ra Y$ is a fixed isomorphism. Then it is a straightforward computation to show also that $\al_{Y'}(h,h')=1$ for all $h,h'\in H$. \bll{Thus $c_{Y'}$ is also a trivial $H$-system at $Y'$.}
\bl With the above chosen sets of isomorphisms for $Y$ and $Y'$ it follows  $\al_{Y'}(g_{1}, g_{2})=\al_{Y}(g_{1}, g_{2})$ for any $g_{1}, g_{2}\in G_{Y}$. Moreover one has that $f\ot \id_{\pi}:Y\ot \pi \xra{\simeq} Y'\ot \pi$ is an isomorphism in $\cc^{G_{Y}}$.
\el
\bpf
One has to check that  for any $h \in H$ the following equation holds
\beqn 
(c_{Y'}^{h}\ot \pi(h))(T^{g}(f)\ot \id_{\pi})=(f \ot 1)(c_{Y}^{h}\ot \pi(h))
\eeqn
This follows directly from the definition of $c_{Y'}^{h}$.
\epf
\noindent
It follows that in this situation \beq S^{c_{Y'}}_{Y',\pi}\simeq \ind_{G_{Y}}^{G}(Y'\ot \pi)\simeq \ind_{G_{Y}}^{G}(Y\ot \pi)\simeq S^{c_{Y}}_{Y,\pi}.\eeq  Therefore, from the above Lemma we may suppose that for any  simple object $Y$ of $\cs_{\cd}$ the set of chosen isomorphisms $c_{Y}$ of $Y$ is given such that $\al_{Y}|_H=1$. 
\mdn
{\bf For shortness, for the rest of this paper we denote by $S_{Y,\pi}$ the simple object of $\cc^{G}$ corresponding to this chosen set of isomorphisms $c_{Y}$ of $Y$.} \mdn

\mdn \bll{$\lam$ still depends on the chosen set of isomorphisms. In order to get the $G$-invariance without $d_{Y}(g,h)$ one needs to assume that in the system $C$ the system of $c_{\ct^{g}(Y)}$ is given by $\;^{g}c_{Y}$. }
\np
{\bf Questions:} 
 \bne
 \item \bll{See Propo 4.56 from dgno, $\cc'$ is not defined only some relations are known, the intersection with $\rep(G)$ and the wedge with $\rep(G)=(\cc_{1}^{G})$.}\item 
  \bll{Omega also might depend on the class of the chosen set of isomorphisms. for $X$ and $Y$.}
 \item show properties of $d_{X,Y}^{m,n}$.
 \item \bll{Is any relation in the Grothendieck group that records when two objects $G$-cross centralize?}
 \item Show that the system from \cite{nnw} is both good, $H$-trivial and closed conjugate.
 \item
 \bll{Thus $X\perp_{G}Y$ if and only if $X\perp \ct^{m}(Y)$ for any $m \in G$.}
 \item Given a fusion subcategory $\cs$ is the  
 abelian subcategory $\cs^{\perp_{G}}$ also a fusion subcategory? It seems that the answer is NO from next Example.
 \item Does  $\mathcal{T}$ consists of all objects of $\cs^{\perp_{G}}$ whose degrees belong to $H$?
 \item
 Suppose that $X \perp Y$, i.e. $d(1,1)_{X, Y}d_{Y, X}(1, 1)=\omega(X, Y)\id_{X \ot Y}$. Does it follow that $\omega(X, Y)$ is a root of unity?
 \item
 Suppose that $X \perp_{G} Y$ with
 \beq
 d_{X, Y}^{m,n}d_{Y, X}^{n,m}=\omega(\ct^{m}([X]),  \ct^{n}([Y]))\id_{X \ot Y}
 \eeq
 Given $\pi  \in \Irr_{\al_{X}}(G_{X})$ define 
 \beq
 \delta(ngn^{-1})=\frac{\pi(1)}{\pi(mhm^{-1})}\omega(\ct^{m}([X]),  \ct^{n}([Y]))
 \eeq
 Is it $\delta$ well defined and a 
 $\al_{Y}$-projective representation of $G_{Y}$?
 \ene
 \bne
\item Make a section with projective representations and their characters, multiplicity?
\item what kind of multiplicity is used in the fusion rules?
\item lifting a projective representation to a linear representation the character is the same for the lifting? so we can get directly the inequality.
\item put a an example with NNW, how to get multiplicative in one side.
\item maybe we should put $K$ first to obtain the same notation as they have.
\ene
\np
\bne
\item use the description of the simple objects from g-funtors as $S(g, M)$. 
\item give a new condition when two such objects cross centralize.
\item
there are gonna be new $d_{M,N}d_{N, M}=\id_{M\otimes N}$
\item use the other property of $\lam$ to get the relations between omegas'
\item One has that $X_{c, -\hat b}\otimes X_{b, -\hat c}=X_{a, -\hat a}$ where $a=b+c$. None of them does cross centralize $Z_{\ro}$ but their product it does.
\item for the beginning I can suppose that the action is strict, i. $\ct^{g}\ct^{h}=\ct^{gh}$.
\item take $\delta$ from Remark 4.15, is it a projective representation on $N(g)$?
\item use corollary 3.6 from \cite{dgno} the rank of $\cd$ is the number of components of $\cd'$.
\item change the function $d_{Y}(g,h)$ to another letter.
\item sync with fsc.tex, put the braided part from smc over the braided part form fsc
\item Here one needes the compatibility with the unit $c_{\1, X}=\id_{X}$. Maybe it results from the other? $\tilde{c}$ is a braiding and it should be the identity
\item square versus squared
\item from the first relation the crossed centralizer relation I don't think it is transitive
\item The question is if all crossed centralizing objects of a given object form a fusion subcategory?
\item two objects cocentralize each other if and only if there are simple objects seating over them which M\"{u}ger centralize each other.
\item From the S-matrix compute explixitly which object cross centralize; then try to derive the formulae for the categories.
\item compute an explict example from their papers to deduce what are the stable fusion subcategories. $A=\mathbb Z_{p}$
\item the simple objects of $\cz(\cc)$ are given in GNN by the equivariant structure.
\item the second cohomology is trivial.
\item isomorphism chosen such that all cocycles are exactly 1.
\item Finish the last two cases so they do not depend on the $\eps$ from the initial category.
\item $Z_{\ro}\perp_{G}Z_{\ro'}$ does not seem to give the same thing as computing directly the $S$-entry. $Z_{\ro, \Delta}$ centralizes $Z_{\ro', \Delta'}$ It possibly explains why the entries in modulus is less or equal $n$ and equality might hold only in the conditions from cross centralizer between $Z_{\ro}$ and $Z_{\ro'}$. 
\item
Try also generalized Tambara-Yamagami categories.
\item Establish completely the detramnination $\eps \mapsto \eps_{a}$ and $\Delta \mapsto \eps_{\Delta}$.
\item Put the other remark with $b=1$ in the formula for $\tilde \lam$.
\item
$|G|=p^{n}$ then $\cc^{G}$ is prime for $p=2$.
\bpf
Suppose that
$\cc(H, K, G)'=\cc(K, H, \mathcal{T}, \tilde \lam)$. It follows that $HK=G$ and $H \cap K=1$
\epf
\item If $G$ is abelian then $\cc^{G}$ is the Deligne product of some fusion categories.
\item $G/H$ acts on $\cs$ since $H \subset G_{V}$ for all $V \in \cs$. Is the equivariantization $\cc(H,K, \cs, \lam)$?
\mdn They have the same Perron Frobenius dimension and the simple objects look very similar. $Stab_{G/H}(V)=G_{V}/H$
\item Cosets of $\cc(H, K, \cs, \lam)$. Apply the result of Drinfeld.
\item In Drinfeld doubles of groups
\beq
(S_{a, \gm})'=\cc(C_{G}(N(a)), N(a), \lam)
\eeq
The conjecture is that the number of components of $\cc(C_{G}(N(a)), N(a), 1)$ equals the number of conjugacy classes of $G$.
\ene
\subsection{Some canonical objects of the center $\rep(D(A))$.}
Suppose that $\cc$ admits a (quasi-)functor, i.e. is the category of representations of a (quasi-)Hopf algebra, $A$. Then $A \in \cz(\cc)$ via the half braiding given by the inverse of:
\beq A \ot M \ra M\ot A,\;\;
a \ot m \mapsto a_{1}Sa_{3}m\ot a_{2}
\eeq

Then  by restriction this map induces isomorphisms $\phi_{g, M}:R_{g} \ot M \ra M \ot R_{h^{-1}gh}$.
\bp
There is  a functor \bb{an embedding} of categories
\beq
\rep(G_{a})\hookrightarrow \rep(D(A))
\eeq
given by $V \mapsto T(a, V):=\oplus_{x \in G/G_{a}}R_{xax^{-1}}\ot V$
 where the commuting structure 
 \beq
 M \ot T(a, V) \ra T(a, V) \ot M
 \eeq
 is given on components by 
 \beq M\ot R_{xax^{-1}}\ot V \xra{\phi_{M, g}\ot \pi_{V}(h)} R_{zaz^{-1}}\ot M \ot V
 \eeq
 \ep
 \bpf
 It is easy to verify that this element is central.
 \epf 
\bp\beq
T(a, V)\ot T(b, W)=\oplus_{x \in D}T(\;^{x}ab, m_{\;^{x}a,b}(\;^{x}V, W))
\eeq
\ep
\bl
The regular characters of each summand is $T(a, \mtr{reg}_{G_{a}})$.
\el
\subsection{Admitting good systems and $H$-trivial systems of chosen isomorphisms}
A system $c_{Y}$ of chosen isomorphisms is called \bll{{\it good }} if the values of the character $\al_{Y}:G_{Y}\times G_{Y}\ra k^{*}$ are all roots of unity.
\mdn
It follows from Equation \eqref{conjcc} that if $c_{Y}$ is good then $\;^{g}c_{Y}$ is also good for any $g \in G$.
\mdn
\bll{I can choose a global system in which $\;^{g}c_{Y}=c_{\ct^{g}(Y)}$ for any $g \in G$ and any $Y$ simple in $\cc$, see the appendix.}
\mdn\bll
{Does the system remain $H$ trivial?}
\bll{This is a fact that wasn't noticed in \cite{nnw}, that you can simplify a little bit.}

%\ed
\subsubsection{Trivial systems}
 Fix an element $Y$ for any isomorphism class of the simple objects of $\cs$. Suppose that $S_{Y,\pi}\in \cd$.  Since each $h \in H_{\cd}$ acts as a scalar on $\pi$ it follows that the cohomology class of 
$\al_{Y}|_{H}$ is trivial. Thus one has that $\al_{Y}(h, h')=t(hh')t(h)^{-1}t(h')^{-1}$ for some function $t:H\ra k^{*}$. Then changing the isomorphisms $c^{h}_{Y}$ to $t(h)c_{Y}^{h}$ we may suppose that $\al_{Y}(h, h')=1$ for all $h,h'\in H$. \mdn We call $c_{Y}$ a \it trivial $H$-system at $Y$ if $\al_{Y}^{c_{Y}}|_{H}=1$. Note that the set of all trivial $H$ systems at $Y$ are in bijection with the set $\hat H:=Hom(H, k^{*})$.
\section{Examples: The Centers of Tambara-Yamagami categories}\label{example}
\bll{consider the action form Turaev-Virelizier and see if from there one can write down the tensor structure}
\subsection{The center of a fusion category} Let $
\cc$ be a  fusion category graded by a finite group $G$ 
with $\cd$ the trivial component of the grading.
Following \cite{gnn} one has that $\cz_{\cd}(\cc)^{G} \cong \cz(\cc)$. 
\subsection{Definition of the Tambara Yamagami category}
 In \cite{ty} D. Tambara
and S. Yamagami completely classified all $\mathbb{Z}_{2}$-graded fusion categories in which
all but one simple object are invertible. They showed that any such category
$T Y(A, \ch, \tau)$ is determined, up to an equivalence, by a finite Abelian group $A$,
a non-degenerate symmetric bilinear form $\ch : A \times A \ra k^{*}$, and a square root
$\tau \in k$ of $|A|^{-1}$. The category $T Y(A, \ch, \tau)$ is described as follows. It is a skeletal
category (i.e., such that any two isomorphic objects are equal) with simple objects
$\{a | a \in A\}$ and $m$, and tensor product
\beq
a \ot b = a + b, a \ot m = m, m \ot  a = m, m \ot m =
\sum_{a \in A}a\eeq

for all $a, b  \in A$, and the unit object $1\in A$. The associativity constraints are given in \cite{ty}. One has that $\cc:=T Y(A, \ch, \tau)$ is $\mathbb{Z}_{2}$-graded fusion category with $\cc_{0}=\cd=\{a\;|\;a \in A\}$ and $\cc_1=\{m\}$.
\mdn \bll{write down the associativity constraints.}
\subsection{Description of the relative center}

\subsubsection{Simple objects of the fusion category $\Z_{\cd}(\cc)$} Clearly $\Z_{\cd}(\cc)=\Z(\cd)\oplus \Z_{\cd}(\cc_{1})$.

Following \cite [ Subsection 4.2]{gnn} simple objects of $\Z(\cd)$ are of the type $X_{a, \phi}$ where $\phi \in \hat{A}:=\mathrm{Hom}(A, k^{*})$. As objects of $\cc$ one has $X_{a, \phi}=a$ and the relative central structure $X_{a, \phi} \ot b \ra b \ot  X_{a, \phi}$ given by $\phi(b)\id_{a+b}$ for all $b \in A$.

Moreover simple objects of  $\Z_{\cd}(\cc_{1})$ of $Z_{\ro}$ with $\ro: A \ra k^{*}$ satisfying $\ro(a+ b)=\ch(a, b)^{-1}\ro(a)\ro(b)$. As objects of $\cc$ one has $Z_{\ro}=m$ and the relative central structure $ Z_{\ro} \ot b \ra b \ot  Z_{\ro}$ given by $\ro(b)\id_{m}$ for all $b \in A$.

Let $A \ra \hat A $ given by $a \mapsto \hat a$ be the homomorphism defined by $\hat a (x) = \ch(x, a)$. Similarly, let $ \hat A \ra A $ given by  $ \phi \mapsto \hat \phi$ be the homomorphism defined by $\phi(x) = \ch(x, \hat \phi)$. Recall that $\ch$ is non-degenerate.
\subsection{The action of $\mathbb{Z}_{2}$ on $\Z_{\cd}(\cc)$} 
One has that $\ct^{\delta}(X_{(a, \phi)})=X_{(-\hat{\phi}, - \hat a)}$ and $\ct^{\delta}(Z_{\ro})=Z_{\ro}$. Moreover the natural transformation $\ct^{\delta}\ct^{\delta}\ra \ct^{1}$ is given on objects by $\phi(a)$ on $X_{(a, \phi)}$ and by $\tau(\sum_{x \in A}\ro^{-1}(x))$ for $Z_{\ro}$.
\mdn
\bab  The tensor structure of $\ct^{\delta}$ is given by the following constants $c^{(a, \phi);(b ,\psi)},\; c^{(a, \phi);\ro}, c^{\ro;(a, \phi)}$ and $c^{\ro,\ro'}$.\eab 
\mdn\bll{Writing that $\ct^{\delta}\ct^{\delta}\ra \ct^{1}$ is a transformation of tensor functors one obtains that:}
\beq
c^{(a, \phi);(b ,\psi)}c^{(\widehat{\phi^{-1}}, \widehat{ a^{-1}});(\widehat{\psi^{-1}}, \widehat{ b^{-1}})}=\psi(a)\phi(b)
\eeq
In particular, for $\phi=\hat{a^{-1}}$ and $\psi=\hat{b^{-1}}$ it follows that 
\beq
(c^{(a, \hat{a^{-1}}),\;(b,\hat{b^{-1}})})^{2}=\ch(a,b)^{-2}
\eeq
\beq
c^{(a,\phi),\;\ro}c^{(\hat{\phi}^{-1}, \hat{a^{-1}}),\;\ro}=\frac{\ro\phi\hat{a^{-1}}(\lam_{A)}}{\ro(\lam_{A})\phi(a)}
\eeq
\beq
(c_{a}^{\ro,\ro'})c^{\ro,\ro'}_{a^{-1}\widehat{(\ro/\bar{\ro'})^{-1}}}(\tau\ro^{-1}(\lam_{A}))=\ch(a,a)\ro(a)(\bar{\ro'}(a))^{-1}
\eeq
\bll{One has that $c^{(a,\phi), (b, \psi)}=\phi(b)^{-1}$. Compute the other ones as a matrix relation.}
\bab 
Verifying the tensor structure relation it follows that
\beq\label{d1}
c^{(a, \phi),\;(b ,\psi)}c^{(a+b,\phi\psi),\;(c, \delta)}=c^{(b,\psi),\;(c, \delta)}c^{(a,\phi),\;(b+c,\psi\delta)}
\eeq
This relation shows that $c:A\times \hat A\ra k^{*}$ is a two cocycle.
\beq\label{d2}
\frac{c^{(a,\phi),\;\ct}c^{(b,\psi),\;\ct}}{c^{(a+b),\phi\psi,\;\ct}}=c^{(a,\phi),\;(b,\psi)}
\eeq
\beq\label{d3}
\ch(\widehat{\phi^{-1}}, \widehat{\psi^{-1}})c^{(a,\phi),\;\ro}c^{\ro\phi\widehat{a^{-1}},\;(b,\psi)}=\ch(a,b)c^{\ro, (b, \psi)}c^{(a,\phi),\;\ro\psi\widehat{b^{-1}}}
\eeq
\beq\label{d4}
c^{\ro,\;(a,\phi)}c^{\ro,\;(-\hat\phi,-\hat a)}=\frac{\ro\phi\hat{a^{-1}}(\lam_{A}}{\ro(\lam_{A})\phi(a)}
\eeq
\beq\label{d5}
c^{\ro,\;(a,\phi)}c^{\ro\phi \widehat{a^{-1}},\;(b,\psi)}=c^{(a, \phi),\;(b,\psi)}c^{\ro, \; (a+b, \phi\psi)}
\eeq
\beq\label{d6}
c^{\ro,\;(a,\phi)}c_{b}^{(\ro\phi\widehat{a^{-1}}),\;\ro'}=c^{(a,\phi),\;\ro'}c_{b}^{\ro,\; \ro'\phi\widehat{a^{-1}}}
\eeq
\beq\label{d7}
c^{(b, \frac{\widehat b\ro}{\bar{\ro'}});\;(a,\phi)}c_{b}^{\ro,\ro'}=c^{(a,\phi)\;\ro}c^{\ro' \ro'\phi\widehat{a^{-1}}}
\eeq
\beq\label{d8}
c^{(a, \frac{\widehat \ro}{\bar{\ro'}});\;\ro''}c_{b}^{\ro,\ro'}=c^{\ro;\;(a, \widehat a\ro'/\bar{\ro''})}c_{b}^{\ro',\ro''}
\eeq
\eab
\subsubsection{Inertia subgroups} It follows from above that $G_{X_{a, \phi}}=1$ if $\phi \neq -\hat{a}$. Moreover 
$G_{X_{a, -\hat{a}}}= \mathbb{Z}_2$ and  $G_{Z_{\ro}}= \mathbb{Z}_2$.
\subsubsection{Chosen isomorphisms and cocycles}
One can choose \bab $c^{\delta}_{X_{(a, - \hat a)}}=\eps_{a}\id_{a}$ \eab for a fixed squared root $\eps_{a}$ of $\ch(a, a^{-1})$ in order to obtain that the cocycle $\al_{ _{X_{(a, - \hat a)}}}(\delta, \delta)=1$. Similarly one can choose \bab $c^{\delta}_{Z_{\ro}}=\Delta_{\ro}\id_{m}$\eab  for a  fixed square root $\Delta_{\ro}$ of $(\tau(\sum_{x \in A}\ro^{-1}(x)) $. In this case one has $\al_{Z_{\ro}}(\delta, \delta)=1$. 
%\mdn  \bab With these chosen isomorphisms both cocycles $\al_{X_{(a, -\hat a)}}$ and $\al_{Z_{\ro}}$ are the trivial cocycles.\eab
\subsection{Crossed braiding on $\cz_{\cd}(\cc)$} \label{crbr}The crossed braiding on $\cz_{\cd}(\cc)$ is given as follows (see \cite{gnn}):
\beq
c_{X_{(a, \phi)}, X_{(b,\psi)}} = \psi(a) \id_{a+b } : X_{(a, \phi) }\ot X_{(b,\psi)} \ra  X_{(b, \psi)}\ot X_{(a, \phi)}
\eeq
\beq
c_{X_{(a, \phi)}, Z_{\ro} }= \ro(a) \id_{m} : X_{(a, \phi)} \ot Z_\ro \ra Z_\ro \ot X_{(a, \phi)}
\eeq
\beq
c_{Z_\ro,X_{(a, \phi)}} = \id_{m} : Z_\ro \ot X_{(a, \phi)} \ra X_{(- \hat \phi, - \hat a)}\ot Z_\ro
\eeq
\blue{\beq
c_{Z_{\ro'},\; Z_\ro} =\bigoplus_{a \in A} \ro(-a)^{-1}\ \id_{a} : Z_{\ro'} \ot Z_\ro \ra Z_\ro  \ot Z_{\ro'}
\eeq}
\subsection{Perpendicularity and cross centralizing objects}
\subsubsection{Perpendicularity}
It is easy to deduce that $X_{a, \phi}\perp X_{b, \psi}$.

Moreover  $Z_{\ro} \perp Z_{\ro'}$ for any $Z_{\ro'}$ since the composition of the homomorphisms from Equation \eqref{perp} is $\ro(a)\ro'(a)$ times identity.

On the other hand from definition of crossed braiding it follows  that $Z_{\ro}\perp X_{a, \phi}$ if and only if
$a=- \hat \phi$.
\subsubsection{Computation of $\omega$}
One has that
\bne
\item
$\omega^{1,1}_{X_{(a, \phi)}, X_{(b, \psi)}}=\phi(b)\psi(a)=\omega^{\delta,\delta}(X_{(a, \phi)}, X_{(b, \psi)})$
\item
$\omega^{1,\delta}(X_{(a, \phi)}, X_{(b, \psi)})=\chi(a, b)^{-1}\psi( \hat \phi)^{-1}=\omega^{\delta,1}(X_{a, \phi}, X_{b, \psi})$
\item
$\omega^{1,1}(X_{(a, - \hat a)}, Z_{\ro})=\omega^{\delta,\delta}(X_{(a, - \hat a)}, Z_{\ro})=\ro(a)\ch(a, a^{-1})^{\frac{1}{2}}$
\bll{change it to $\omega_{p}$}.
\bab
It seems that $Z_{\ro'}\perp_{G} Z_{\ro}$ if and only if $\ro'=\ro^{-1}$. This happens if and only if $\ch(a,b)^{2}=1$ for all $a,b \in A$.
 \eab
 \blue{ \item$\omega^{1,1}(Z_{\ro'}, Z_{\ro})=\omega^{\delta,\delta}(Z_{\ro'}, Z_{\ro})=\tau((\sum_{x \in A}\ro^{-1}(x))(\sum_{x \in A}\ro^{'-1}(x)))^{\frac{1}{2}}$}
\blue{ \item
$\omega^{1,\delta}(Z_{\ro'}, Z_{\ro})=\omega^{\delta,1}(Z_{\ro'}, Z_{\ro})=\tau((\sum_{x \in A}\ro^{-1}(x))(\sum_{x \in A}\ro^{'-1}(x)))^{\frac{1}{2}}$}
\ene

\subsubsection{Crossed centralizing} It can be checked that the following pairs of objects cross centralize each other: $Z_{\ro}\perp_{G} Z_{\ro'}$, $Z_{\ro}\perp_{G}X_{a, -\hat a}$ and  $X_{a, \phi}\perp_{G}X_{b, \psi}$ for any $a, b \in A$ and any $\phi, \psi \in \hat A$.

\br This shows that $\cs^{\perp_{G}}$ in general is not a fusion subcategory. For example, following the formulae for the tensor product $Z_{\ro}\ot Z_{\ro'}$ from \cite{gnn} it follows that this product has constituents $X_{a, \phi}$ with $\phi \neq - \hat a$.
\er

\subsubsection{Fusion rules for $\cz_{\cd}(\cc)$}\label{frules}
\beq\label{xx}
X_{(a, \phi)}\ot X_{(b,\psi)}=X_{(a+b, \phi+\psi)}
\eeq
\beq\label{xzro}
X_{(a, \phi)}\ot Z_{\ro}=Z_{\ro\phi(-\hat a)}
\eeq
\beq\label{zxro}
Z_{\ro} \ot X_{(a, \phi)}=Z_{\ro\phi (- \hat a)}
\eeq
%\blue{It should also be additive notation? $\ro+\phi - \hat a$.}
\beq\label{zzro}
Z_{\ro'} \ot Z_{\ro}=\bigoplus_{a \in A}X_{(a, \hat a \ro'/\bar \ro)}
\eeq
It follows that $X_{a, \phi}^{*}=X_{-a,\; - \phi}$ and $Z_{\ro}^{*}=Z_{\bar{\ro}}$ where $\bar{\ro}(x)=\ro(-x)$.
\subsubsection{Fusion subcategories of $\cz_{\cd}(\cc)$ stable under the action of $\mathbb{Z}_{2}$} Let 
\beqn
\ro(A):=\{\ro :A \ra k^{*}\;|\; \ro(a+b)=\ch(a,b)^{-1}\ro(a)\ro(b)\}.
\eeqn
Define an action of $A \times \hat A$ on $\ro(A)$ by $(a, \phi).\ro=\ro\phi(-\hat a)$. Clearly this action is transitive since $\ro/\ro' \in \hat A$ for any $\ro, \ro' \in A$. From Equation \eqref{xzro} it follows that this action is well defined. \mdn Thus fusion subcategories of $\cz_{\cd}(\cc)$ are parameterized by subgroups $\Gm \subset A \times \hat{A}$ and sets $M \subset \ro(A)$ that are stable under the left action of $\Gm$ on $\ro(A)$. \mdn \bab Moreover if $\ro, \ro'\in M$ then $X_{a, \widehat {a}\ro/\bar{ \ro'}}\in \Gm$. If $M$ is not empty then $\ro^{*} \in M$ and therefore $(a, \hat a)\in \Gm$. Need a classification of all these fusion categories. 
 Denote this fusion category by $\cs(\Gm, M)$. \eab \mdn It follows that $S(\Gm, M)$ is stable under the action of $\mathbb Z_2$ if and only if  $(-\hat{\phi}, -\hat{a})\in \Gm$ whenever  $(a, \phi)\in \Gm$. The support of the fusion category $\cs(\Gm, M)$ is $G=\mathbb Z_2$ unless $M$ is trivial.
 \bab Note that $\fp(S(\Gamma, M))=|\Gamma|+|M||A|$. \eab
\subsubsection{The order two automorphism on $A \times \hat A$} The braided tensor auto equivalence $\ct^{\delta}$ of $\cz(\mathrm{Vec}_{A})$ determines an order $2$ automorphism $\delta$ of $A \times \hat A$ given by
\beq
\delta((a, \phi))=(-\hat \phi, - \hat a), (a, \phi)\in A \times \hat A.
\eeq

For a subgroup $\Gm \subset A \ts \hat A$ denote by $\Gm_{0}$ the largest subgroup of $\Gm$ stable under $\delta$. Note that the subgroups stable under $\delta$ are in bijection with the subgroups of $A$. To any subgroup $B \subset A$ one associates $\Gm_{B}\subset A \ts \hat A$ given by $\Gm_{B}=\{(b, \; -\hat b)\:|\; b \in B\}$.
\subsubsection{Preinverse image} Denote by $F:\cz_{\cd}(\cc)^{G}\ra \cz_{\cd}(\cc)$ the forgetful functor. Note that in this case 
\beq
F^{-1}(\cs(\Gm, M))=\{X_{a, \eps}\;|\; (a, -\hat a)\in \Gm \} \cup \{Y_{a, b}\;|\; (a, - \hat b)\in \Gm\} \cup \{Z_{\ro, \Delta}\;|\;\ro \in M\}
\eeq
\subsection{On the center of a TY category}
\subsubsection{Simple objects of the center of $\cc=TY(A, \ch, \tau)$}
Following \cite [ Section 2]{buna} we will enumerate below the simple objects of the equivariantization $\cz_{\cc}(\cd)^{G}$ and identify them with the description given in 
\cite [ Proposition 4.1]{gnn}.
\mdn
{\bf Case 1:} If $\phi \neq - \hat a$ one has that $G_{X_{(a, \phi)}}=1$ and therefore $\pi = 1$ for the simple objects 
$S_{X_{a, \phi}, \pi}$. In this case, $S_{X_{a, \phi}, \pi}$  is just the sum of the two objects
\beqn
S_{X_{a, \phi}, 1}=X_{a, \phi} \oplus X_{-\hat{\phi}, -\hat{a}}.
\eeqn
The equivariant structure $T^{\delta}(S_{X_{a, \phi}, 1})\xra{} S_{X_{a, \phi}, 1}$ of $S_{X_{a, \phi}, 1}$  is given by  $ \phi(a)\id_{X_{a,\phi}}\oplus\id_{X_{-\hat \phi, -\hat a}}$. 
\mdn
 If $\phi=-\hat b$ then the equivariant structure of $S_{X_{a,-\hat b},\;1}=X_{a,-\hat b}\oplus X_{b, -\hat a}$ is given by  $\ch(a,b)^{-1}\id_{X_{a,-\hat b}}\oplus\id_{X_{b, -\hat a}}$.
\mdn\bab
Then $S_{X_{a,-\hat b},\;1}=Y_{b, a}$ with the notations from \cite [ Proposition 4.1]{gnn}.
\eab\mdn
{\bf Case 2:} 
If $\phi = - \hat a$ one has that $G_{X_{(a, \phi)} }=\mathbb Z_2$ and therefore $\pi \in \irr(\mathbb{Z}_{2})= \{1, \sg\}$. In this case, $S_{X_{(a, - \hat {a})}, \pi}=a$ as an object of $\cc$ with the equivariant structure given by
$\mu^{\delta}:\ct^{\delta}(S_{X_{(a, - \hat {a})}})\ra S_{X_{(a, - \hat {a})}}$
given by $\eps_{a}\pi(\delta) \id_{S_{X_{(a, - \hat {a})}}}$. They correspond to the $2n$ invertible objects $X_{a, \eps}$ from \cite [ Proposition 4.1]{gnn} where $\eps^{2}=\ch^{-1}(a, a)$. Recall we have fixed a square root $\eps_{a}$ of $\ch(a, a)^{-1}$ and chosen a fixed isomorphism $c^{\delta}_{X_{(a, \;-\hat a)}}=\eps_{a}\id_{a}$. \mdn \bab With notations from \cite [ Proposition 4.1]{gnn} one has that $S_{X_{(a, -\hat a)}, \eps}=X_{a, \eps_{a}}$ and $S_{X_{(a, -\hat a)}, \sg}=X_{a, -\eps_{a}}$.\eab
\mdn
{\bf Case 3:} 
For the simple object $Z_{\ro}$ one has $G_{Z_{\ro}}=\mathbb Z_2$ and therefore $S_{Z_{\ro}, \eps}$ has  the equivariant structure $T^{\delta}(S_{Z_{\ro}, 1} )\ra S_{Z_{\ro}, 1}$  given by $\Delta_{\ro}\id_{Z_{\ro}}$ for the fixed square root $\Delta_{\ro}$ of $(\tau(\sum_{x \in A}\ro^{-1}(x))$.  Also $S_{Z_{\ro}, \sg}=Z_{\ro, \;\Delta_{\ro}}$ has the equivariant structure  $T^{\delta}(S_{Z_{\ro}, \sg} )\ra S_{Z_{\ro}, \sg}$  given by $-\Delta_{\ro}\id_{Z_{\ro}}$.  
\mdn\bab With the notations from \cite[Proposition 4.1]{gnn} one has that $Z_{\ct, \Delta_{\ct}}=S_{Z_{\rho}, \eps}$ and $Z_{\ro, -\Delta_{\ro}}=S_{Z_{\ro}, \sg}$.\eab
%\mdn\bll{These objects correspond to the $2n$ objects indexed by $Z_{\ro, \Delta}$ where  $\Delta^{2} =\tau(\sum_{x \in A}\ro^{-1}(x))$.}
\subsection{Fusion rules for the center, i.e., $\cz_{\cd}(\cc)^{\mathbb{Z}_{2}}$} In the sequel we identify $\eps$ and $\Delta$ with the corresponding irreducible representations of $\mathbb{Z}_{2}$. Using \cite[Corollary 3.11]{buna} it follows that the fusion rules of the center are given by the following equations:
\beqn
X_{a, \eps}\ot X_{b, \eps'}=X_{a+b, \eps''}
\eeqn
where $\eps''=\eps\eps'\tau_{a+b, - \hat{ (a+b)}}^{(a ,-\hat a), (b, -\hat b)}$.
\beqn
X_{a, \eps}\ot Z_{\ro, \Delta}=Z_{\ro\hat a^{2}, 
\Delta'}\eeqn $\Delta'=\eps\Delta\tau^{(a, -\hat a), \ro}_{\ro\hat a^{2}}$.

If $\ro/ \overline{\ro'}\neq  \hat a ^{-2}$ for all $a \in A$ then 
\beqn
 Z_{\ro',\; \Delta'}\ot Z_{\ro, \;\Delta}=\oplus_{a \in A/\sim_{\ro, \ro'}}Y_{a, \hat a \ro/ \bar \ro'}
\eeqn
where they equivalence relation $\sim_{\ro, \ro'}$ is given by.

If $\ro/ \overline{\ro'}=\hat a ^{-2}$ for some $a$ then
\beqn
 Z_{\ro',\;\Delta'}\ot Z_{\ro, \Delta}=X_{a, \eps}\oplus (
\bigoplus_{b \in A- \{a\}/\sim_{\ro, \ro'}}Y_{a, \hat a \ro/ \bar \ro'})
\eeqn
where $\eps=\Delta\Delta'\tau^{\ro', \ro}_{a, \hat a}$.

\beqn
Y_{a, b}\ot Z_{\ro, \Delta}=Z_{(a, -\hat b).\ro,\;\Delta_{0}}\oplus Z_{(b, -\hat a).\ro, \;\Delta_{1}}
\eeqn
Note that $(a, -\hat b).\ro=(b, -\hat a).\ro=\ro \hat{a}^{-1} \hat{b}^{-1}$.

\beqn
Y_{a, b}\ot Y_{c,d}=Y_{a+d, b+c}\oplus Y_{a+c, b+d}
\eeqn
\subsection{Dual objects}It follows that 
\beqn
X_{a, \eps}^{*}=X_{-a, \eps^{-1}[\tau^{(a, -\hat a), (-a, \hat a)}_{(0, \hat 0)}]^{-1}}
\eeqn
\beqn
Z_{\ro, \Delta}^{*}=Z_{\bar \ro,\; (\Delta\tau^{\ro, \bar \ro}_{0, \hat 0})^{-1}}
\eeqn
\beqn
Y_{a, b}^{*}=Y_{-a, -b}
\eeqn
\subsection{On the projective representations $\tau^{Y,Z}_{U}$}
One has \beqn \tau_{X_{b+b', \;-\widehat{b+b'}}}^{X_{b, \;- \hat b}, X_{b', \;-\hat b'}}:=\Hom_{\cc}(X_{b+b', \;-\widehat{b+b'}},\; X_{b, \;- \hat b}\ot X_{b', \;-\hat b'})\eeqn is a one dimensional $\mathbb Z_{2}$-representation \bab which is trivial since {\bf $(T^{\delta})_{2}^{X_{b, -\hat b},\; X_{b'\; - \hat b^{'-1}}}$ is trivial}.\eab 

Since $Z_{\ro} \ot X_{a, - \hat a}=Z_{\ro}$ one has that 
\beq
\tau^{Z_{\ro}, \;X_{a, -\hat a}}_{Z_{\ro}}=\Hom_{\cc}(Z_{\ro},\; Z_{\ro} \ot X_{a, -\hat a})
\eeq
and 
\beq
\tau^{X_{a, -\hat a},\; Z_{\ro}}_{Z_{\ro}}=\Hom_{\cc}(Z_{\ro},X_{a, -\hat a}\ot Z_{\ro}\;)
\eeq
are one dimensional representations of $\mathbb{Z}_{2}$.

On the other hand  since 
\beq
\tau^{Z_{\ro},\; Z_{\ro'}}_{X_{a, \phi}}=\Hom_{\cc}(X_{a, \phi},\;Z_{\ro}\ot Z_{\ro'}\;)\eeq
it follows that this is also a one dimensional representation of $\mathbb{Z}_{2}$ if $(a, \phi).\ro=\ro'$.

It follows that all representations are one dimensional and we denote them by $\tau^{(a, \phi), (b, \psi)}_{a+b ,\phi+\psi}$ respectively...
\mdn \bll{These objects are necessary to deduce the fusion rules. Maybe they can be written directly by knowing the fusion rules of the relative center.}
\subsection{The $S$-matrix of $\cz(\cc)$}\label{smat} Following \cite [ Section 4.3 ]{gnn} one has that
\beq\label{aepsaeps}
S_{X_{a, \; \eps},\; X_{a', \; \eps'}}=\ch(a, a')^{2}
\eeq
\beq\label{aepszro}
S_{X_{a, \; \eps},\;Z_{\ro, \Delta}}=\eps\sqrt{n}\ro(a)
\eeq
\beq\label{yabzro}
S_{Y_{a, b}, \; Z_{\ro, \Delta}}=0 
\eeq
\beq\label{aepsyab}
S_{X_{a, \; \eps},\;Y_{a, b}}=2\ch(a, b+c)
\eeq
\beq\label{yabycd}
S_{Y_{a, b},\; Y_{c, d}}=2(\ch(a, d)\ch(b, c)+\ch(a, c)\ch(b, d))
\eeq
\beq\label{zrozro'}
S_{Z_{\ro, \Delta}, \; Z_{\ro', \Delta'}}=\frac{1}{\Delta \Delta'}\sum_{a \in A}\ch(a, a)^{2}\ro(a)\ro'(a)
\eeq
Note that $Z_{\ro, \Delta}$ cannot centralize any $Y_{a, b}$.\mdn 
\bll{ Recall that the $S$-matrix is obtained by using fusion rules and Verlinde formula.}
\subsubsection{Crossed centralizing objects}
From the $S$-matrix one can deduce that \bne
\item
$
X_{(a, -\hat a)}\perp X_{(b, -\hat b)}$ if and only if $\ch(a,b)^{2}=1$.
\item $
X_{(a, -\hat a)}\perp Z_{\ro}$ if and only if $\ro(a)^{2}=\ch(a, a)$.
\item $X_{(a,-\hat b)}$ with $b\neq a$ never cross centralizes $Z_{\ro}$.
\item $X_{(a,-\hat a)}\perp X_{(b, -\hat c)}$ with $b\neq c$  if and only if $\ch(a,b+c)=1$.
\item $X_{(a,-\hat b)}\perp X_{(c, -\hat d)}$ with $a\neq b$ and $c \neq d$  if and only if $\ch(a,d)\ch(b,c)=\ch(a,c)\ch(b,d)=1$.
\item $Z_{\ro}\perp Z_{\ro'}$ fi and only if 
\beq
(\sum_{a \in A}\ch(a,a)^{2}\ro(a)\ro'(a))^{2}=n^{3}(\sum_{a\in A}\ro(a)^{-1})(\sum_{a\in A}\ro'(a)^{-1})
\eeq
\ene
\subsection{Fusion subcategories of the center of Tambara-Yamagami categories}
Using the above description for fusion subcategories of an equivariantization it follows that one has four types of fusion subcategories of the center. \mdn   
\bab write in each case what are the irreducible objects of the fusion subcategory\eab
\mdn
\bll{The centralizer is obtained by mixing the $S$-matrix with the results written here.}
\mdn
{\bf Case 1:}
$\cd=\cc(1,1,  \cs(\Gm, M), \lam)$. In this case since the support of $\cs(\Gm, M)$ is trivial it follows that $M$ is the empty set. Moreover since $H=1$ it follows by Remark \ref{trivh} that $\lam$ is also trivial. Then since $K=1$ it also follows that $\tilde \lam$ is trivial. Applying Theorem \ref{switch} now one has that
\beq
\cc(1,1,  \cs(\Gm), 1)'=\cc(1,1,  \cs(\Gm^\perp), 1)
\eeq
where $\Gm^{\perp}:=\{ (b, \psi)\;|\;\phi(b)\psi(a)=1 \text{\;for all \;} (a, \phi)\in \Gm\}$.
\mdn
\bab Note that these are fusion subcategories of $\cz(\cd)$ whose centralizer is also contained in $\cz(\cd)$.\eab % was also described in \cite{nnw}.
\mdn
{\bf Case 2:}
$\cd=\cc(\mathbb{Z}_{2},1, \cs(\Gm, M),  \lam)$. Since $H=\mathbb Z_2\subset \bigcap_{Y \in \irr(\cs(\Gm, M))}G_{Y}$ it follows that any $(a, \psi) \in \Gm$ is of the type $(a, -\hat{a})$. Thus 
\beqn
\Gm=\Gm_{B}=\{(b , -\hat b) \;|\; b \in B\}
\eeqn
for some subgroup $B \subset A$.
Also in this case since the support of $\cs(\Gm, M)$ is trivial it follows that $M$ is the empty set. \mdn Note that in this case is defined by $\lam: B\times $ \mdn  \blue{It follows that $\lam([X_{b, -\hat b}], \;-)$ is a character of $\mathbb Z_{2}$ and 
and $\lam(bb' -)=\lam(b, -)\lam(b' - )$.} Therefore
\beq
\irr(\cc(\mathbb Z_2,1,  \cs(\Gm(B)), \lam))=\{X_{b, \eps} \in \cz(\cc)\;|\;\lam(b ,\delta)=\eps/\ch(b, b)^{\frac{1}{2}}\}.
\eeq
Using the above $S$-matrix from subsection \ref{smat} it follows that
\beq
\cd'=\{X_{c, \eps}\;|\;c \in B^{\perp}\}\cup \{Y_{b, c}\;|\;b+c \in B^{\perp}\}\cup \{Z_{\ro, \Delta}\;|\; \ro(b)=\lam(b, \delta)^{-1}\ch(b, b)^{\frac{1}{2}}\}
\eeq
where $B^{\perp}$ is taken with respect to $\ch$. 
\mdn Thus applying Theorem \ref{switch} one has that
\beq
\cc(\mathbb Z_2,1,  \cs(\Gm(B)), \lam)'=\cc(1,\mathbb Z_2,  \cs(\Gm', M), 1)
\eeq
where \blue{$\Gm'=\{(c, -\hat c)\;|\; \ch(b, c)^{2}=1, b \in B\}\cup \{(d, -\hat c)\;|\;d+c \in B^{\perp}\}$ and $M=\{\ro\in \ro(A)\;|\; \ro(b)=\lam(X_{b , \;-\hat b}, \delta)^{-1}\ch(b, b)^{\frac{1}{2}}, \; b \in B\}.$ Here the orthogonal subgroup $B^{\perp}$ is taken with respect to $\ch$.}
\mdn
{\bf Case 3:}
$\cd=\cc(1,\mathbb{Z}_{2}, \cs(\Gm, M),  1)$. Here for the centralizer $\cd'$ one has that $H=\mathbb{Z}_{2}$ and $K=1$. It follows that $\cc(1,\mathbb{Z}_{2}, \cs(\Gm, M),  \lam)'$ consists only on invertible objects $X_{a, \eps} \in \cz(\cc)$. Thus this case is dual to the previous case and one has that
\beq
\cc(1, \mathbb Z_2,  \cs(\Gm, M), 1)'=\cc(\mathbb Z_2, 1,  \cs(\Gm'_{B}), \tilde \lam)
\eeq
where $B=\{\Gm_{0}, b+c\:|\; (b, -\hat c) \in \Gm\}^{\perp}$ and $\tilde \lam(b, \delta)=\ro^{-1}(b)\ch(b,b)^{1/2}$ for any $\ro \in M$. Recall that $\Gm_{0}$ denotes the largest subgroup of $\Gm$ stable under the automorphism $\delta$.
\mdn
{\bf Case 4:}
 $\cd=\cc(\mathbb{Z}_{2},\mathbb{Z}_{2}, \cs(\Gm, M),  \lam)$.  Since $H=\mathbb Z_2\subset \bigcap_{Y \in \irr(\cs(\Gm, M))}G_{Y}$ it follows as in case 2 that any $(a, \psi) \in \Gm$ is of the type $(a, -\hat{a})$. Thus 
\beqn
\Gm=\Gm_{B}=\{(b , -\hat b) \;|\; b \in B\}
\eeqn
for some subgroup $B \subset A$. Thus one can write that
\beq
\cc(\mathbb Z_2, \mathbb Z_2,  \cs(\Gm_{B}, M), \lam)'=\cc(\mathbb Z_2, \mathbb Z_2,  \cs(\Gm'_{B'}, M'), \tilde \lam)
\eeq
for another subgroup $B' \subset A$ and another subset $M' \subset \ro(A)$.
\onh\; one can write that
\beqarn
\cd&=& \{X_{c, \eps}\;|\;(c, -\hat c) \in \Gm, \lam([X_{c, -\hat c}], \delta)=\eps_{c}(\delta) \}\cup  \{Z_{\ro, \Delta}\;|\; \lam([Z_{\ro}], \delta)=\eps_{\Delta}(\delta)\}
\eeqarn
Using the $S$-matrix from subsection \ref{smat} it follows that
 \beqarn
 \cd' &=& \{X(c, \eps)\;|c,\in B'\;|\; \eps\ro(c)=1\;\text{for all}\; \ro \in M \}\cup \\&\cup&  \{Z_{\ro, \Delta}\;|\; \lam([Z_{\ro}], \delta)=\eps_{\Delta}(\delta)\}
 \eeqarn
where $\B'=\{a \in A\;|\; \ch(a,c)^{2}=1 \text{\; for all \;} c \in B\; \text {and } \ro(a)^{2}=\ch(a,a)\text{\;for all\;} \ro \in M\}$, and $M'=\{\ro\in \ro(A)\;|\; \ro(c)\eps=1 \;\text{for all} \; X_{c, \eps} \in \cd \;\text{and}\;\sum_{a \in A}\ch(a,a)^{2}\ro(a)\ro'(a)=\Delta\Delta'\}$. Equation \eqref{f1tilde} implies that \beqn
\tilde \lam([X_{(a, - \hat a)}], \delta)=\omega^{1,1}([Z_{\ro}], [X_{(a, - \hat a)}]).
\eeqn
for any $\ro \in M$. Moreover
\beqn
\tilde \lam([Z_{\ro'}], \delta)=\lam^{-1}([Z_{\ro}], \delta)\omega^{1,1}([Z_{\ro}],  [Z_{\ro'}]).
\eeqn
for any $\ro \in M$ and any $\ro' \in M'$.
\mdn
\bab Decide when the center of TY is a prime category. Decide which fusion subcategories are nondegenrate.\eab
\subsection{Example $A=\mathbb{Z}_{p}$ and $\ch$ coming from a quadratic form}
\mdn
\bll{need to count how many $\rho$'s are.}
\mdn
\bll{verify the conjecture in this case, I basically know delta on $K$, it should give a projective representation since the M\"{u}ger centraliser is not empty. All the others are verifying this property? need a dependence of $\omega$ and $cocycle$.}
\mdn{If I can make it one then somehow all the troubles disappear.}
\bibliographystyle{amsalpha}

\np
{\bf Questions:} 
 \bne
 \item Given a fusion subcategory $\cs$ is the  
 abelian subcategory $\cs^{\perp_{G}}$ also a fusion subcategory? It seems that the answer is NO from next Example.
 \item Does  $\mathcal{T}$ consists of all objects of $\cs^{\perp_{G}}$ whose degrees belong to $H$?
 \item
 Suppose that $X \perp Y$, i.e. $d(1,1)_{X, Y}d_{Y, X}(1, 1)=\omega(X, Y)\id_{X \ot Y}$. Does it follow that $\omega(X, Y)$ is a root of unity?
 \item
 Suppose that $X \perp_{G} Y$ with
 \beq
 d_{X, Y}^{m,n}d_{Y, X}^{n,m}=\omega(\ct^{m}([X]),  \ct^{n}([Y]))\id_{X \ot Y}
 \eeq
 Given $\pi  \in \Irr_{\al_{X}}(G_{X})$ define 
 \beq
 \delta(ngn^{-1})=\frac{\pi(1)}{\pi(mhm^{-1})}\omega(\ct^{m}([X]),  \ct^{n}([Y]))
 \eeq
 Is it $\delta$ well defined and a 
 $\al_{Y}$-projective representation of $G_{Y}$?
 \ene
 \bne
\item Make a section with projective representations and their characters, multiplicity?
\item what kind of multiplicity is used in the fusion rules?
\item lifting a projective representation to a linear representation the character is the same for the lifting? so we can get directly the inequality.
\item put a an example with NNW, how to get multiplicative in one side.
\item maybe we should put $K$ first to obtain the same notation as they have.
\ene
\np
\bne
\item use the description of the simple objects from g-funtors as $S(g, M)$. 
\item give a new condition when two such objects cross centralize.
\item
there are gonna be new $d_{M,N}d_{N, M}=\id_{M\otimes N}$
\item use the other property of $\lam$ to get the relations between omegas'
\item One has that $X_{c, -\hat b}\otimes X_{b, -\hat c}=X_{a, -\hat a}$ where $a=b+c$. None of them does cross centralize $Z_{\ro}$ but their product it does.
\item for the beginning I can suppose that the action is strict, i. $\ct^{g}\ct^{h}=\ct^{gh}$.
\item take $\delta$ from Remark 4.15, is it a projective representation on $N(g)$?
\item use corollary 3.6 from \cite{dgno} the rank of $\cd$ is the number of components of $\cd'$.
\item change the function $d_{Y}(g,h)$ to another letter.
\item sync with fsc.tex, put the braided part from smc over the braided part form fsc
\item Here one needs the compatibility with the unit $c_{\1, X}=\id_{X}$. Maybe it results from the other? $\tilde{c}$ is a braiding and it should be the identity
\item square versus squared
\item from the first relation the crossed centralizer relation I don't think it is transitive
\item The question is if all crossed centralizing objects of a given object form a fusion subcategory?
\item two objects cocentralize each other if and only if there are simple objects seating over them which M\"{u}ger centralize each other.
\item From the S-matrix compute explixitly which object cross centralize; then try to derive the formulae for the categories.
\item compute an explict example from their papers to deduce what are the stable fusion subcategories. $A=\mathbb Z_{p}$
\item the simple objects of $\cz(\cc)$ are given in GNN by the equivariant structure.
\item the second cohomology is trivial.
\item isomorphism chosen such that all cocycles are exactly 1.
\item Finish the last two cases so they do not depend on the $\eps$ from the initial category.
\item $Z_{\ro}\perp_{G}Z_{\ro'}$ does not seem to give the same thing as computing directly the $S$-entry. $Z_{\ro, \Delta}$ centralizes $Z_{\ro', \Delta'}$ It possibly explains why the entries in modulus is less or equal $n$ and equality might hold only in the conditions from cross centralizer between $Z_{\ro}$ and $Z_{\ro'}$. 
\item
Try also generalized Tambara-Yamagami categories.
\item Establish completely the detramnination $\eps \mapsto \eps_{a}$ and $\Delta \mapsto \eps_{\Delta}$.
\item Put the other remark with $b=1$ in the formula for $\tilde \lam$.
\item
$|G|=p^{n}$ then $\cc^{G}$ is prime for $p=2$.
\bpf
Suppose that
$\cc(H, K, G)'=\cc(K, H, \mathcal{T}, \tilde \lam)$. It follows that $HK=G$ and $H \cap K=1$
\epf
\item If $G$ is abelian then $\cc^{G}$ is the Deligne product of some fusion categories.
\item $G/H$ acts on $\cs$ since $H \subset G_{V}$ for all $V \in \cs$. Is the equivariantization $\cc(H,K, \cs, \lam)$?
\mdn They have the same Perron Frobenius dimension and the simple objects look very similar. $Stab_{G/H}(V)=G_{V}/H$
\item Cosets of $\cc(H, K, \cs, \lam)$. Apply the result of Drinfeld.
\item In Drinfeld doubles of groups
\beq
(S_{a, \gm})'=\cc(C_{G}(N(a)), N(a), \lam)
\eeq
The conjecture is that the number of components of $\cc(C_{G}(N(a)), N(a), 1)$ equals the number of conjugacy classes of $G$.
\ene
\subsection{Some canonical objects of the center $\rep(D(A))$.}
Suppose that $\cc$ admits a (quasi-)functor, i.e. is the category of representations of a (quasi-)Hopf algebra, $A$. Then $A \in \cz(\cc)$ via the half braiding given by the inverse of:
\beq A \ot M \ra M\ot A,\;\;
a \ot m \mapsto a_{1}Sa_{3}m\ot a_{2}
\eeq

Then  by restriction this map induces isomorphisms $\phi_{g, M}:R_{g} \ot M \ra M \ot R_{h^{-1}gh}$.
\bp
There is  a functor \bb{an embedding} of categories
\beq
\rep(G_{a})\hookrightarrow \rep(D(A))
\eeq
given by $V \mapsto T(a, V):=\oplus_{x \in G/G_{a}}R_{xax^{-1}}\ot V$
 where the commuting structure 
 \beq
 M \ot T(a, V) \ra T(a, V) \ot M
 \eeq
 is given on components by 
 \beq M\ot R_{xax^{-1}}\ot V \xra{\phi_{M, g}\ot \pi_{V}(h)} R_{zaz^{-1}}\ot M \ot V
 \eeq
 \ep
 \bpf
 It is easy to verify that this element is central.
 \epf 
\bp\beq
T(a, V)\ot T(b, W)=\oplus_{x \in D}T(\;^{x}ab, m_{\;^{x}a,b}(\;^{x}V, W))
\eeq
\ep
\bl
The regular characters of each summand is $T(a, \mtr{reg}_{G_{a}})$.
\el
%\ed
\br
Relation \eqref{po} shows that the first Property of $\lam_{\cd}$ does not depend on the collection of isomorphisms $c$. Indeed, if $\lam_{\cd}'$ denotes the corresponding function $\lam$ for another chosen set of isomorphisms $c'$ then one has that
\bqarn
\lam_{\cd}'([Y], h)=\frac{{(f_Y}^{-1}\pi)(h)}{\pi(1)}={f_Y}^{-1}(h)\lam_{\cd}([Y], h)
\eqarn
where $f_{Y}:=f_Y^{c,c'}$.
Then 
\bqarn
\lam'_{\cd}([Y], g)\lam_{\cd}'([Y], h)&=&{f_Y}^{-1}(g)\lam_{\cd}([Y], g){f_Y}^{-1}(h)\lam_{\cd}([Y], h)\\ &=&\frac{\lam_{\cd}([Y], gh)}{\al_Y(g,h)}f_Y^{-1}(g)f_Y^{-1}(h)=\frac{\lam_{\cd}'([Y], gh)}{\al'_Y(g,h)}
\eqarn
Using Equation \eqref{yz} the same relation \eqref{po}  gives also the following dependency of $\tau_U^{Y, Z}$ on the change of the isomorphisms from $c$ to $c'$:
\bq
{(\tau_{c})}_U^{Y, Z}(h)=\frac{f_{Y}(h)f_{Z}(h)}{f_{U}(h)}{(\tau_{c'})}_U^{Y, Z}(h)
\eq
\er

\bll{\bl
Suppose that $X\perp_{G}Y$. Then one can modify a good chosen systems of isomorphisms to another good one such that in this new system one has $\omega_{X, Y}^{1, p}=1$ for any $p \in G$.
\el
}
\subsubsection{Expressing the conjugate projective representation} Let $\al_{X}$ be a cocycle on $G_{X}$.
Note that inside the twisted group algebra $k_{\al_{X}}G_{X}$  one has \beqn \ovr{php^{-1}}=\al_{X}^{-1}(ph, p^{-1})\al_{X}^{-1}(p,h)\al_{X}(p, p^{-1})\ovr{p}\ovr{h}\ovr{p^{-1}}.\eeqn
Denote by \beq\omega_{\al_{X}}(p,h):=\frac{\al_{X}(p, p^{-1})}{\al_{X}(ph, p^{-1})\al_{X}(p,h)}\eeq
It follows that
\beq\label{conjpp}
\pi(php^{-1})=\omega_{\al_{X}}(p,h)\pi(p)\pi(h)\pi(p^{-1})
\eeq
\bll{for any $p \in G_{Y}$.}

\subsection{An example from group extensions}
Suppose that $A=kG$ and $B=kN$ for some normal subgroup $N$ of $G$. Thus we have the following exact sequence
\beq
1 \ra N \ra G \xra{\pi} F \ra 1
\eeq
where $F:=G/N$. Then it is well known that $F$ acts on $\Rep(kN)$ and $\Rep(kG)=\Rep(kN)^{F}$.

Consider $s:F \ra G$ a section of $\pi$ with $s(g^{-1})=s(g)^{-1}$ and $s(1)=1$. The action of $F$ on $N$ is given by $f.a=s(f)^{-1}as(f)$ for all $f \in F$ and all $a \in N$. The two cocycle $\sg$ is given by $\sg(f, g)=s(f)s(g)s(fg)^{-1}$ for all $f,g \in F$. With these notations the isomorphism $G \xra{\phi} N \#_{\sg} G/N$ is given by $\phi(g)=gs(g)^{-1}\#_{\sg} \ovr{g}$. Its inverse is given by $\pi^{-1}(a \#_{\sg}f)=as(f)$ for all $f \in F$ and all $a \in N$.

The action of $F$ on $\rep(N)$ is given as following. If $\ovr{g} \in F$ then $T^{\ovr{g}}(Y)=ks(g) \ot Y$ as vector space and $a.(s(g) \ot y)=s(g)\ot (s(g)^{-1}as(g)).y$ for all $a \in N$ and all $y \in Y$.

Suppose that $G_{Y}$ is the stabilizer in $G$ of a simple $kN$-module $Y$. Then $F_{Y}:=\pi(G_{Y})$. Then the induction functor from $G_{Y}$ to $G$ gives the bijection between simple $G$-modules containing $Y$ and projective representations of $F_{Y}$ with factor set $\tilde \al_{Y}$.

The structure of $G_{Y}$-module is given by
\beq\label{rel}
g.(y \ot v)=gs(g)^{-1}c_{Y}^{\ovr{g}}(y) \ot \ovr{g}.v
\eeq

Define $\tilde c^{g}:Y \ra Y$ given by
\beq
\tilde c^{g}(y)=gs(g)^{-1}c_{Y}^{\ovr{g}}(y)
\eeq
for all $g \in G_{Y}$ and all $y \in Y$. It follows that $
g.(y \ot v)=\tilde c^{\ovr{g}}(y) \ot \ovr{g}.v$. Note that $\tilde c^{g}_{Y}$ is invertible with inverse $y \mapsto s^{-1}(g)g(c_{Y}^{g})^{-1}(y)$
\subsection{The datum associated to a fusion subcategory of $\rep(G)$. }One has that $H_{\rep(G/M)}=MN/N$. Moreover 
$S_{Y, \pi} \in \rep(G/M)$ if and only if $M \subseteq \ker_{G}(\Ind_{G_Y}^{G}(Y\ot V_{\pi}))$. Applying Lemma \ref{krnind} it follows that $S_{Y, \pi} \in \rep(G/M)$ if and only if $M \subseteq \ker_{G_{Y}}(Y\ot V_{\pi})$. Note that by the definition of the $G_{Y}$-module $Y \ot \pi$ one has $M \subseteq \ker_{G_{Y}}(Y\ot V_{\pi})$ if and only if $\tilde c^{m}_{Y}=\omega(Y, m)\id_{Y}$ and $\pi(m)=\omega(Y, m)^{-1}\id_{V_{\pi}}$ for some root of unity $\omega(Y, m)\in k^{*}$. This define 
\beqarn
\cs_{\rep(G/M)}&=&\{Y \in \rep(N), M \subseteq \core_{G}(G_{Y})\;|\;\tilde c^{gmg^{-1}}_{Y}=\omega(Y, gmg^{-1})\id_{Y} \;\;\\&& \text{for\; all \;}m \in M, g \in G_{Y}\}\eeqarn
Alternatively this can be written as
\beq
\cs_{\rep(G/M)}=\{Y \in \rep(N)\;|\;\tilde c^{gmg\inv}_{Y}\ot \pi(gmg\inv)=\id_{Y} \ot \id_{V}  \;\; for\; all \;m \in M\}\eeq
Morever 
\beq
\lam:\gr(\cs_{\rep(G/M)})\times MN/N\ra k^{*}
\eeq
is defined by $\lam([Y], \ovr{mn}):=\omega\inv(Y, m)$.
\subsection{On the fusion subcategory associated to a datum}
Suppose there is a datum $(\cs, H, \lam)$ with $H=\ovr{H}/N$ and $\lam: \irr(\cs) \times H \ra k^{*}$. Then
\beqn
	\cc(\cs, H, \lam)=\{ \Ind_{G_{Y}}^{G}(Y\ot V_{\pi})\;|\; Y \in \irr(\cs), \pi \in \irr_{\al_{Y}}(G_{Y}/N)\; \pi|_{H}=\lam([Y], \;-)\}
\eeqn
his the subgroup $M$ associated to $\cc(\cs, H, \lam)$ is given by the intersection of kernels of all induced modules $\Ind_{G_{Y}}^{G}(Y\ot V_{\pi})$. Applying again Lemma \ref{krnind} it follows that 
\beqarn
M &=&\bigcap_{(Y, \pi) \in \irr(\cs) \times \irr_{\al_{Y}}(G_{Y}/N)\; \pi|_{H}=\lam(Y, \;-)}\core_{G_{Y}}(Y\ot V_{\pi})\\&=&\bigcap_{Y, \pi }\{m \in G_{Y}\:|\tilde c^{gmg^{-1}}_{Y}\ot \pi(gmg^{-1})=\id_{Y} \ot \id_{V_{\pi}} \text{for\; all} \;g \in G_{Y}\}
\eeqarn

\red{\bnex
Suppose that $A$ is a semisimple Hopf algebra that fits to a cocentral extension 
\beq
k \ra kG \ra A  \xra{\pi} kF \ra k
\eeq
with $G$ and $F$ finite groups. Then $H:=A^{*\cop}\bwt kG$ is a normal Hopf subalgebra of $D(A)$ \cite{bnda} fitting the cocentral extension
\beq\label{c22}
k \ra H \ra D(A ) \xra{\tilde \pi} kF \ra k.
\eeq
It follows that $\rep(H)$ is a $F$-crossed braided fusion category and $\rep(D(A))=\rep(H)^{F}$. Indeed note that $k^{F} $ is a central Hopf subalgebra of $A^{*}$ via the inclusion $\pi^{*}$ and therefore a central Hopf subalgebra of $D(A)^{*}$ via $\tilde \pi^{*}$. Moreover following \cite{gnn} since $\rep(D(A))$ is braided it follows that $\rep(H)$ is a $F$-crossed braided fusion category. Note that the case $A=k^{L}$ for a finite group $L$ was considered in \cite{mns}.
\enex}

\begin{definition}
We say that $X\perp Y$ if $\deg(X) \in G_{Y}$,  $\deg(Y) \in G_{X}$ and the following composition
\beq\label{perp}
X\ot Y \xra{c_{X, Y}} \ct^g(Y)\ot X\xra{c^g_Y\ot \id} Y\ot X \xra{c_{Y, X}} \ct^h(X)\ot Y\xra{c^{h}_{X} \ot X} X\ot Y
\eeq
 is a scalar multiple $\omega([X], [Y])$ of the identity morphism.  
 \end{definition}
%\ed
\section{\blue{Examples and applications}}
%\textcolor[rgb]{1.00,0.00,0.00}{

\subsection{Cocentral extensions of Hopf algebras}
%{\bf Check that this gives an action of $G$ on the fusion $\Rep(A)$.}
Let $A$ be a semisimple Hopf algebra and $G$ be a finite group. Suppose that $H:=A\;^{\tau}\#_{\sg}kG$ is a cocentral extension of $A$ by $kG$. Therefore $H$ fits into a short exact sequence of Hopf algebras.
$$k \ra A \ra H\ra kG \ra k.$$ \\ In particular, the cocycles $\sg: kG\ot kG \ra U(A)$ and $\tau:kG \ra A\ot A$ verify the following compatibility conditions:

\bn{enumerate}
\item $\sg(1,g)=\sg(g,1)=1$
\item $(\eps_A\ot \id )(\tau(g))=1=(\id \ot \eps_A)(\tau(g))$
\item  \lb{cocyl} $[h.\sg(l,m)]\sg(h,lm)=\sg(h,l)\sg(hl,m)$
\item  \lb{twm} $h.(l.a)=\sg(h,l)[(hl).a]\sg^{-1}(h,l)$
\item $(\tau(g)_j)_1\tau(g)_l\ot (\tau(g)_j)_2\tau(g)^l\ot \tau(g)^j=\tau(g)_j\ot (\tau(g)^j)_1\tau(g)_r\ot  (\tau(g)^j)_2\tau(g)^r$

This means that $\tau(g)$ is a twist for $A^{op}$ or equivalently $\tau(g)^{-1}$ is a twist for $A$.
\item \lb{A} $\Delta_A(g.a) = \tau (g)(g.a_1 \ot g.a_2)\tau(g)^{-1}$ \\ (This is condition $(A)$ from \cite{AD})
\item \lb{D} $\Delta(\sg(g,h))\tau(gh)=\tau(g)(g.\tau(h)_p\ot g.\tau(h)^p)(\sg(g,h)\ot \sg(g,h))$
(This is condition $(D)$ from \cite{AD})  \end{enumerate}
Note that conditions $(B)$ and $(C)$ from \cite{AD} are automatically satisfied.

We used the notations $\tau(g)=\tau(g)_r\ot \tau(g)^r$ from \cite{AD}.
\md
\blue{Note that if $\tau$ trivial then $\sg(g,h)$ is a group like element.}
The multiplication on $A$ is given by
\beq
(a\#_{\sg}h)(b\#_{\sg}l)=a(h.b)\sg(h,l)\#_{\sg}hl
\eeq

and the multiplication on $A$ is given by:
\beq
\Delta(a\#g)=(a_1\tau(g)_j\#_{\sg}g)\ot (a_2\tau(g)^j\#_{\sg}g)
\eeq
Moreover the antipode is given by (see Formula (2.25) from \cite{AD}) $$S(a\#_{\sg}g)=g^{-1}S(a)=(g^{-1}.S(a))\#_{\sg}g^{-1}.$$
\bb{ According to \cite{AD} one may suppose that
\bq
\sum \tau(g)_iS(\tau(g)^i)=1=\sum S(\tau(g)_i)\tau(g)^i
\eq for all $g \in G$.
}
One may suppose that for all $g \in G$ one has $\sg(g,g^{-1})=1$ and therefore $\bar{g}^{-1}=\bar{g^{-1}}$.

\bl\lb{sgg} Assume that for all $g \in G$ one has $\sg(g,g^{-1})=1$. This implies that $$(gh)^{-1}.\sg(g, h)=\sg^{-1}(h^{-1}, g^{-1}).$$ Thus $(gh)^{-1}.\sg^{-1}(g, h)=\sg(h^{-1},g^{-1})$.
\el
\bpf
Indeed one has that

$$
(\bar{g}\bar{h})^{-1}= (\bar{h})^{-1} (\bar{g})^{-1}=\bar{h^{-1}}\bar{g^{-1}}=\sg(h^{-1}, g^{-1})\#\bar{h^{-1}}\bar{g^{-1}}.
$$

On the other hand $$(\bar{g}\bar{h})^{-1}=(\sg(g,h)\#gh)^{-1}=\bar{gh}^{-1}\sg^{-1}(g, h)=((gh)^{-1}.\sg^{-1}(g,h)))\#\bar{gh}^{-1}$$
\epf
It is known from \cite{natale-repg} that $G$ acts by tensor autoequivalences on $\rep(H)$. For an element $a \in A$ and $V\in A$-mod let $a|_V$ be the operator on $V$ given by left multiplication by $a$, namely $v \mapsto a.v$.
\bn{enumerate}
\item
The action of $G$ is given as following: $\ct^g(V)=V$ as vector spaces with the action of $A$ given by $$av:=(g^{-1}.a )v$$
\item
Moreover the tensor structure $(\ct^2)^{g,h}_V:\ct^g(\ct^h(V))\ra \ct^{gh}(V)$ of $\ct$ is given by $$(\ct_2)^{g,h}_V=\sg^{-1}(h^{-1}, g^{-1})|_V.$$
\item
The tensor structure of $\ct^g$ is given by $$(\ct^g)_2^{V,W}:\ct^g(V \ot W)\ra \ct^g(V)\ot \ct^g(W)$$ is given by $\tau^{-1}(g^{-1})|_{V\ot W}$.
\item
Note that $\ct^e = \id_\C$ and that all $\ct_0$, $\ct^{g, e}_2$, $\ct^{e, g}_2$ are identities.
\end{enumerate}

In order to verify that this defines an action of $G$ on $\rep(A)$ note the following:
\bn{enumerate}
\item
$\ct^g(V)$ is an $A$-module since $g.(ab)=(g.a)(g.b)$.
\item The functor $\ct$ is a tensor functor if the compatibility condition \ref{ro-2} of the $G$-actionis satisfied. This is equivalent to the cocycle condition \ref{cocyl} for $\sg$ and the triple $\{h,l,m\}=\{c^{-1}, b^{-1}, a^{-1}\}$.
\item
The fact that $(\ct_2)^{g,h}_V:\ct^g(\ct^h(V))\ra \ct^{gh}(V)$ is a morphism of $A$-modules is equivalent to $\sg^{-1}(h^{-1},g^{-1})[h^{-1}.(g^{-1}.a)]=(h^{-1}g^{-1}.a)\sg^{-1}(h^{-1},g^{-1})$ which is the twisted module condition \ref{twm}. %It is true for $A$ commutative.}
\item
The tensor structure $$(\ct^g)_2^{V,W}:\ct^g(V \ot W)\ra \ct^g(V)\ot \ct^g(W)$$ of $\ct^g$ is given by $\tau^{-1}(g^{-1})|_{V\ot W}$ and it is a morphism of $A$-modules by condition \ref{A}.
\item
 The fact that $\ct^g$ is a tensor functor is equivalent to the co-cocycle condition for $\tau$ and $g^{-1}$ from \cite{AD}.
%\item
%{\bf The fact that $\ct$ is a tensor functor is equivalent to the condition $(D)$ from \cite{AD}.}

\end{enumerate}

Also from the same paper it follows the following Proposition.

\bp
Assume that $H=A\;^{\tau}\#_{\sg}kG$ is a cocentral extension of $A$ by $kG$. Then $\rep(H)=\rep(A)^G$.
\ep
%}

\bpf
One defines an equivalence of categories $$\mtc{F}:\rep(A)^G \ra \rep(H)$$ by $(V, \mu_V^g)\mapsto V$ where the action of $kG$ on $V$ is given by $\bar{g}.v=\mu_V^g(v)$. The morphisms are sent to themselves.

One has to verify that $V$ becomes an $A\#_{\sg}kG$-module in this way. The equivariant condition \ref{deltau} is equivalent to $\bar{g}(\bar{h}v)=\sg(g,h)(\bar{gh}.v)$. The compatibility between the $A$-action and the $G$-action follows from the fact that $\mu_V^g$ is a morphism of $A$-modules. Indeed, one has that

\begin{eqnarray*}
 \bar{g}.(a.v)  =  u_V^g(a.v) & = & u_V^g((g^{-1}.(g.a)).v)\\ & = & u_V^g((g.a).^gv)\\ & = & (g.a)u_V^g(v)\\ & = & (g.a).(\bar{g}.v)
\end{eqnarray*}

It is easy to verify that any morphism $f:(V, \mu^g_V)\ra (V',\mu_{V'}^g)$ in $\rep(A)^G$ is a morphism of $A\#_{\sg}kG$-modules.

\blue{Is this a tensor functor?}
\md
Define also $$\mtc{E}: \rep(H)\ra \rep(A)^G$$ by $W\mapsto (W, u_W^g)$ where $u^g_W(w):=(\bar{g}.w)$. The fact that $u_W^g$ is a morphism of $A$-modules is equivalent to $\bar{g}(g^{-1}.a)=a\#\bar{g}$.  The fact that $(W, \mu^g_W)$ is an equivariant object is equivalent to the identity $\bar{gh}=\bar{g}\bar{h}\sg(h^{-1}, g^{-1})$.

It is easy to see that $\mtc{F}$ and $\mtc{E}$ are inverses one to the other.

The fact that $\mtc{E}$ is tensor equivalence is equivalent to the identity

$$\tau(g)_i\ot \tau(g)^i=g. \tau^{-1}(g^{-1})_j\ot g. \tau^{-1}(g^{-1})^j$$

Put $h=g^{-1}$ in the equality \ref{D} and take the inverse of $\tau$.

%We say that $G$ measures $A$ via $g.a$ if the following three properties are satisfied:

%Then $G$ acts on $\rep(A)$ via $\ct^g(V)=^gV$ where $^gV=V$ as vector spaces but the $A$-module structure is given by $a.^gv=(g^{-1}.a)v$.

%\textcolor[rgb]{1.00,0.00,0.00}{ Moreover $(\ct_2^{g,h})|_V$ is given by left multiplication with the scalar $\sg(h^{-1}, g^{-1})^{-1}$ and the tensor structure $(\ct^g_2)_{V,W}:\ct^g(V\ot W)\ra \ct^g(V )\ot \ct^g(W)$ is given by left multiplication by $\tau(g^{-1})^{-1}$. }
\epf
\subsubsection {Connections with Clifford theory for cocentral extensions}\lb{cliff}
Suppose that we have a cocentral extension as above. Applying the results from \cite{clifth} it follows that any irreducible $A$-module containing the $B$-module $M$ is induced from an irreducible $Z_A(M)$-module containing $M$. We will give explicitly this description:

Using Lemma 3.6 from \cite{clifth} it follows that the stabilizer $Z_A(M)$ is given by
\beq
Z_A(M)=B\#_{\sg}\; kG_M=\oplus_{g \in G_M}B \#_{\sg}k\bar{g}
\eeq
where $G_M$ is the stabilizer of $M$ in $G$, i.e
\beq
G_M=\{g \in G\;|\; \ct^g(M)\cong M\}
\eeq

Thus an induced module from $Z_A(M)$ is of the form $$A\ot_{Z_A(M)}V\cong \oplus_{t \in G/G_M}k\bar{t}\ot V$$

and the functor $L_Y$ from Theorem \ref{correspondence} coincides with the induction functor of modules from $Z_A(M)$ to $A$.

Moreover the cocycle $\al_{[Y]}$ is constructed the same as in \cite{MW}.
\subsection{Sonia's setings}

Define an action of $G^{\op}$ on $\rep(A)$.
\bne
\item $T^g(V)=V$ with $\bar{g}.v={\mu_V^g}\inv(v)$.
\item ${T_2}^{g, h}:T^g(T^h(V))\ra T^{hg}(V)$ given by $ v \mapsto \sg\inv(h,g)v$.
\item
The tensor structure of $T^g$ is given by $(\ct^g)_2^{V,W}:\ct^g(V \ot W)\ra \ct^g(V)\ot \ct^g(W)$ is given by $\tau^{-1}(g)|_{V\ot W}$.
\item Define an equivalence of categories $\mtc{F}:\rep(A)^G \ra \rep(H)$ by $(V, \mu_V^g)\mapsto V$ where the action of $kG$ on $V$ is given by $\bar{g}.v=(\mu_V^g)\inv(v)$. The morphisms are sent to themselves.

\item Define also $\mtc{E}: \rep(H)\ra \rep(A)^G$ by $W\mapsto (W, u_W^g)$ where $u^g_W(w):=(\overline{g^{-1}}.w)$
\ene

\subsubsection{Identification of the cocycle $\al_{[Y]}$}

Suppose that $Y \in \Irr(B)$ and let $c_f:\ct^f(Y) \ra Y$ be a chosen isomorphism of $B$-modules for any $f \in F_Y$.
Following Equation \ref{alfa} the cocycle $\al_{[Y]}$ is determined by the relation
\bq
c_f(c_g(\sg(g^{-1}, f^{-1})c_{fg}^{-1}(m)))=\al_{[Y]}(f,g)m
\eq
for all $m \in Y$. \bb{Since $c_g$ is a $B$-module isomorphism it follows that $c^g_Y(am)=(g.a)m$ for all $m \in Y$. }Then the above condition can be written as:
\bq
L_{\sg^{-1}(f, g)}\circ c_f\circ c_g \circ c_{fg}^{-1}=\al_{[Y]}^{-1}(f,g)\mtr{Id}_{M}
\eq
Therefore
\bb{\bq\lb{cco}
c_f\circ c_g \circ c_{fg}^{-1}=L_{\sg(f, g)}\al_{[Y]}^{-1}(f,g)\mtr{Id}_{M}
\eqï¿½}
\subsubsection{Identification of the projective $F_Y$-representation $\Hom_B(Y, S)$}
Let $S$ be an $A$-module and $Y$ be a constituent of $S|_B$. Then applying Formula \ref{pr} one has that $\Hom_B(Y, S)$ is a projective $F_Y$ -representation via
\bq\lb{cpr}
(g.u)(m)=\bar{g}.u(c_g^{-1}(m))
\eq
for all $g \in F_Y$ and all $m \in Y$.
To check this is a $F_Y$-projective representation with cocycle $\al_{[Y]}$ one has to use Lemma \ref{sgg}. Indeed,
\beqarn
[\bar{f}(\bar{g}.u)](m) & = & \bar{f}.[(\bar{g}.u)(c_f^{-1}(m))]\\ &=& \bar{f}.[\bar{g}.u(c_g^{-1}(c_f^{-1}(m)))]
\\ &= & \al_{[Y]}(f, g) \bar{f}\bar{g}.[u(c_g^{-1}(c_f^{-1}(m)))]\\ & = & (\sg(f,g)\#\bar{fg}). [u(\sg(g^{-1},f^{-1})c_{fg}^{-1}(m))]
\\ &= & \al_{[Y]}(f, g)\sg(f,g)((fg).\sg(g^{-1},f^{-1}))(\bar{fg}.[u(c_{fg}^{-1}(m))])\\ &=& \al_{[Y]}(f, g)[\bar{fg}.u](m)
\eeqarnï¿½
\subsubsection{Identification of the module $Y\ot V_{\pi}\in \cc^{F_Y}$}
 Note that one has $\cc^{F_Y}=\rep(B\#_{\sg}kF_Y)$ and there is a cocentral extension of Hopf algebras
 \bq
 k \ra B \ra B\#_{\sg}kF_Y\ra k.
 \eq
Following Definition \ref{pr} one can define the following $F_Y$-equivariant object associated to $Y$ and $V_{\pi}$ with equivariant structure:
\bq
u_{Y\ot V_{\pi}}^g:Y\ot V_{\pi}\ra Y\ot V_{\pi}
\eq
given by
\bq\lb{yl}
m \ot v \mapsto c_Y^g(m)\ot \pi(g)(v)
\eq
for all $g \in F_Y$ and $m \in Y$.

Thus the $B\ot kF_Y$-module structure on $Y\ot V_{\pi}$ is given by:
\bq\lb{syl}
[b\#_{\sg} \bar{g}](m \ot v)=bc_Y^g(m)\ot g.v
\eq
\subsubsection{Identification of the simple module $S_{Y,\pi}$}  By definition $S_{Y,\pi}$ is the induced module
\bq
B\#_{\sg}kF\ot_{B\dzs kF_Y}(Y\ot V_{\pi})\cong kF\ot_{kF_Y}(Y\ot V_{\pi})
\eq

Thus the simple $S_{Y, \pi}$ module can be identified to $kF\ot_{kF_Y}(Y\ot V_{\pi})$ where the action of an element $b\#_{\sg} g \in A$ is given by:
\beq\lb{smc}
(b\dzs g)(t\ot_{F_Y}(m \ot v))=s\ot_{kF_Y}((s^{-1}.b)\sg^{-1}(t^{-1},h^{-1})\sg({f'}^{-1},s^{-1})c_{f'}(m)\ot f'.v)
\eeq
where $s \in F/F_Y$ is determined by $gt=sf'$ with $f'=s^{-1}gt \in F_Y$.
\subsubsection{Identification of the projective representation $H_{t,s}^{\pi, \dt}$}
One has
\bq H_{t,s}^{\pi, \dt}:=\Hom_B(U, \ct^t(Y) \ot V_{\pi}\ot \ct^s(Z) \ot V_{\dt}) \eq
with the action of $T:=G_U\cap G_{\ct^t(Y)}\cap G_{\ct^s(Z) }$ given by:
\beqarn
(g.f)(x)& = & (c_Y^{t^{-1}gt}\ot \pi(t^{-1}gt)\ot c^{s^{-1}gs}_Z\ot \dt(s^{-1}gs))\\ & & (\sg(tg^{-1}t^{-1},t^{-1})\sg^{-1}(t^{-1}, g^{-1})\ot \id_V \ot \sg(sg^{-1}s^{-1},s^{-1})\sg^{-1}(s^{-1}, g^{-1}))  \ot \id_{W})\\ & &  \tau^{-1}(g^{-1}) f((c_U^{g})^{-1}(x))
\eeqarn
\subsubsection{Identification of the projective representation $\tau^{\ct^t(Y), \ct^s(Z)}_U$}
Define $m(g;t,s) \in B\ot B$ by $$m(g;t,s):=[\sg(tg^{-1}t^{-1},t^{-1})\sg^{-1}(t^{-1}, g^{-1})\ot \sg(sg^{-1}s^{-1},s^{-1})\sg^{-1}(s^{-1}, g^{-1})] \tau^{-1}(g^{-1}).$$
Then using Formula \ref{pr'} the projective representation of $G_U\cap G_{\ct^t(Y)}\cap G_{\ct^s(Z)}$ is given by
\beqarn
(g.f)(x) & = & (c^{t^{-1}gt}_Y\ot c^{s^{-1}gs}_Z)(m(g;t,s)f((c_U^g)^{-1}(x)))
\eeqarn

\blue{It is projective with cocycle $\al_{[U]}\;^t\al_{[Y]}^{-1}\;^s\al_Z^{-1}$.}

\subsubsection{Identifcation of the projective representation of $\tau_{Y, Y^*}^{\one}$or $H^{\pi,\dt}_{Y,Y^*}$}
One has
\bq H_{t,s}^{\pi, \dt}:=\Hom_B(\one, \ct^t(Y) \ot V_{\pi}\ot \ct^t(Y^*) \ot V_{\dt}) \eq
with the action of $T:=G_U\cap G_{\ct^t(Y)}\cap G_{\ct^s(Z) }$ given as following:
If $f$ is given by $\one \mapsto y_i\ot y^{*}_i \ot v_i \ot w_i$ then $g.f$ is given by
\beqarn
1\mapsto (c^{t^{-1}gt}_Y\ot c^{t^{-1}gy}_{Y^*})[m(g,t,t)(y_i\ot y^{*}_i)] \ot \pi(t^{-1}gt)v_i \ot \dt(s^{-1}gs) w_i
\eeqarn
\bb{This has to be real representation of $G$.}
\subsubsection{Identifcation of the isomorphism $d_{Y, t}:\ct^t(Y)^*\ra \ct^t(Y^*)$}
It is easy to see that this isomorphism can be taken to be identity.
\subsubsection{Comments}
\blue{ Compare this definition if one replaces $\tau^{Y,Z}_U$ with $\ct^t(Y), \ct^s(Z)$. This generalizes the results obtained by Goof.}
%%%%%%%%%%%%%%%%%%%%%%%%%%%%%%%%%%%%%%%%%%%%%%%%%%%
\subsection{Cocentral abelian extensions of Hopf algebras}%; the case of Goof}
Let $A=k^F$ in the previous subsection. Then denote $$\sg(g,h)=\sum_{x\in F}\sg_x(g,h)p_x$$ and $$\tau(g)=\sum_{x,y\in F}\tau_g(x,y)p_x\ot p_y$$

The above identities become as follows:

\bne
\item $\sg_x(1,g)=1=\sg_x(g,1)$
\item $\tau_g(1,x)=1=\tau_g(x,1)$
\item $\sg_x(h,l)\sg_x(hl,m)=\sg_{h^{-1}.x}(l,m)\sg_x(h,lm)$
\item $h.(l.x)=(hl).x$
\item $\tau_g(mn,r)\tau_g(m,n)=\tau_g(m,nr)\tau_g(n,r)$
\item This identity is automatically satisfied.
\item This identity becomes the pentagon equation:
$$
\frac{\sg_{mn}(g,h)}{\sg_{m}(g,h)\sg_{n}(g,h)}=\frac{\tau_g(m,n)\tau_h(g^{-1}.m, g^{-1}.n)}{\tau_{gh}(m,n)}
$$
\ene
%\blue{this should come from the previous section}\blue{ Explain how to get the formula of Goof}

\subsubsection{Identification of the cocycle $\al_{[Y]}$}
In this situation $Y=k_g \cong kp_g$ is a one dimensional $B$-module.
For any $x\in F_g$ choose a scalar $c^g_x \in k^*$. Then following Equation \ref{cco} the cocycle $\al_g$ on $F_g$ is determined by:
\beq
\al_g(x, y)=(c^g_{xy})^{-1}c_x^gc_y^g\sg_g(y^{-1}, x^{-1})
\eeq
Thus $\al_g$ is cohomologus to the cocylce $(x,y)\mapsto \sg_g(y^{-1}, x^{-1})$.
\br
Note that all this scalars $c^g_x$ can be taken to be $1$ and in this case one gets equality between the two cocycles.
\er
\subsubsection{Identification of the projective representation $\Hom_B(Y, S)$}
Using Equation \ref{cpr} it follows that the projective representation $\Hom_{k^G}(k_y, S)$ becomes
\bq
[x.u](p_g)=(c^g_x)^{-1}\bar{x}.u(p_g)
\eq
\subsubsection{Identification of the module $Y\ot V_{\pi}\in \cc^{F_Y}$}
Following Equation \ref{syl} it follows that $1_g \ot V_{\pi}$ is a $k^G\dzs kF_g$-module via:
\beq
(p_g\dzs \bar{f})(1_m\ot v)=\delta_{g, m}c^m_f(1_m\ot \pi(f)m)
\eeq
\subsubsection{Identification of the simple module $S_{g,\pi}$}
The simple $A$-module $S_{g,\pi}=S_{Y, \pi}$ is identified with
\bq
kF\ot_{k{F_g}}( k_g \ot V_{\pi})
\eq
as in Equation \ref{smc}. Moreover this can be identified with $kF\ot_{kF_g}V_{\pi}$ with the following structure following from Equation \ref{smc}:
\bq
(p_g \dzs f)(t\ot_{kF_m}(1_m \ot v))=\delta_{g, t\rhd m}\sg_m(f,t)\sg^{-1}(s, f')c^m_{f'}(s\ot_{kF_m}(1_m\ot \pi(f')v)
\eq
\subsubsection{}
In this situation $m(g;t,s)$ is given by
\bq
\sum_{x,y\in G}\sg_x(tg^{-1}t^{-1},t^{-1})\sg_x^{-1}(t^{-1}, g^{-1})p_x\ot \sg_y(sg^{-1}s^{-1},s^{-1})\sg_y^{-1}(s^{-1}, g^{-1}) \tau^{-1}_{g^{-1}}(x,y)p_y
\eq
\subsubsection{Identification of the projective representation $\tau^{Y, Z}_U$}
Suppose that $Y=k_g$, $Z=k_h$ and $U=k_{gh}$.
The space $\Hom_{B}(k_{gh}, k_g\ot k_h)$ is one dimensional.  It is a projective representation of $F_g\cap F_h\cap F_{gh}=F_g\cap F_h$ via the following structure:
\beq
x.f=(c_x^{gh})^{-1}c^g_xc^h_x\tau_x^{-1}(g,h)f
\eeq

%\blue{Is this equivalent to a real representation? Since it projective it should be!}
In particular $\tau_y:=\tau^{y,y^*}_{\one}$ is given by multiplication with $\tau_y(g,g^{-1})$.
\subsubsection{Identification of the projective representation $\tau^{Y,Z;t,s}_U$}
The space $\Hom_B(\ct^t(y)\ct^s(z), \ct^t(y)\ct^s(z))$ is a one dimensional projective representation by
\bq
g.f=(c_g^{\ct^t(y)\ct^s(z)})^{-1}c^g_{\ct^t(x)}c^g_{\ct^s(x)}\mu_g(\ct^t(x), \ct^s(y))\gm^{-1}(g,t)\gm(g,s)\gm^{-1}(t, t^{-1}gt)\gm^{-1}(s, s^{-1}gs)
\eq
\subsubsection{Fusion rules for  abelian cocentral extensions}
One has the following formula for fusion rules:
\bq
S_{g,\pi}\ot S_{h, \delta}\cong\bigoplus_{t \in D} S_{g\;^th, \;(\pi\dw_{G_x\cap \;^tG_y}\ot ^t\dt \dw_{G_x\cap \;^tG_y})\uw^{G_{x\ct^t(y)}}}
\eq
where $D$ is a set of double coset representatives of $F_x\backslash F/ F_y$.
\bpf
Looking at restriction to $k^X$ one sees that the possible simple objects of $\cc$ are $x\ct^t(y)$ with $t \in D$. Then it can be seen that for a fixed $t \in D$ the set
$\Hom_{\cc}(x\;^ty, x\;ty\ot V_{\pi}\ot V_{\dt}$ as $T$-projective representation coincides to
$$(\pi\dw^{G_x}_{G_x\cap \;^tG_y}\ot \;^t\dt \dw^{G_{\;^ty}}_{G_x\cap \;^tG_y})\uw^{G_{x\ct^t(y)}}_{G_x\cap \;^tG_y}.$$
\epf

\subsubsection{Fusion rules for  quasi-abelian cocentral extensions}
One has the following formula for fusion rules:
\bq
S_{g,\pi}\ot S_{h, \delta}\cong\bigoplus_{t \in D} S_{g\;^th, \;(\pi\dw_{G_x\cap \;^tG_y}\ot ^t\dt \dw_{G_x\cap \;^tG_y})\uw^{G_{x\ct^t(y)}}}
\eq
where $D$ is a set of double coset representatives of $F_x\backslash F/ F_y$.
\bpf
Looking at restriction to $k^X$ one sees that the possible simple objects of $\cc$ are $x\ct^t(y)$ with $t \in D$. Then it can be seen that for a fixed $t \in D$ the set
$\Hom_{\cc}(x\;^ty, x\;ty\ot V_{\pi}\ot V_{\dt}$ as $T$-projective representation coincides to
$$(\pi\dw^{G_x}_{G_x\cap \;^tG_y}\ot \;^t\dt \dw^{G_{\;^ty}}_{G_x\cap \;^tG_y})\uw^{G_{x\ct^t(y)}}_{G_x\cap \;^tG_y}.$$

There is an element $g_i(t)$ such that $g\;^th=\;^s(g_{i(t)}).$ Therefore
$$S_{g\;^th, (\pi\dw^{G_x}_{G_x\cap \;^tG_y}\ot \;^t\dt \dw^{G_{\;^ty}}_{G_x\cap \;^tG_y})\uw^{G_{x\ct^t(y)}}_{G_x\cap \;^tG_y}}=S_{g_i(t), d_s^{-1}(\pi\dw^{G_x}_{G_x\cap \;^tG_y}\ot \;^t\dt \dw^{G_{\;^ty}}_{G_x\cap \;^tG_y})\uw^{G_{x\ct^t(y)}}_{G_x\cap \;^tG_y}}$$

It follows that $$\lam(a\ct^t(b), h)=d_s^{-1}(h)\pi(h)\;^t\dt(h).$$

Thus
$$\lam(a_ia_j, h)=d_s^{-1}(h)d^t(h)\lam(a_i, h)\lam(a_j, h).$$
\epf

\subsubsection{The character $\tau_{y,z}^{yz}(g)$} It can easily be computed that $$\tau_{y,z}^{yz}(g)=\mu_g(y,z).$$ for all $g \in G_y \cap G_z$
\subsubsection{Fusion subcategories}

Following Theorem \ref{cls-by-fusiondata} it follows that fusion subcategories of $\cc(\omega, \gm,\mu, c)^G$ are parameterized by
\bne
\item Fusion subcategories $\cs$ of $\bvec^{\omg}_X$; thus $\cs=\bvec^{\omg}_Y$ with $Y$ a subgroup of $X$.
\item Subgroups of $H \leq G$ acting trivially on $\Irr(\cs)=Y$.
\item A twisted bicharacter $\lam:Y \times H \ra k^*$ i.e, verifying the following properties:
\bne
\item $\lam(hh', y)=\lam(h,y)\lam(h',y)\gm\inv_{y, y'}(h)$
\item $\lam(yy',h)=\lam(y,h)\lam(y', h)\mu_h(y,y')$
\item The character of $k_{\gm_{yz} |_H}[H]$ given by $\lam(-,\;y) \dw^{G_y}_H\lam(-,\;z)\dw^{G_z}_H\tau_{yz}^{y, z}\dw^{G_y \cap G_z}_H$ is a $G_{yz}$-stable character of $H$ with respect to the 2 cocyle $\gm_{yz}$
\item $\lam(^gy,\;h)=d_y(g,h)\lam(y,ghg^{-1})$
\ene
\ene
Denote by $\cd(Y, H, \lam)$ the corresponding fusion subcategory.
%\bb{Then make $\omega=1$ to get the abelian cocentral extension case of Hopf algebras.}
\subsubsection{}It seems that the cocentral abelian quasi-Hopf algebra should be
\bq
k^X\;^{\mu^{-1}}\#_{\gm^{-1}\tau}kG.
\eq
where $\gm^{-1}\tau(g,h)=\gm(h^{-1}, g^{-1})$.

Maybe this quasi-Hopf algebra fits to the twisted Drinfeld Double.

\bl
One has that $\cc^{G_y}=\rep(k^X\#_{\gm}kG_y)$.
\el

\bl With the above notations one has that:
$S_{a, \pi(g)}=S_{\;^ta, \;\frac{\gm_{t, t\inv gt}(a)}{\gm_{g,t}(a)}\;^t\pi(g)}$.
\el
\bpf
One has to look at there restriction of $S_{a,\pi}$ at $k^X\#_{\al_{\;^ta}}kG_{\;^ta}$. Following Equation \ref{eqvr} for $g \in G_{\;^ta}$ it follows that
\bq
\ct^t(a) \ot v \xra{\mu_{S_{a,\pi}}^g} \ct^t(a) \ot \frac{\gm_{t, t\inv gt}(a)}{\gm^{\inv}_{g,t}(a)}\pi(t^{-1}gt)(v)
\eq
\epf
Thus $d_a(t, g)=\frac{\gm_{g,t}(a)}{\gm_{t, t\inv gt}(a)}$.
%\bb{Find exactly the scalar in the general proof.}
\subsubsection{M\"{u}ger centralizer of these fusion subcategories}
\subsubsection{Lagrangian subcategories} %examples from 3 qabelian from naidu
%\textcolor[rgb]{1.00,0.00,0.00}{
\section{\blue{The twisted Drinfeld double $D^{\omega}(G)$.}}
%\textcolor[rgb]{1.00,0.00,0.00}{
%%%%%%%%%%%%%%%%%
\subsection{The quasi-Hopf algebra structure of $D^{\omega}(G)$}
Consider the quasi Hopf algebra $D^{\og}(G)$ fas defined in \cite{dpr}. One has that the associator is given by the formula
 $$\phi=\sum_{a,b,c}\omega^{-1}(a,b,c)p_a\ot _b \ot p_c.$$
 %%%%%%%%%%%%%%%%%%%%%%%
 \subsubsection{The action of $G$ on ${\mtr{Vec}^{\omega\inv}}_G$}\lb{actio}
Using Lemma 6.3 from \cite{Na} it follows that $$\rep(D^{\omega}(G))\cong \cc(\omg^{-1}, \gm\inv,\mu\inv, 1)^G.$$
%}
The action of $G$ on $\bvec^{\omg^{-1}}_G$ is given by the adjoint action $$g.k_x=k_{gxg^{-1}}$$ and the following data:
$$(\ct_2^{g,h})|_x: \ct^g(\ct^h(k_x))\ra \ct^{gh}(k_x),\;\; v \mapsto \gamma\inv_x(g,h)v$$
and $$(\ct^g_2)_{x,y}:\ct^g(k_x\ot k_y)\ra \ct^g(k_x)\ot \ct^g(k_y),\;\;v\ot w\mapsto \mu_g
\inv (x,y)(v \ot w).$$
%\textcolor[rgb]{0.10, .00, 0.90}{
\subsubsection{Chosen isomorphisms} On can choose the isomorphisms $c_Y^g$ al to be identities. Denote by $C=\mtr{1}$ the above collection of chosen isomorphisms.
\subsubsection{The cocycle $\al_y$}From formula \ref{alfa} it follows that in this case the cocylce $\al_y$ on $G_y$ is given by $$\al_y(g,h)=\gm_{g,h}(y)=\eta^{-1}_y(g,h)=\beta\inv_y(g,h)$$
since $g,h \in C_G(y)$.
%%%%%%%%%%%%%%%%%%%%%%%%
\subsubsection{} We identify the notions from \cite{Na} with those from \cite{nnw}. With their notations one has $$\gamma_{g,h}(x)=\eta_{x}(g,h)$$ and $$\mu_g(x,y)=\nu_g^{-1}(x, y).$$ On the other hand,  remark that for any $g,h \in C_G(x)$ one has $\gamma_{g,h}(x)=\beta_x(g,h)$%; maybe $\beta$ appears on the centralizers; Alos on the centralizers one has
and $\mu_g(x,y)=\gamma_g(x,y)^{-1}$.

\subsubsection{On the projective representation $\tau_{U}^{Y, Z}$} With the above action and using Formula \ref{pr'} one finds the projective representation $\tau_{yz}^{y,z}$ given on the space $\Hom_k(yz, yz)$ as a multiplication with a scalar, namely  $g.f=\mu_{g}\inv(y,z)f$. Thus
\bq
\tau_{y,z}^{yz}(g)=\mu_g\inv(y,z)=\nu_g(y,z)=\beta_g(y,z)
\eq
since $y, z \in C_{G}(g)$
Note that in this case $\tau^{y,z}_{yz}(1)=1$.
\subsubsection{}
Following Equation \ref{dyt} it follows that
\bq
d_a(t, h)=\frac{\gm_{tht^{-1},t}(a)}{\gm_{t,h}(a)}=\frac{\eta_a(tht^{-1},t)}{\eta_a(t,h)}
\eq
for a fixed $t \in G/G_a$ all $h \in H$.
%%%%%%%%%%%%%%%%%%%%
\subsection{A fusion datum for $D^{\omg}(G)$}
Suppose that $G$ acts on $Vec^{\omega\inv}_G$ as in subsection \ref{actio}. Then as in the definition from Section \ref{intro} a fusion datum will consists of a triple $(\cs, H, \lam)$. \\The fusion subcategory $\cs$ of $\bvec^{\omg}_G$ should be $\cs=\mtr{Vec}^{\omega}_K$ where $K$ is a subgroup of $G$. Then $K$ is a normal subgroup of $G$ since $\cs$ is stable under the action of $G$.  \\ Since $H$ acts trivially on $\cs$ it follows that $H$ and $K$ commute element wise, i.e $[H,\;K]=1$. On the other hand since $\gr(\cs)=k[K]$ it follows that $\lam$ can be viewed as a function $\lam_{\cd}: K\times H \ra U(1)$ satisfying the following properties:
$$\lam_{\cd}(y,hh')=\lam_{\cd}(y,h)\lam_{\cd}(y, h')\beta\inv_y(h, h')$$ and $$\lam_{\cd}(yz, h)=\lam_{\cd}(y,h)\lam_{\cd}(z,h)\beta_h(y,z).$$
The $G$-invariance of $\lam_{\cd}$ with respect to $d$ can be written now as:
\bq\lb{ginvl}
\lam_{\cd}(txt^{-1}, h)=d_x(t,h)\lam_{\cd}(x, t\inv ht)
\eq
\subsection{Comparison with the results of \cite{nnw}} Let $\cd$ be a fusion subcategory of $\Rep(D^{\og}(G))$. We will show that in this case our results recover the results of \cite{nnw}. More precisely we will show that the fusion datum associated to $\cd$ coincides with the fusion datum from \cite{nnw}. One has that $S_{\cd}=\mtr{Vec}^{\og}_{K_{\cd}}$ and $H_{\cd}$ coincide since in both cases it is found by $\Rep(G/H_{\cd})=\cd \cap \rep(G)$. It remains to show that $\lam_{\cd}=B_{\cd}$ from formula (22) of \cite{nnw}.
%\textcolor[rgb]{.00, .00, 1.00}{
\subsubsection{The equality bewteen simple modules}
Using Lemma \ref{chf} it follows that $S_{a,\ch}=S_{x^{-1}ax, \delta}$ where
$$
\delta(xhx^{-1})=\frac{\beta_a(x^{-1}hx,x^{-1})}{\beta_a(x^{-1},h)}\ch(h)
$$

\subsubsection{Coincidence of bicharacters}
Since $$S_{a, \ch}=S_{x\inv a x, \;^xd_a(x\inv, -)\;^{x\inv}\pi}$$
one has $\lam_{\cd}(x^{-1}ax, h)=d_a(x\inv, h)\frac{\ch(xhx^{-1})}{\ch(1)}$. On the other hand note that by formula (22) from \cite{nnw} one has that
$$
B_{\cd}(x^{-1}ax,h)=\frac{\beta_a(x,h)\beta_a(xh, x^{-1})} {\beta_a(x^{-1}, x)}\frac{\ch(xhx^{-1})}{\ch(1)}
$$

Since $h \in \mtr{core}_G(G_a)$ one has that $h \in G_{x^{-1}ax}$ and therefore
$$
\delta(h)=\delta(x^{-1}(xhx^{-1})x)=\beta_a(h,x^{-1})\beta_a^{-1}(x^{-1},xhx^{-1})\ch(xhx^{-1})
$$
Therefore if $S_{a, \ch} \in \cd$ then
\bn{eqnarray*}
\lam_{\cd}(x^{-1}ax, h) & = &\frac{\delta(h)}{\delta(1)}=\frac{\delta(x^{-1}(xhx^{-1})x)}{\delta(1)}
\\ &= &\frac{ \beta_a(h, x^{-1})}{\beta_a(x^{-1},xhx^{-1})}\frac{\ch(xhx^{-1})}{\ch(1)}
\end{eqnarray*}
Applying Lemma \ref{idnt} one obtains that
$$
\lam_{\cd}(x^{-1}ax, h) =\frac{\beta_a(x,h)\beta_a(xh, x^{-1})} {\beta_a(x, x^{-1})}\frac{\ch(xhx^{-1})}{\ch(1)}
$$
Note that by formula (22) from \cite{nnw} one has that
$$
B_{\cd}(x^{-1}ax,h)=\frac{\beta_a(x,h)\beta_a(xh, x^{-1})} {\beta_a(x^{-1}, x)}\frac{\ch(xhx^{-1})}{\ch(1)}
$$
Since $\beta_a(x,x^{-1})=\beta_a(x^{-1},x)$ it follows that
$$
\lam_{\cd}(x^{-1}ax, h) =B_{\cd}(x^{-1}ax,h)
$$
\subsection{Verifying the properties of the bicharacters}
\subsubsection{Properties of the character $B_{\cd}$}
The $G$-invariance of $B_{\cd}$ can be written as:
\bq\lb{ginv}
B_{\cd}(yky\inv, h)=c_{\beta_h}(y,k)B_{\cd}(k, y\inv hy)
\eq
On the other hand Lemma 5.6 from \cite{nnw} implies that
\bq
c_{\beta_h}(y,k)=c_{\beta_k}(y\inv, h)
\eq
since $hk=kh$.
Thus the invariance can be written as
\bq\lb{ginv2}
B_{\cd}(yky\inv, h)=c_{\beta_k}(y\inv,h)B_{\cd}(k, y\inv hy)
\eq
\subsubsection{Veryfing the first property}
One has
\bq
\al_y(g, h)=\gamma^{-1}_{g, h}(y)=\beta^{-1}_y(g,h)=\eta^{-1}_y(g,h)
\eq
since $g, h\in C_G(y)$. In this way the first property \ref{yz} of $\lam_{\cd}$ coincides to the property $i$ from Definition 5.4 from \cite{nnw}.
\subsubsection{Verifying the second property}
On the other hand it can be checked by the definition that $\tau^{yz}_{y,z}(g)=\mu_g(y,z)$ for all $g \in C_G(y)\cap C_G(z)$.
Thus $\tau^{yz}_{y,z}(g)=\beta_g(y,z)$ on centralizers and Property \ref{hh'} coincides to the second property of Definition 5.4 from \cite{nnw}.
\subsection{Veryfying the $G$-invariance}
To verify the $G$-invariance one has to verifyy that
\bq
c_{\beta_k}(y\inv, h)=d_k(y,h)
\eq

From the proof of $G$-invariance of $\lam_{\cd}$ it is enough to verify the above identity for $k=a$ chosen in $\Gm$ and $y=t$ chosen in $\Gm_a$. In this case the above identity reduces to $$c_{\beta_a}(t\inv,h)=d_a(t,h).$$
By Lemma \ref{idnt} one has:
\bq
c_{\beta_a}(t\inv,h)=\frac{\beta_a(h, t)}{\beta_a(t, t\inv ht)}
\eq
On the other hand one has
\bq
d_a(t,h)=\frac{\eta_a(tht\inv, t)}{\eta_a(t,h)}
\eq
Thus one has to verify the following identity:
\bb{\bq\lb{tvf}
\frac{\beta_a(t\inv, h)\beta_a(t\inv h, t)}{\beta_a(t\inv, t)}=\frac{\eta_a(tht\inv, t)}{\eta_a(t,h)}.
\eq}
%he third condition is automatically satisfied and the fourth coincides to the $G$-invariance from \cite{nnw} of $\lam_{\cd}$.
%\input{applications}

%\bb{Move it to the appendix}% \el
\bibliographystyle{amsalpha}

\end{document}